\documentclass[10pt]{amsart}
\usepackage[top=0.5in, bottom=0.5in, left=0.5in, right=0.5in]{geometry}

\usepackage{framed,multirow}

\usepackage{amsmath,amssymb,amsthm,amsrefs,a4wide}
\usepackage{latexsym}

\usepackage{mathtools}
\usepackage[foot]{amsaddr}
\usepackage{url}
\usepackage{xcolor}
\usepackage{graphicx,color}
\usepackage{diagbox}
\usepackage{multirow}
\usepackage{subcaption}
\usepackage{epstopdf}
\usepackage{comment}
\usepackage{braket}
\usepackage{caption}
\usepackage{adjustbox}
\usepackage{xcolor}
\usepackage{bbm}
\usepackage{enumerate}
\usepackage{bm}
\usepackage{algorithm,algpseudocode}

\algrenewcommand\textproc{}

\newcommand\pfrac[2]{\frac{\partial #1}{\partial #2}}

\newcommand\figref[1]{Fig.\ \ref{#1}}
\newcommand\secref[1]{Section \ref{#1}}

\newcommand{\bvec}[1]{\bm{#1}}

\newcommand{\vx}{\bvec{x}}

\newcommand{\bbN}{\mathbb{N}}

\newcommand{\bbR}{\mathbb{R}}

\newcommand{\ccG}{\mathcal{G}}

\global\long\def\Rd{\mathbb{R}^{d}}

\numberwithin{equation}{section}
\numberwithin{figure}{section}
\newtheorem{rem}{\protect\remarkname}
\newtheorem{lemma}{Lemma}
\newtheorem{theorem}{Theorem}
\newtheorem{assumption}{Assumption}
\newtheorem{defn}{Definition}

\providecommand{\remarkname}{Remark}

\title[Combining monte carlo and tensor-network methods for partial differential equations]{Combining monte carlo and tensor-network methods for partial differential equations by sketching}%

\author{Yian Chen}
\address{Department of Statistics, University of Chicago}

\author{Yuehaw Khoo}
\address{Department of Statistics, University of Chicago}

\author{Ziang Yu}
\address{Committee on Computational and Applied Mathematics, University of Chicago}

\email{yianc@uchicago.edu, ykhoo@uchicago.edu, ziangyu@uchicago.edu}

\begin{document}
\maketitle
\begin{abstract}
    In this paper, we propose a general framework for solving high-dimensional partial differential equations with tensor networks. Our approach uses Monte Carlo simulations to update the solution and re-estimates the new solution from samples as a tensor network using a recently proposed tensor train-sketching technique. We showcase the versatility and flexibility of our approach by applying it to two specific scenarios: quantum imaginary-time evolution via auxiliary-field quantum Monte Carlo and simulating the Fokker-Planck equation through Langevin dynamics. We also provide convergence guarantees and numerical experiments to demonstrate the efficacy of the proposed method.
\end{abstract}

\section{Introduction}

High-dimensional partial differential equations (PDEs) describe a wide range of phenomena in various fields, including physics, engineering, biology, and finance. However, the traditional finite difference and finite element methods scale exponentially with the number of dimensions. To circumvent the curse of dimensionality, researchers propose to pose various low-complexity ansatz on the solution to control the growth of parameters. For example, \cites{han2018solving,yu2018deep,zhai2022deep} propose to parametrize the unknown PDE solution with deep neural networks and optimize their variational problems instead. \cites{ambartsumyan2020hierarchical,chen2021scalable,chen2023scalable,ho2016hierarchical} approximates the differential operators with data-sparse low-rank and hierarchical matrices. \cite{cao2018stochastic} considers a low-rank matrix approximation to the solution of time-evolving PDE.

The matrix product state (MPS), also known as the tensor train (TT), has emerged as a popular ansatz for representing solutions to many-body Schr\"odinger equations \cites{white1992density,chan2011density}. Recently, it has also been applied to study statistical mechanics systems where one needs to characterize the evolution of a many-particle system via Fokker-Planck type PDEs \cites{dolgov2011fast,chertkov2021solution,tt_committor}. Despite the inherent high-dimensionality of these PDEs, the MPS/TT representation mitigates the curse-of-dimensionality challenge by representing a $d$-dimensional solution through the contraction of $d$ tensor components. Consequently, it achieves a storage complexity of $O(d)$.

To fully harness the potential of MPS/TT in solving high-dimensional PDEs, it is crucial to efficiently perform the following operations:
\begin{enumerate}
    \item Fast applications of the time-evolution operator $P$ to an MPS/TT represented solution $\phi$.
    \item Compression of the MPS/TT rank after applying $P$ to $\phi$, as the rank of $P\phi$ can be larger than that of $\phi$.
\end{enumerate}
While these operations can be executed with high numerical precision and $O(d)$ time complexity when the PDE problem exhibits a specialized structure (particularly in 1D-like interacting many-body systems), general problems may necessitate exponential running time in dimension $d$ to perform these tasks.

On the other hand, Monte Carlo methods employ a representation of $\phi$ as a collection of random walkers. The application of a time-evolution operator $P$ to such a particle representation of $\phi$ can be achieved inexpensively through short-time Monte Carlo simulations. As time progresses, the variance of the particles, or random walkers, may increase, accompanied by a growth in the number of walkers. To manage both the variance and computational cost of the random walkers, it is common to use importance sampling strategies or sparsification of random walkers \cites{greene2019beyond,lim2017fast}. Furthermore, quantum Monte Carlo methods \cite{Ceperley1995} suffer from an exponentially large variance, resulting from having negative or complex sample weights (a problem commonly known as ``sign problem''). Most quantum Monte Carlo methods require adding certain constraints to restrict the space of the random walks, reducing the variance at the expense of introducing a bias. This gives rise to another challenge: how to unbias quantum Monte Carlo to systematically improve the results  \cite{qin2016coupling}?

Our contribution is to combine the best of both worlds. Specifically, we adopt the MPS/TT representation as an ansatz to represent the solution, and we conduct its time evolution by incorporating short-time Monte Carlo simulations. This integration of methodologies allows us to capitalize on the advantages offered by both approaches, leading to improved performance and broader applicability in solving high-dimensional PDEs. The improvement is two-folds,
\begin{enumerate}
\item From the viewpoint of improving tensor network methods, we simplify the application of a semigroup $P$ to an MPS/TT $P \phi$ via Monte Carlo simulations. While the application of $P$ using Monte Carlo is efficient, one needs a fast and accurate method to estimate an underlying MPS/TT from the random walkers. To address this requirement, we employ a recently developed parallel MPS/TT sketching technique, proposed by one of the authors, which enables estimation of the MPS/TT from the random walkers without the need for any optimization procedures.
\item From the viewpoint of improving Monte Carlo methods, we propose a novel approach to perform walker population control, by reducing the sum of random walkers into a low-rank MPS/TT. In the context of quantum Monte Carlo, we observe a reduced variance in simulations. Furthermore, when it comes to probabilistic modeling, the MPS/TT structure possesses the capability to function as a generative model \cite{ruthotto2021introduction}, enabling the generation of fresh samples and conditional samples from the solution. 
\end{enumerate}
We demonstrate the success of our algorithm in both statistical and quantum mechanical scenarios for determining the transient solution of parabolic type PDE. In particular, we use our method to perform Langevin dynamics and auxiliary-field quantum Monte Carlo for systems that do not exhibit 1D orderings of the variables, for example, 2D lattice systems. We further provide convergence analysis in the case of solving the Fokker-Planck equation, which demonstrates variance error that does not suffer from the curse of dimensionality.

The rest of the paper is organized as follows. First, we introduce some preliminaries for tensor networks, in \secref{sec:background}. We discuss the proposed framework that combines Monte Carlo and MPS/TT in \secref{sec:framework}. In \secref{sec:applications}, we demonstrate the applications of our proposed framework on two specific evolution systems: the ground-state energy problem for quantum many-body system and the density evolution problem by solving the Fokker-Planck equation. In Section~\ref{section:convergence}, we prove the convergence of the proposed method for solving the Fokker-Planck equation. The corresponding numerical experiments for the two applications are provided in \secref{sec:experiments}. We conclude the paper with discussions in \secref{sec:discussions}.

\section{Background and Preliminaries}
\label{sec:background}

In order to combine particle-based simulations and tensor-network-based approaches, we need to describe a few basic tools regarding tensor-network. 


\subsection{Tensor networks and notations}

Our primary objective in this paper is to obtain an MPS/TT representation of the solution of the initial value problem \eqref{eq:ivp} at any time $t\geq 0$. Since the technique works for $u(x,t)$ at any given $t\geq 0$, we omit $t$ from the expression and use $u(x)$ to denote an arbitrary $d$-dimensional function:
\begin{align}
    u:\ X_1\times X_2 \times \cdots \times X_d \rightarrow \bbR,\ \text{where } x_1\in X_1,\ x_2\in X_2,\cdots, x_d\in X_d.
\end{align}
And the state spaces $X_1,\cdots,X_d\subseteq \bbR$. 

\begin{defn}
    \label{def:mps}
    We say function $u$ admits an MPS/TT representation with ranks or bond dimensions $(r_1,\dots,r_{d-1})$, if one can write
    \begin{align}
        u(x_1,x_2,\dots,x_d)\  =& \ \  \sum_{\alpha_1 = 1}^{r_1} \sum_{\alpha_2 = 1}^{r_2} \cdots \sum_{\alpha_{d-1} = 1}^{r_{d-1}} \ccG_1(x_1,\alpha_1)\ccG_2(\alpha_1, x_2, \alpha_2)\dots\ccG_d(\alpha_{d-1}, x_d)
        \label{eq:TT}
    \end{align}
    for all $(x_1,x_2,\dots,x_d)\in X_1\times \cdots \times X_d$. We call the $3$-tensor $\ccG_k$ the $k$-th tensor core for the MPS/TT.
\end{defn}

\begin{figure}[htb]
    \centering
    \begin{subfigure}{0.32\textwidth}
        \centering
        \includegraphics[width=\textwidth]{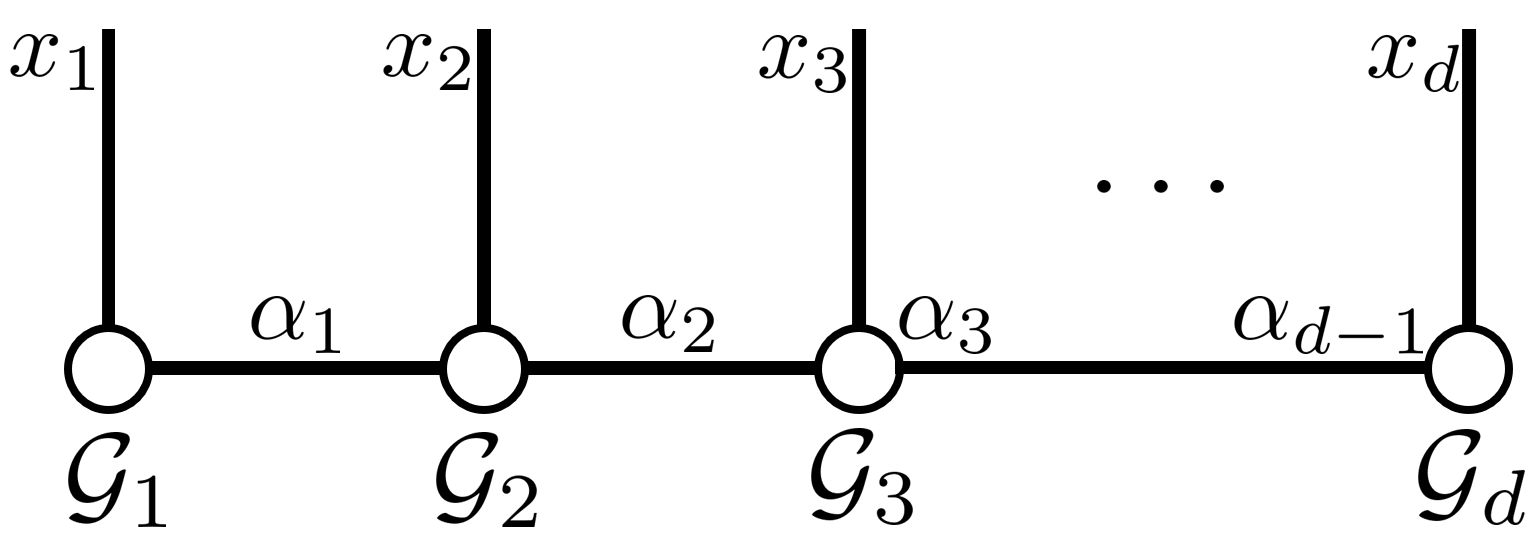}
        \caption{MPS/TT}
        \label{fig:MPS/TT}
    \end{subfigure}\quad
    \begin{subfigure}{0.32\textwidth}
        \centering
        \includegraphics[width=\textwidth]{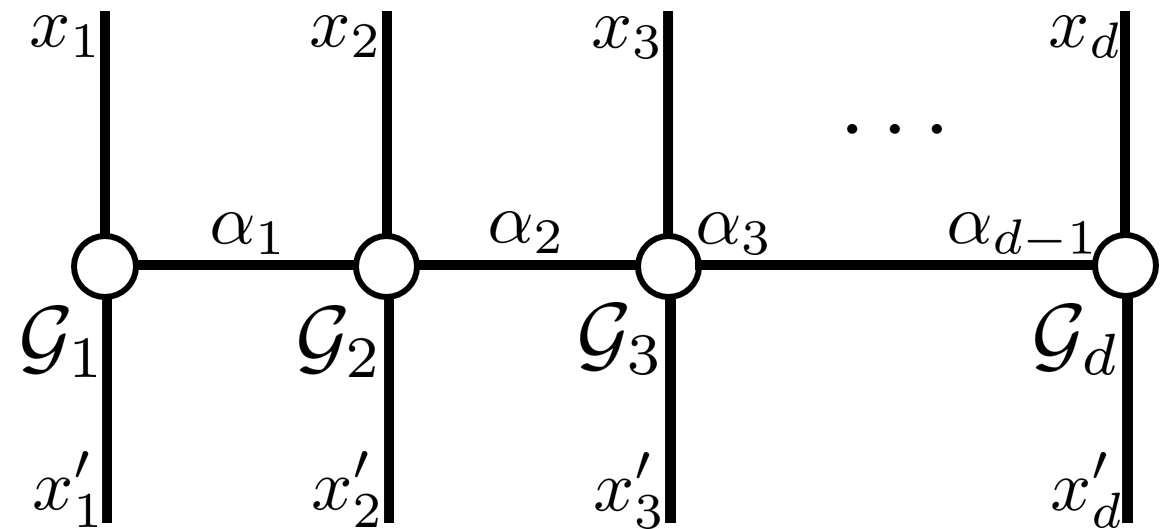}
        \caption{MPO}
        \label{fig:MPOs}
    \end{subfigure}
    \caption{Tensor diagram for a $d$-dimensional MPS/TT and MPO. (a): An MPS/TT representing $u$ in Definition~\ref{def:mps}. The exposed legs indicates the MPS/TT takes $x_1,\ldots,x_d$ as inputs, and the connected legs indicate the summation over $\alpha_1,\ldots,\alpha_{d-1}$. (b) An MPO representing $O$. Each tensor core has two exposed legs pointing upwards and downwards, respectively, indicating two free dimensions.}
\end{figure}

We present the tensor diagram depicting the MPS/TT representation in \figref{fig:MPS/TT}. In this diagram, each tensor core is represented by a node, and it possesses one exposed leg that represents the degree of freedom associated with the corresponding dimension. 
\begin{defn}
    \label{def:mpo}
    $O$ admits a matrix product operator (MPO) representation with ranks or bond dimensions $(r_1,\dots,r_{d-1})$, if one can write
    \begin{align}
        O(x_1,\dots,x_d;x_1',\dots,x_d') \ = & \ \  \sum_{\alpha_1 = 1}^{r_1} \sum_{\alpha_2 = 1}^{r_2} \cdots \sum_{\alpha_{d-1} = 1}^{r_{d-1}} \ccG_1(x_1,x_1',\alpha_1) \ccG_2(\alpha_1,x_2,x_2',\alpha_2) \dots\ccG_d(\alpha_{d-1}, x_d,x_d') 
        \label{eq:MPO}
    \end{align}
    for all $(x_1,x_2,\dots,x_d)\in X_1\times \cdots \times X_d$, $(x_1',x_2',\dots,x_d')\in X_1'\times \cdots \times X_d'$. We call the $4$-tensor $\ccG_k$ the $k$-th tensor core for the MPO.
\end{defn}

The tensor network diagram corresponding to the MPO representation is depicted in \figref{fig:MPOs}. For more comprehensive discussions on tensor networks and tensor diagrams, we refer interested readers to \cite{tt_committor}.

Finally, we introduce some indexing conventions. For two integers $m,n\in\bbN$ where $n> m$, we use the \texttt{MATLAB} notation $m:n$ to denote the set $\{m,m+1,\cdots,n\}$. When working with high-dimensional functions, it is often convenient to group the variables into two subsets and think of the resulting object as a matrix. We call these matrices unfolding matrices. In particular, for $k=1,\cdots,d-1$, we define the $k$-th unfolding matrix by $u(x_1,\cdots,x_k;x_{k+1},\cdots,x_d)$ or $u(x_{1:k};x_{k+1:d})$, which is obtained by grouping the first $k$ and last $d-k$ variables to form rows and columns, respectively. 

\subsection{Tensor-network operations}

In this subsection, we introduce several tensor network operations of high importance in our applications using the example of MPS/TT. Similar operators can be extended to more general tensor networks. 

\subsubsection{Marginalization}
\label{sec:marginalization}

To marginalize the MPS/TT representation of $u$ defined in Definition~\ref{def:mps}, one can perform direct operations on each node $\ccG_k$. For instance, if the goal is to integrate out a specific variable $x_k$, the operation can be achieved by taking the summation:

\begin{equation}
    \sum_{x_k} \ccG_k(\alpha_{k-1},x_k,\alpha_{k}).
\end{equation}
The overall computational cost of the marginalization process is at most $O(d)$, depending on the number of variables that need to be integrated out.

\subsubsection{Normalization}
\label{sec:normalization}
Often, one needs to compute the norm $$\|u\|_2^2 = \sum_{x_1,\ldots,x_d} u(x_1,\ldots,x_d)^2,$$ for a MPS/TT. One can again accomplish this via operations on each node. In particular, one can first form the Hadamard product $u\odot u$ in terms of two MPS/TTs, and then integrate out all variables $x_1,\ldots,x_d$ of $u\odot u$ to get $\|u\|_2^2$ using marginalization of MPS/TT as described in \secref{sec:marginalization}. The complexity of forming the Hadamard product of two MPS/TTs is $O(d)$, as mentioned in \cite{oseledets2011tensor}. Therefore, the overall complexity of computing $\|u\|_2^2$ using the described approach is also $O(d)$.

\subsubsection{Sampling From MPS/TT Parametrized Probability Density}
\label{sec:conditional_sampling}

If given a density function $u(x_1,x_2,\dots,x_d)$ in MPS/TT format, it is possible to exploit the linear algebra structure to draw independent and identically distributed (i.i.d.) samples in $O(d)$ time \cite{dolgov2020approximation}, thereby obtaining a sample $(y_1,y_2,\dots,y_d)\sim u$. This approach is derived from the following identity:
\begin{align}\label{eq:conditionals}
    u(x_1,x_2,\dots,x_d) = u(x_1) u(x_2|x_1) u(x_3|x_2,x_1) \dots u(x_d|x_{d-1},\dots,x_{1}),
\end{align}
where each $$u(x_i|x_{i-1},\dots,x_1) = \frac{u(x_{1:i})}{u(x_{1:i-1})},$$  is a conditional distribution of $u$. We note that it is easy to obtain the marginals $u(x_1),\ldots u(x_{1:d-1})$ (hence the conditionals) in $O(d)$ time. If we have such a decomposition \eqref{eq:conditionals}, we can draw a sample $y$ in $O(d)$ complexity as follows. We first draw component $y_1\sim u(x_1)$. Then we move on and draw $y_2\sim u(x_2\vert x_1=y_1)$. We continue this procedure until we draw $y_d \sim u_d(x_d|x_{1:d-1}=y_{1:d-1})$. 


\subsection{Tensor-network sketching}
\label{sec:sketching}

In this subsection, we present a parallel method for obtaining the tensor cores for representing a $d$-dimensional function $u(x)$ that is discretely valued, i.e. each $x_j$ takes on finite values in a set $X_j$, as an MPS/TT. This is done via an MPS/TT sketching technique proposed in \cite{hur2022generative} where the key idea is to solve a sequence of core determining equations. Let  $u_{k}:X_{1}\times \cdots \times X_{k-1} \times \Gamma_{k-1} \rightarrow \bbR, k=2,\ldots,d$ be a set of functions such that $\text{Range}(u_k(x_{1:k};\gamma_k))=\text{Range}(p(x_{1:k};x_{k+1:d}))$. In this case, a representation of $u$ as Definition~\ref{def:mps} can be obtained via solving  $\ccG_k$ from the following set of equations
\begin{align}\label{eq:CDE}
    u_1(x_1,\gamma_1) &= \ccG_1(x_1,\gamma_1),\\
    u_k(x_{1:k}, \gamma_k) &= \sum_{\gamma_{k-1} \in \Gamma_{k-1}} u_{k-1}(x_{1:k-1}, \gamma_{k-1}) \ccG_k(\gamma_{k-1},x_k,\gamma_k),\cr
    u(x) &=  \sum_{\gamma_{d-1} \in \Gamma_{d-1}} u_{d-1}(x_{1:d-1}, \gamma_{d-1}) \ccG_k(\gamma_{d-1},x_d), \notag
\end{align}
based on knowledge of $u_k$'s. 

However \eqref{eq:CDE} is still inefficient to solve since each $u_k$ is exponentially sized, moreover, such a size prohibits it to be obtained/estimated in practice. Notice that \eqref{eq:CDE} is over-determined, we further reduce the row dimensions by applying a left sketching function to \eqref{eq:CDE},
\begin{align}\label{eq:left_sketch_CDE}
    &\sum_{x_1} \cdots \sum_{x_{k-1}} S_{k-1}(x_{1:k-1},\xi_{k-1}) u_k(x_{1:k}, \gamma_k) \\
    = &\sum_{\gamma_{k-1} \in \Gamma_{k-1}} \left( \sum_{x_{1:k-1}}  S_{k-1}( x_{1:k-1},\xi_{k-1}) u_{k-1}(x_{1:k-1}, \gamma_{k-1}) \right) \ccG_k(\gamma_{k-1},x_k,\gamma_k), \notag
\end{align}
where $S_{k-1}:X_{1}\times \cdots \times X_{k-1} \times \Xi_{k-1} \rightarrow \bbR$ is the left sketching function which compresses over variables $x_{1},\cdots,x_{k-1}$.

Now to obtain $u_k$ where $\text{Range}(u_k(x_{1:k};\gamma_k))=\text{Range}(u(x_{1:k};x_{k+1:d}))$, we use a right-sketching  by sketching the dimensions $x_{k+1:d}$, i.e.
\begin{align}
    u_k (x_{1:k}, \gamma_k) = \sum_{x_{k+1:d}}  u(x_{1:k}, x_{k+1:d}) T_{k+1} (x_{k+1:d}, \gamma_{k}),
    \label{eq:right_sketch}
\end{align}
where $T_{k+1}:X_{k+1}\times \cdots \times X_d \times \Gamma_{k} \rightarrow \bbR$ is the right sketching function which compresses $u$ by contracting out variables $x_{k+1},\cdots,x_d$. Plugging such a $u_k$ into \eqref{eq:CDE}, we get
\begin{align}
     B_k[u](\xi_{k-1},x_k,\gamma_{k}) = \sum_{\gamma_{k-1} \in \Gamma_{k-1}} A_{k}[u](\xi_{k-1},\gamma_{k-1}) {\ccG_k} (\gamma_{k-1},x_k,\gamma_{k}),
    \label{eq:both_sketch_CDE}
\end{align}
where
\begin{align}\label{eq:system_A}
A_{k}[u](\xi_{k-1},\gamma_{k-1}) &=\sum_{x_{1:k-1}} \sum_{x_{k:d}}  S_{k-1}( x_{1:k-1},\xi_{k-1}) u(x_{1:k-1}, x_{k:d})  T_{k}(x_{k:d}, \gamma_{k-1}), 
\end{align}
\begin{align}\label{eq:system_B}
B_k[u](\xi_{k-1},x_k,\gamma_{k}) &= \sum_{x_{1:k-1}}  \sum_{x_{k+1:d}} S_{k-1}( x_{1:k-1},\xi_{k-1}) u(x_{1:k-1},x_k, x_{k+1:d}) T_{k+1}(x_{k+1:d}, \gamma_{k}),
\end{align}
and we can readily solve for $\ccG_k$.


Many different types of sketch functions can be used, e.g. random tensor sketches or cluster basis sketches \cite{ahle2019almost,wang2015fast}. Take a single $S_k(x_{1:k},\xi_k)$ as an example, we choose a separable form  $S_k(x_{1:k},\xi_k) = h_1(x_1)\cdots h_k(x_k)$ for some $h_1,\ldots,h_k$ (and similarly for $T_k$'s) in order to perform fast tensor operations. When the state space is discrete, random tensor sketch amounts to taking $h_1,\cdots,h_k$ to be random vectors of size $\vert X_1 \vert,\ldots, \vert X_k\vert$. 

We give special focus to cluster basis sketch, defined as the following:
\begin{defn}\label{def:cluster}
Let $\{b_l\}_{l=1}^n$ be a set of single variable basis. A cluster basis sketch with $c$-cluster consists of choosing $S_k(x_{1:k},\xi_k)\in \left\{b_{l_1}(x_{i_1})\cdots b_{l_c}(x_{i_c})\ \vert\ (l_1,\ldots,l_c)\in [n]^c, \ \{x_{i_1},\ldots,x_{i_c}\}\subseteq \{x_1,\ldots,x_k\}\right\}$, and also   $T_k(x_{k:d},\gamma_{k-1})\in \left\{b_{l_1}(x_{i_1})\cdots b_{l_c}(x_{i_c})\ \vert\ (l_1,\ldots,l_c)\in [n]^c, \ \{x_{i_1},\ldots,x_{i_c}\}\subseteq \{x_k,\ldots,x_d\}\right\}$. For convenience, $S_k, T_k$'s are chosen to be orthogonal basis.
\end{defn}
A similar construct is used in \cite{peng2023generative}. In this case, $A_k$ is of size estimate ${k-1 \choose c}n^c\times {d-k+1 \choose c}n^c$, and $B_k$ is of size ${k-1 \choose c}n^c\times \vert X_k\vert \times {d-k+2 \choose c}n^c$. In principle, we only need ${d-k+2 \choose c}n^c>r$ for determining a rank-$r$ MPS/TT, therefore in practice $c\leq 2$ meets the need. \emph{The most important property we look for is that the variance $A_k[\hat u], B_k[\hat u]$ is small, if $\hat u$ is an unbiased estimator of $u$.} As we shall see in Section~\ref{section:convergence}, when $u$ is a density and $\hat u$ is an empirical distribution with $N$ samples approximating $u$, the choice of such a cluster basis results in the standard deviation of $A_k[\hat u],B_k[\hat u]$ being $O\left(\frac{n^{2c+1}}{\sqrt{N}}\right)$. In contrast, if one chooses a randomized function that involves all $ d$ variables as a sketch, the variance would be $O(n^d/\sqrt{N})$.

\begin{figure}[htb]
    \centering
    \begin{subfigure}{0.25\textwidth}
        \centering
        \vbox{\includegraphics[width=\textwidth]{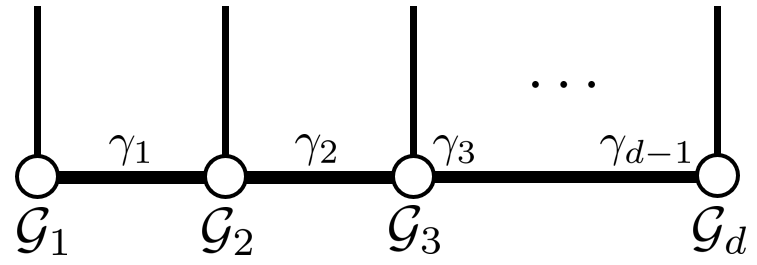}
        \vspace{-0.6em}}
        \caption{Original MPS/TT}
        \label{fig:SVT_1}
    \end{subfigure}
    \begin{subfigure}{0.34\textwidth}
        \centering
        \vbox{\includegraphics[width=\textwidth]{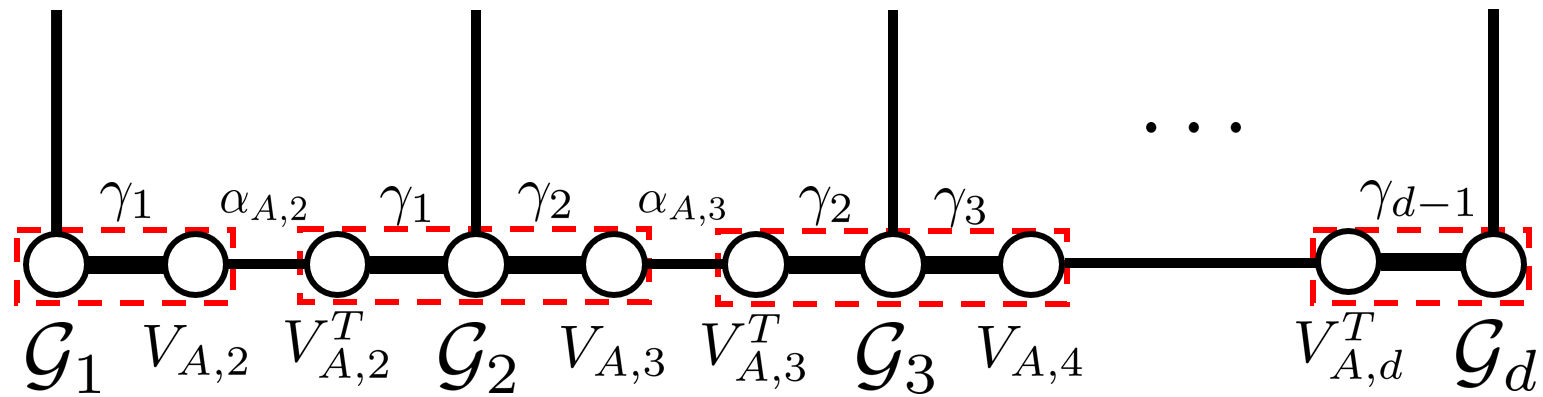}
        \vspace{-0.6em}}
        \caption{Insert projectors and regrouping}
        \label{fig:SVT_2}
    \end{subfigure}
    \begin{subfigure}{0.25\textwidth}
        \centering
        \vbox{\includegraphics[width=\textwidth]{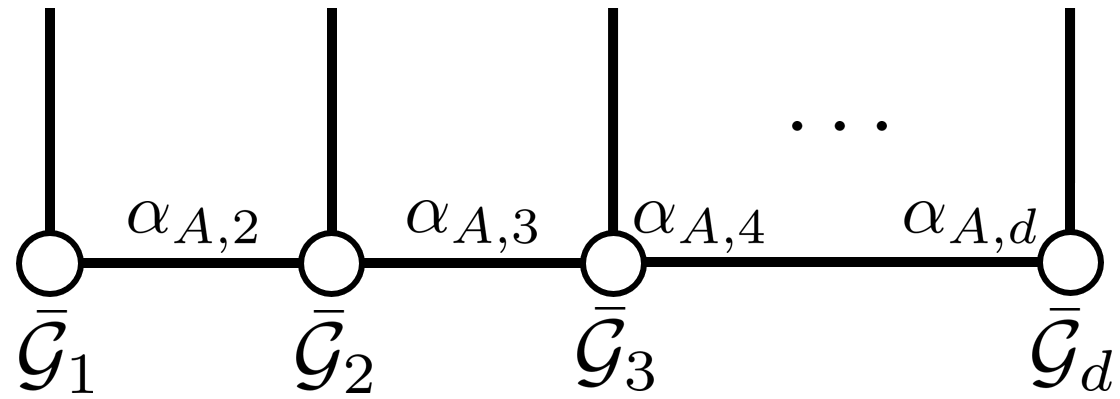}
        \vspace{-0.6em}}
        \caption{Reduced MPS/TT}
        \label{fig:SVT_3}
    \end{subfigure}
    \caption{Tensor diagrams to reduce the bond dimensions of MPS/TT via truncated SVD. The regrouping operation for each reduced tensor $\bar{\ccG}_k$ is highlighted by red dashed boxes.}
    \label{fig:SVT}
\end{figure}

\begin{rem}[Rank of the MPS/TT representation]
   The rank of the MPS/TT may be too large due to oversketching. For instance, when using a cluster basis sketch with a fixed cluster size, the core determining matrices can grow polynomially with the total number of dimensions. However, the intrinsic rank of the MPS/TT may be small. To address this issue, we can utilize truncated singular value decompositions (SVD) of the matrices $\{A_k[u]\}_{k=2}^d$ to define projectors that reduce the rank of the MPS/TT representation. Let
    \begin{align}\label{eq:Adef}
        &A_{k}[u](\xi_{k-1},\gamma_{k-1}) = \sum_{\alpha_{A,k-1}=1}^{r_{A,k-1}} U_{A,k}(\xi_{k-1}, \alpha_{A,k-1}) S_{A,k}(\alpha_{A,k-1}, \alpha_{A,k-1}) V_{A,k}^T(\alpha_{A,k-1}, \gamma_{k-1}),\ \ k=2,\dots,d
    \end{align}
    be the truncated SVD with rank $r_{A,k-1}$ for $A_{k}[u]$. By solving \eqref{eq:both_sketch_CDE} with such a truncation, we obtain an MPS/TT with tensor cores $\{\ccG_k\}_{k=1}^d$ (see \figref{fig:SVT_1}). We can further insert projectors $\{V_{A,k} V_{A,k}^T\}_{k=2}^d$ between all cores to get the ``trimmed'' MPS/TT, as shown in \figref{fig:SVT_2}. We use thick legs to denote dimensions with a large number of indices ($\gamma_k$'s) and thin legs to denote dimensions with a small number of indices ($\alpha_{A,k}$'s). Then we can redefine the reduced tensor cores by grouping the tensor nodes 
    \begin{align}
        &\bar{\ccG}_k \coloneqq \sum_{\gamma_{k-1}}  \sum_{\gamma_{k}} V_{A,k}^T(\alpha_{A,k-1}, \gamma_{k-1}) \ccG_k(\gamma_{k-1},x_k,\gamma_{k}) V_{A,k+1}(\gamma_{k}, \alpha_{A,k}), \ \ k=2,\dots,d-1, \notag\\
        &\text{and}\  \bar{\ccG}_1 \coloneqq \sum_{\gamma_{1}} \ccG_1(x_1, \gamma_{1}) V_{A,2}(\gamma_{1}, \alpha_{A,1}), \ \bar{\ccG}_d \coloneqq \sum_{\gamma_{d-1}}  V_{A,d}^T(\alpha_{A,d}, \gamma_{d}) \ccG_d(x_d,\gamma_{d-1}).
    \end{align}
    The regrouping operations are highlighted with red dashed boxes in \figref{fig:SVT_2}. Now the new tensor core $\bar{\ccG}_k$ is of shape $r_{A,k-1}\times \vert X_k \vert \times r_{A,k}$ (\figref{fig:SVT_3}). We reduce the bond dimensions of the original MPS/TT $(\vert \Gamma_1\vert,\dots,\vert \Gamma_{d-1}\vert)$ to $(r_{A,1},\dots,r_{A,{d-1}})$.
    \label{rem:SVT}
\end{rem} 

\section{Proposed Framework}
\label{sec:framework}

In this section, we present a framework for solving time-evolving systems that arise in both statistical and quantum mechanical systems. In many applications, the evolution of a $d$-dimensional physical system in time is described by a PDE of the form,
\begin{align}
    \frac{\partial \phi(x,t)}{\partial t} = -A \phi(x, t),\ t\geq 0,\ \phi(x,0)=\phi_0(x),
    \label{eq:ivp}
\end{align}
where $A$ is a positive-semidefinite operator and has a semigroup $\{\exp{(-At)}\}_t$ \cite{clifford1961algebraic,hollings2009early,howie1995fundamentals}. $x=(x_1,x_2,\cdots,x_d)$ is a $d$-dimensional spatial point. Depending on the problem, $\phi(x, t)$ might be constrained to some sets. The solution of \eqref{eq:ivp} can be obtained by applying the semigroup operator $\exp{(-At)}$ to the initial function, i.e. $\phi(x,t) = \exp(-At) \phi_0(x)$ for all $t\in\bbR$. We are especially interested in obtaining the stationary solution $\phi^* = \phi(\cdot,t), t\rightarrow\infty$. When \eqref{eq:ivp} is a Fokker-Planck equation, $\phi^*$ corresponds to the equilibrium distribution of Langevin dynamics. When \eqref{eq:ivp}  is the imaginary-time evolution of a Schr\"odinger equation, $\phi^*$ is the lowest energy state wavefunction.

Our method alternates between the two steps detailed in Alg.~\ref{alg:the_alg}. Notice that we introduce three versions of the state $\phi$: $\phi_t$ is the ground truth state function at time $t$; $\hat{\phi}_t$ is a particle approximation of $\phi_t$; $\phi_{\theta_t}$ is a tensor-network representation of $\phi_t$. The significance of these two operations can be described as follows. In the first step, we use a particle-based simulation to bypass the need of applying $\exp(-A\delta t)\phi_{\theta_t}$ exactly, which may have a computational cost that scales unfavorably with respect to the dimensionality. In particular, we generate $f(x;x^i)$ that can be represented easily as a low-rank TT. In the second step, the use of the sketching algorithm \cite{hur2022generative} bypasses a direct application of recursive SVD-based compression scheme \cite{oseledets2011tensor}, which may run into $O(dN^2)$ complexity. From a statistical complexity point of view, the variance of an empirical distribution $\text{Var}(\hat \phi_{t+1})$ scales exponentially in $d$, therefore one cannot apply a standard tensor compression scheme to $\hat \phi_{t+1}$ directly since it preserves the exponential statistical variance. Hence it is crucial to use the techniques proposed in \cite{hur2022generative}, which is designed to control the variance of the sketching procedure. 

The method to obtain a Monte Carlo representation of the time-evolution in Step 1 of Alg.~\ref{alg:the_alg} is dependent on the applications, and we present different schemes for many-body Schr\"odinger and Fokker-Planck equation in \secref{sec:applications}. In this section, we focus on the details of the second step. 

\begin{algorithm}
  \caption{Combining tensor-network and Monte Carlo method by sketching}
  \label{alg:the_alg}
  \begin{algorithmic}[1]
    \State Apply the semigroup operator $\exp(-A\delta t)$ to the tensor-network approximation of the solution $\phi_{\theta_t}(x)$ using particle simulations for $\delta t > 0$, i.e.
\begin{align}
    \phi_{t+1}(x) = \exp{(-A\delta t)}\ \phi_{\theta_t}(x) = \mathbb{E}_{y\sim \mu}[f(x;y)] \approx \frac{1}{N} \sum_{i=1}^N f(x;x^i)=:\hat \phi_{t+1}(x),
    \label{eq:sampled_ivp}
\end{align}
where $\{x^i\}_{i=1}^N\subset \mathbb{R}^d$ is a collection of $N$ i.i.d. samples according to a distribution $\mu$ (depending on the application, see Section~\ref{sec:applications} and $f(x;x^i)$ is a $d$-dimensional function with parameters $x^i$. In the traditional Monte Carlo simulation, $f$ is simply the Dirac delta function at the given sample point, i.e. $f(x;x^i)=\delta (x-x^i)$.
    \State Estimate $\phi_{t+1}(x)$ as a tensor-network $\phi_{\theta_{t+1}}(x)$ from its particle approximations $\hat \phi_{t+1}(x)$ via the parallel TT-sketching method with linear time complexity with respect to the number of samples and constant time with respect to the dimension if distributed computing is used (\secref{sec:sketching}). Often certain normalization constraints (with respect to a certain norm $\|\cdot\|$) need to be enforced for $\phi_{\theta_{t+1}}$. In this case one simply adds an extra step, letting $\phi_{\theta_{t+1}}\leftarrow \frac{\phi_{\theta_{t+1}}}{\|\phi_{\theta_{t+1}}\|_2}$, which can be done with $O(d)$ complexity (\secref{sec:normalization}).
  \end{algorithmic}
\end{algorithm}





\subsection{MPS/TT sketching for random walkers represented as a sum of TT}

The concept of MPS/TT sketching or general tensor network sketching for density estimations has been extensively explored in \cite{hur2022generative,ren2022high,tang2022generative}. The choice of tensor network representation in practice depends on the problem's structure. In this work, we demonstrate the workflow using MPS/TT, but the framework can be readily extended to other tensor networks.

A crucial assumption underlying our approach is that each particle $f(x;x^i)$ exhibits a simple structure and can be efficiently represented or approximated in MPS/TT format. For instance, in the case of a vanilla Markov Chain Monte Carlo (MCMC), we have $f(x;x^i)=\delta(x-x^i)=\prod_{j=1}^d \delta(x_j - x_j^i)$, which can be expressed as a rank-1 MPS/TT. Consequently, the right-hand side of \eqref{eq:sampled_ivp} becomes a summation of $N$ MPS/TTs. The tensor diagram illustrating this general particle approximation is shown in \figref{fig:workflow_step1}. A naive approach for representing this sum as an MPS/TT consists of directly adding these MPS/TTs, which yields an MPS/TT with a rank of $N$ \cite{oseledets2011tensor}. This growth in rank with the number of samples can lead to a high computational complexity of $O(\text{poly}(N))$ when performing tensor operations.

On the other hand, depending on the problem, if we know that the intrinsic MPS/TT rank of the solution is small, it should be possible to represent the MPS/TT with a smaller size. To achieve this, we apply the method described in \secref{sec:sketching} to estimate a low-rank tensor directly from the sum of $N$ particles. Specifically, we construct $A_k[\hat \phi_{t+1}]$ and $B_k[\hat \phi_{t+1}]$ from the empirical samples $\hat \phi_{t+1}$ and use them to solve for the tensor cores through \eqref{eq:both_sketch_CDE}, resulting in $\phi_{\theta_{t+1}}$. This process is illustrated in \figref{fig:workflow_step2}, \figref{fig:workflow_step3}, and \figref{fig:workflow_step4}. It is important to note that we approximate the ground truth $\phi_{t+1}$ using stochastic samples $\hat{\phi}_{t+1}$. The estimated tensor cores form an MPS/TT representation of $\phi_{t+1}$, denoted as $\phi_{\theta_{t+1}}$.

\subsection{Complexity analysis}

In this section, we assume we are given $\hat \phi_{t+1}$ in \eqref{eq:ivp} as $N$ constant rank MPS/TT.  In terms of samples, forming $A_k[\hat \phi_{t+1}])$ and $B_k[\hat \phi_{t+1}]$ in \eqref{eq:system_A} and \eqref{eq:system_B} can be done as 
\begin{align}\label{eq:sampled A B}
    &A_k[ \hat \phi_{t+1}](\xi_{k-1},\gamma_{k-1}) \\
    = &\frac{1}{N}\sum_{i=1}^N \left[ \sum_{x_{1:k-1}} \sum_{x_{k:d}}  S_{k-1}( x_{1:k-1},\xi_{k-1}) f(x_{1:k-1}, x_{k:d};x^i) T_{k}(x_{k:d}, \gamma_{k-1}) \right],\notag \\
    &B_k[\hat \phi_{t+1}](\xi_{k-1},x_k,\gamma_{k}) \notag\\
    = &\frac{1}{N}\sum_{i=1}^N \left( \sum_{x_{1:k-1}}  \sum_{x_{k+1:d}} S_{k-1}( x_{1:k-1},\xi_{k-1}) f(x_{1:k-1},x_k, x_{k+1:d}; x^i)  T_{k+1}(x_{k+1:d}, \gamma_{k}) \right). \notag
\end{align}
$B_k$ is a $3$-tensor of shape $|\Xi_{k-1}|\times |X_k| \times |\Gamma_{k}|$. $A_k$ is a matrix of shape $|\Xi_{k-1}| \times |\Gamma_{k}|$. The size of the linear system is independent of the number of dimensions $d$ and the total number of particles $N$. The number of sketch functions $|\Xi_{k}|$ and $|\Gamma_k|$ are hyperparameters we can control. Let $n=\max_{k} |X_k|$, $\tilde r = \max_k \{|\Xi_{k}|, |\Gamma_k|\}$. First we consider the complexity of forming the core determining equations, i.e. evaluating $A_k$ and $B_k$. We note that the complexity is dominated by evaluating the tensor contractions between sketch functions $S_k$'s, right sketch functions $T_k$'s and samples $f(\cdot;x^i)$'s. If the sample $f(\cdot;x^i)$  and sketch functions $S_k(\cdot,\xi_k)$, $T_k(\cdot,\gamma_{k-1})$ are in low rank MPS/TT representations, then the tensor contractions in the square bracket in \eqref{eq:sampled A B} can be evaluated in $O(d)$ time. Each term in the summation $\sum_{i=1}^N$ can be done independently, potentially even in parallel, giving our final $A_k$'s and $B_k$'s. Taking all $d$ dimensions into account, the total complexity of evaluating $A_k$'s and $B_k$'s is $O(n \tilde r^2 N d)$. 

Next, the complexity of solving the linear system \eqref{eq:both_sketch_CDE} is $O(n\tilde r^3)$. For all $d$ dimensions, the total complexity for solving the core determining equations is $O(n\tilde r^3 d)$. Combing everything together, the total computational complexity for MPS/TT sketching for a discrete particle system is 
\begin{align}
    O(n \tilde r^2 N d) + O(n\tilde r^3 d).
\end{align}

\begin{figure}[htb]
    \centering
    \begin{subfigure}{0.25\textwidth}
        \centering
        \includegraphics[width=\textwidth]{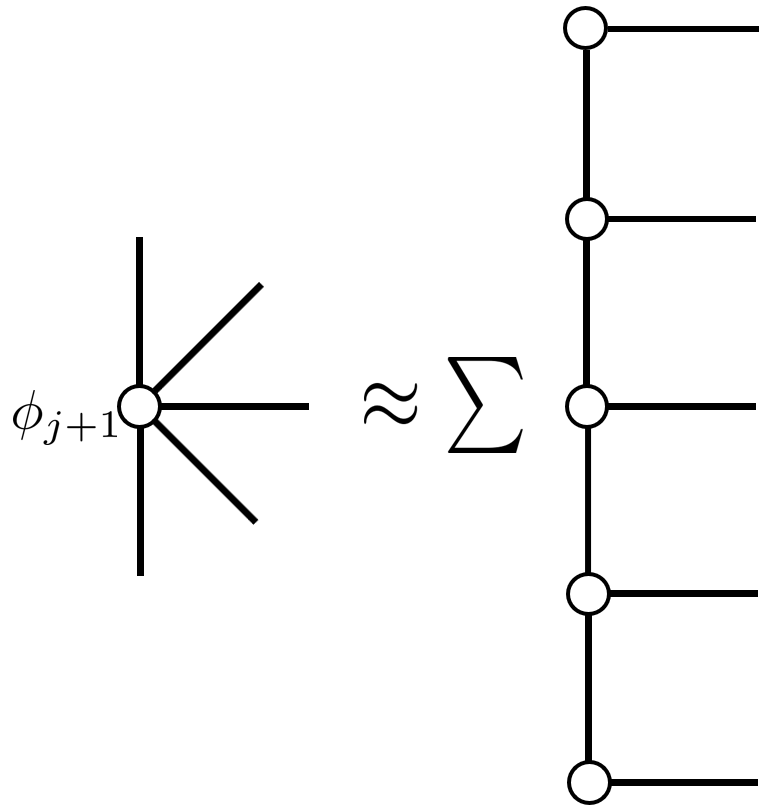}
        \caption{Particle approximation of $\phi_{t+1}$}
        \label{fig:workflow_step1}
    \end{subfigure}\quad \quad
    \begin{subfigure}{0.378\textwidth}
        \centering
        {\vbox{
        \includegraphics[width=\textwidth]{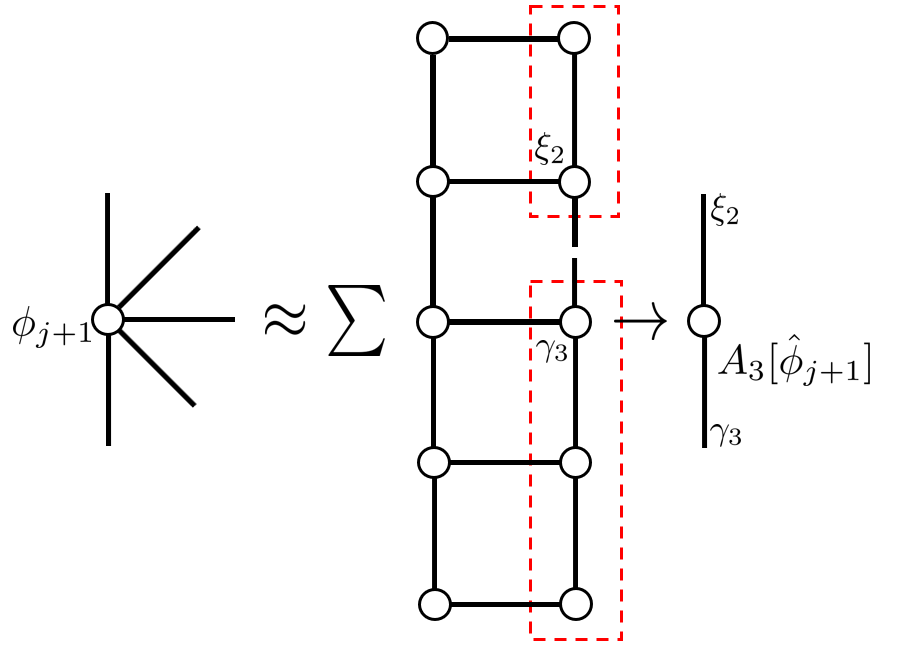}
        \vspace{-0.6em}}}
        \caption{Form $\{A_k\}_{k=1}^d$.}
        \label{fig:workflow_step2}
    \end{subfigure}\\
    \begin{subfigure}{0.358\textwidth}
        \centering
        \includegraphics[width=\textwidth]{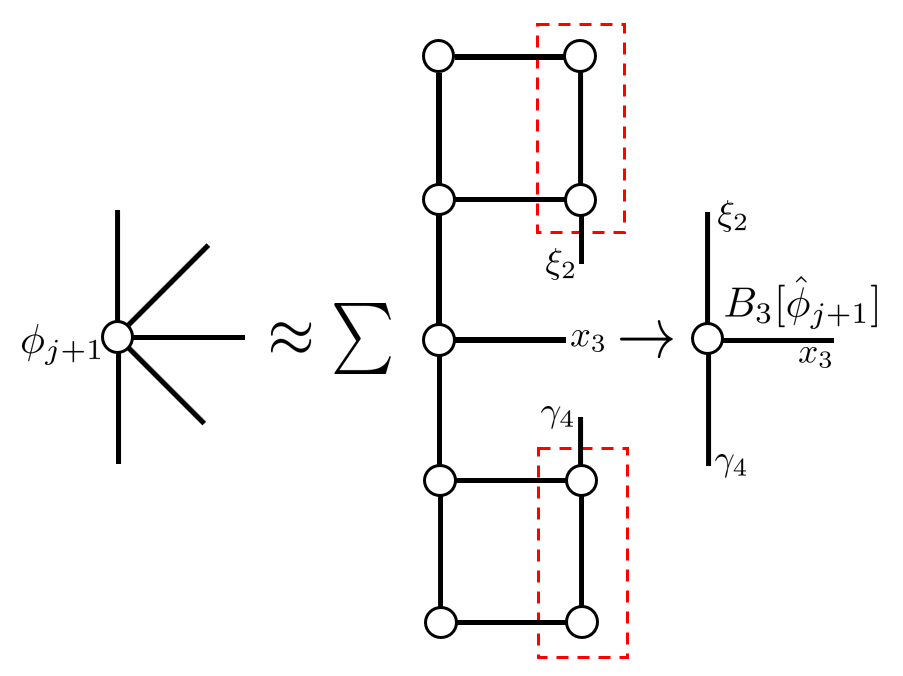}
        \caption{Form $\{B_k\}_{k=1}^d$.}
        \label{fig:workflow_step3}
    \end{subfigure}\quad \quad \quad
    \begin{subfigure}{0.238\textwidth}
        \centering
        \hbox{\includegraphics[width=\textwidth]{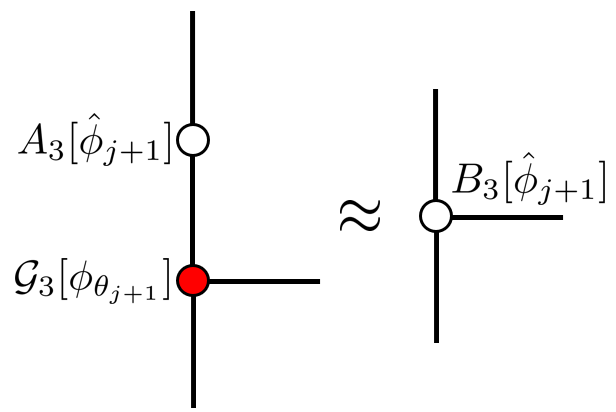}
        \hspace{2.5em}}
        \caption{Solving core determining equations}
        \label{fig:workflow_step4}
    \end{subfigure}
    \caption{Tensor diagram for the workflow of estimating an MPS/TT from particles. Step (a) shows how $\phi_t$ is represented as empirical distribution $\hat \phi_{t+1}$. Step (b), (c), (d) shows how to form $A_k[\hat \phi_{t+1}], B_k[\hat \phi_{t+1}]$ in \eqref{eq:sampled A B} and use them to solve for $\ccG_k$. Here we use the determination of $\ccG_k$ where $k=3$ as an example.}
    \label{fig:workflow}
\end{figure}

\begin{rem}
MPS/TT sketching is a method to control the complexity of the solution via estimating an MPS/TT representation in terms of the particles. Another more conventional approach to reduce the rank of MPS/TT is TT-round \cite{tt-decomposition}. With this approach, one can first compute the summation in \figref{fig:workflow_step1} to form a rank $O(N)$ MPS/TT and round the resulting MPS/TT to a constant rank. There are two main drawbacks of this approach. Firstly, the QR decomposition in TT-round has complexity $O(N^2)$. Secondly, TT-round tries to make $\phi_{\theta_{t+1}}\approx \hat \phi_{t+1}$, which at the same time makes  $\text{Var}( \phi_{\theta_{t+1}}))\approx\text{Var}(\hat  \phi_{t+1})$, and such an error can scale exponentially. 
\end{rem}

\section{Applications}
\label{sec:applications}

In this section, to demonstrate the generality of the proposed method, we show how it can be used in two applications: quantum many-body ground-state problem (Section~\ref{sec:quantum-evolution}) and solving Fokker-Planck equation (Section~\ref{sec:fokker-planck-evolution}). For these applications, we focus on discussing how to approximate $\exp(-\beta A)\phi_{\theta_t}$ as a sum of MPS/TT, when $\phi_{\theta_t}$ is already in an MPS/TT representation, as required in \eqref{eq:ivp}.

\subsection{Quantum many-body ground-state problem}
\label{sec:quantum-evolution}

In this subsection, we apply the proposed framework to ground-state energy estimation problems in quantum mechanics for the spin system. Statistical sampling-based approaches have been widely applied to these types of problems, see e.g. \cite{lee2022unbiasing}. In this problem, we want to solve
\begin{equation}\label{eq:quantum pde}
\frac{\partial \phi(x,t)}{\partial t} = -H \phi(x,t),\quad \|\phi(\cdot,t)\|_2=1
\end{equation}
where $\phi(\cdot,t):\{\pm 1\}^d\rightarrow \mathbb{C}$, and $H$ is the Hamiltonian operator. Let $O:=[O_i]^d_{i=1}$, where
\begin{equation}
    O_i = I_2 \otimes I_2 \otimes \cdots \otimes \underbrace{ \tilde O}_{\text{$i$-th term}} \otimes \cdots \otimes I_2,
\end{equation}
for some $\tilde O \in \mathbb{C}^{2\times 2}$, $O_i^2 = I_{2^d}$, then $H$ usually takes the form of $H = H_1[O] + H_2[O]$, where $H_1[O] = -h\sum_{i=1}^d O_i$, and $H_2[O] = \sum_{i,j=1}^d J_{ij} O_i O_j$. When $t\rightarrow \infty$, $\phi(\cdot,t)=\frac{\exp(-Ht)\phi_0}{\|\exp(-Ht)\phi_0\|_2}$ gives the lowest eigenvector of $H$. This can be done with the framework detailed in Section~\ref{sec:framework} with $A=H$. In Section~\ref{sec:stochastic quantum ite}, we detail how $\exp(-\delta t H)\phi_{\theta_t}$ can be approximated as a sum of $N$ functions $f(x;x^i)$ via a specific version of quantum Monte Carlo, the auxiliary-field quantum Monte Carlo (AFQMC). Alternatively, instead of applying $\exp(-\delta t H)$ via sampling, one can deterministically apply the operator $I - \delta t H$ as the propagator for imaginary-time evolution. For a detailed discussion, see Section~\ref{sec: quantum evolution deterministic}.

\subsubsection{Stochastic Quantum Imaginary-Time Evolution} \label{sec:stochastic quantum ite}

The AFQMC method \cite{afqmc1,afmc2} is a powerful numerical technique that has been developed to overcome some of the limitations of traditional Monte Carlo simulations. AFQMC is based on the idea of introducing auxiliary fields to decouple the correlations between particles by means of the application of the Hubbard–Stratonovich transformation \cite{hubbard1959calculation}. This reduces the many-body problem to the calculation of a sum or integral over all possible auxiliary-field configurations. The method has been successfully applied to a wide range of problems in statistical mechanics, including lattice field theory, quantum chromodynamics, and condensed matter physics \cite{zhang2013auxiliary,carlson2011auxiliary,shi2021some,qin2016benchmark,lee2022twenty}. However, an issue of AFQMC is that the random walkers are biased towards a mean-field solution \cite{qin2016coupling} in order to reduce the variance caused by the ``sign problem'', giving rise to bias when determining the ground-state energy. To solve this issue, \cite{qin2016coupling} alternates between running AFQMC and determining a new mean-field solution in order to remove such a bias. Our approach can be regarded as a generalization to the philosophy in \cite{qin2016coupling}, where the mean-field is now replaced by a more general MPS/TT representation. 

The idea of applying $\exp(-\delta t H)$ on an MPS/TT $\phi_{\theta_t}$ is that we write it as the expectation of a random rank-$1$ MPO:
\begin{equation}\label{eq:discrete hs transformation}
    \exp(-\delta t H) = \mathbb{E}_{\vec{\sigma}\sim P(\vec{\sigma})} B(\vec{\sigma}), 
\end{equation}
where $\vec{\sigma}$ is referred to as the auxiliary fields, and $P(\vec{\sigma})$ is some probability density function. Naively, the propagator $\exp(-\delta t H)$ can be approximated by
\begin{equation}
    \exp(-\delta t H) \approx \frac{1}{N} \sum_{i=1}^N B(\vec{\sigma}^i), 
\end{equation}
where $N$ is the total number of Monte Carlo samples, and $\vec{\sigma}^i$ represents the $i$-th sample. All samples $\{\vec{\sigma}^i\}_{i=1,\cdots,N}$ are drawn from $P(\vec{\sigma})$, i.i.d. 

In practice, we further incorporate importance sampling to reduce the variance. We introduce a new probability function:
\begin{equation}
    \tilde{P}_{\phi_{\theta_t}}(\vec{\sigma}) \propto \max\left\{\braket{B(\vec{\sigma})\phi_{\theta_t},\phi_{\theta_t}},0\right\} P(\vec{\sigma}). 
\end{equation}
The probability function $\tilde{P}_{\phi_{\theta_t}}(\vec{\sigma})$ favors those samples that lead to large overlap with the current wavefunction $\phi_{\theta_t}$. 
With this new probabilty function $\tilde{P}_{\phi_{\theta_t}}$, the propagator $\exp(-\delta t H)$ is then approximated by
\begin{equation}
    \exp(-\delta t H) \approx \frac{1}{N} \sum_{i=1}^N \frac{P(\vec{\sigma}^i)}{\tilde{P}(\vec{\sigma}^i)}B(\vec{\sigma}^i), 
\end{equation}
where $\{\vec{\sigma}^i\}_{i=1,\cdots,N}$ are sampled from $\tilde{P}_{\phi_{\theta_t}}$, i.i.d. 

To construct the decomposition \eqref{eq:discrete hs transformation}, we use a Suzuki-Trotter approximation \cite{suzuki1976relationship,trotter1959product}:
\begin{equation}\label{eq:trotter decomposition}
\begin{split}
    \exp(-\delta t H) &= \exp(-\delta t (H_1[O]+H_2[O])) \\
    &= \exp(-\delta t H_1[O]/2) \exp(-\delta t H_2[O]) \exp(-\delta t H_1[O]/2) + O(\delta t^3). 
\end{split}
\end{equation}
For the term $\exp(-\delta t H_1[O]/2)$, where $H_1[O]=-h\sum_{i=1}^d O_i$, we simply have
\begin{multline}\label{eq:one body prop}
    \exp(-\delta t H_1[O]/2) = \exp(\delta t h\sum_{i=1}^d O_i/2) = \prod_{i=1}^d \exp(-\delta t h O_i/2) \\
    = \exp(\delta t h \tilde{O}/2) \otimes \cdots \otimes \exp(\delta t h \tilde{O}/2), 
\end{multline}
which is already a rank-$1$ MPO. For the term $\exp(-\delta t H_2[O])$, where $H_2[O] = \sum^d_{i,j=1} J_{ij}O_i O_j$, the trick is to use a discrete Hubbard-Stratonovich transformation \cite{ulmke2000auxiliary}:
\begin{equation}
    \exp(-\delta t J_{ij} O_iO_j) = \frac{1}{2} e^{-\delta t J_{ij}} \sum_{\sigma_{ij}=\pm 1} \exp(\delta t \lambda_{ij} \sigma_{ij} (O_i-O_j)), \ \text{for } J_{ij} > 0, 
\end{equation}
and
\begin{equation}\label{eq:discrete hs repulsive}
    \exp(-\delta t J_{ij} O_iO_j) = \frac{1}{2} e^{\delta t J_{ij}} \sum_{\sigma_{ij}=\pm 1} \exp(\delta t \lambda_{ij} \sigma_{ij} (O_i+O_j)), \ \text{for } J_{ij} < 0, 
\end{equation}
where the constants $\lambda_{ij}$ are given by $\cosh{2\lambda_{ij}} = e^{2\delta t}$. For both cases, we let $p_{ij}=1/2$ being the probability function, and 
\begin{equation}
    b_{ij}(\sigma_{ij}) = e^{-\delta t J_{ij}} \exp(\delta t \lambda_{ij} \sigma_{ij} (O_i-O_j)), \ \text{for } J_{ij} > 0
\end{equation}
or 
\begin{equation}
    b_{ij}(\sigma_{ij}) = e^{\delta t J_{ij}} \exp(\delta t \lambda_{ij} \sigma_{ij} (O_i+O_j)), \ \text{for } J_{ij} < 0
\end{equation}
being the rank-$1$ MPO, we can then rewrite the decomposition as
\begin{equation}
    \exp(-\delta t J_{ij} O_iO_j) = \sum_{\sigma_{ij}=\pm1} p_{ij}(\sigma_{ij}) b_{ij}(\sigma_{ij}). 
\end{equation}
For the whole term $\exp(-\delta t H_2[O])$, we then have
\begin{multline}\label{eq:two body prop}
    \exp(-\delta t H_2[O]) = \exp(-\delta t \sum^d_{i,j=1} J_{ij}O_i O_j)
    = \prod_{i,j=1}^d \exp(-\delta t J_{ij} O_iO_j) + O(\delta t^2) \\
    = \prod_{i,j=1}^d \sum_{\sigma_{ij}=\pm1} p_{ij}(\sigma_{ij}) b_{ij}(\sigma_{ij}) + O(\delta t^2)
    = \sum_{\vec{\sigma}} P(\vec{\sigma}) B_{H_2}(\vec{\sigma}) + O(\delta t^2), 
\end{multline}
where $\vec{\sigma}$ denotes the collection of all the auxiliary fields $\{\sigma_{ij}\}_{i,j=1,\cdots,d}$. The probability function is
\begin{equation}
    P(\vec{\sigma}) = \prod_{i,j=1}^d p_{ij}(\sigma_{ij}), 
\end{equation}
and the rank-$1$ MPO is
\begin{equation}
    B_{H_2}(\vec{\sigma}) = \prod_{i,j=1}^d b_{ij}(\sigma_{ij}). 
\end{equation}
In general, an extra Trotter error $O(\delta t^2)$ is introduced. Plugging \eqref{eq:one body prop} and \eqref{eq:two body prop} into \eqref{eq:trotter decomposition} gives
\begin{equation}
    \exp(-\delta t H) = \sum_{\vec{\sigma}} P(\vec{\sigma}) B(\vec{\sigma}) + O(\delta t^2) = \mathbb{E}_ {{\vec{\sigma}}\sim P(\vec{\sigma})} B(\vec{\sigma}) + O(\delta t^2), 
\end{equation}
which is exactly the desired form, with the merged rank-$1$ MPO being
\begin{equation}
    B(\vec{\sigma}) = \exp(-\delta t H_1[O]/2) B_{H_2}(\vec{\sigma}) \exp(-\delta t H_1[O]/2). 
\end{equation}

\begin{figure}[!htbp]
\centering
\includegraphics[width=0.6\linewidth]{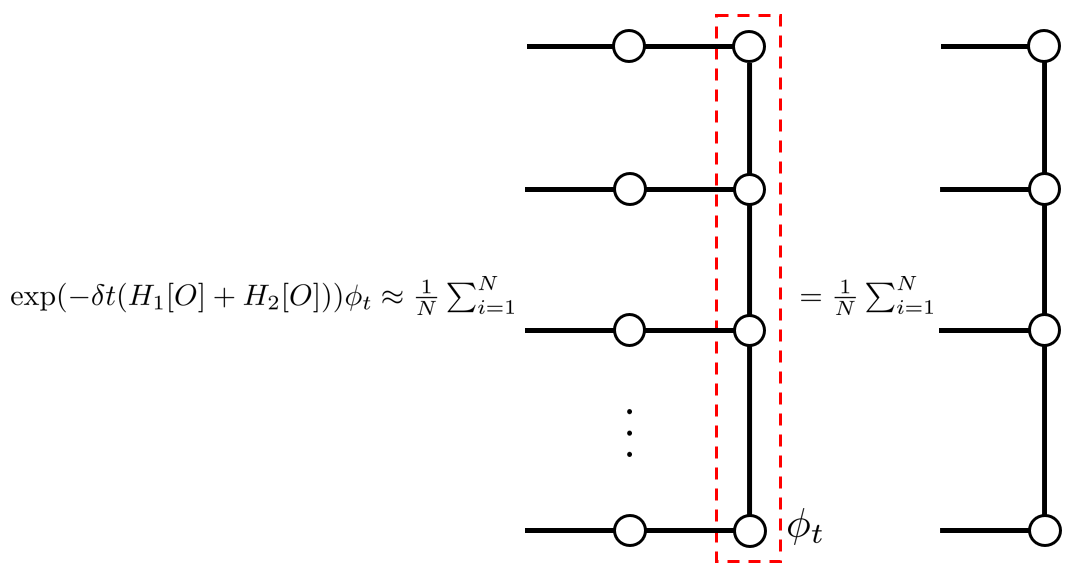}\\
\vspace{-0.25cm}
\caption{Tensor diagram for approximating evolution equation as MPO-MPS products. We remove the internal legs connecting the tensor cores in MPO to indicate the MPO has rank $1$. }
\label{fig:afqmc_MPO}
\end{figure}

To summarize, we can approximate the imaginary time propagator $\exp(-\delta t (H_1[O] + H_2[O]))$ with a summation of $N$ rank-$1$ MPOs. Applying the propagator to many-body wavefunction $\phi_t$ reduces to MPO-MPS contractions. The corresponding tensor diagram is shown in \figref{fig:afqmc_MPO}.

\subsubsection{Deterministic quantum imaginary-time evolution}\label{sec: quantum evolution deterministic}

Instead of applying the propagator $\exp(-\delta t H)$, an alternative approach to solving \eqref{eq:quantum pde} is to apply the operator $I - \delta t H$. The two operators become equivalent in the limit as $\delta t \to 0$. 
We recall that $H$ takes the form
\begin{equation}\label{eq:hamiltonian explicit}
    H = H_1[O] + H_2[O] = -h\sum_{i=1}^d O_i + \sum_{i,j=1}^d J_{ij} O_i O_j, 
\end{equation}
of which each term in the right-hand-side is a rank-$1$ MPO. The ideneity operator $I$ can also be treated as a rank-$1$ MPO. Putting them together, we can rewrite the propagator $I-\delta t H$ as a sum of some rank-$1$ MPOs: 
\begin{equation}
    I - \delta t H = \sum_{i=1}^L H^{(i)},  
\end{equation}
where each $H^{(i)}$ represents an MPO. Applying $I - \delta t H$ to the wavefunction $\phi_{\theta_t}$ again amounts to MPO-MPS contractions. 
\begin{equation}\label{eq:deterministic sum}
    (I - \delta t H)\phi_{\theta_t} = \sum_{i=1}^L H^{(i)} \phi_{\theta_t}. 
\end{equation}
On the one hand, deterministic evolution offers a clear advantage over stochastic evolution, as the propagation is exact. On the other hand, however, it suffers from high computational cost. Specifically, the number of MPOs in the sum \eqref{eq:deterministic sum} scales with the number of sites, leading to a computational cost that is at least quadratic in the system size. Moreover, in cases where the system is highly connected, i.e., when the adjacency matrix $J$ in \eqref{eq:tfi hamiltonian} is dense, the number of terms in the sum \eqref{eq:deterministic sum} increases even further, making the algorithm significantly more expensive.

\subsection{Fokker-Planck equation}
\label{sec:fokker-planck-evolution}

In this subsection, we demonstrate the application of the proposed framework to numerical simulations of parabolic PDEs, specifically focusing on the overdamped Langevin process and its corresponding Fokker-Planck equations. We consider a particle system governed by the following overdamped Langevin process,
\begin{align}
    d x_t = -\nabla V(x_t) \, dt + \sqrt{2\beta^{-1}} \, d W_t,
    \label{eq:Langevin}
\end{align}
where $x_t \in \Omega\subseteq\Rd$ is the state of the system, $V:\Omega\subset\Rd\rightarrow\bbR$ is a smooth potential energy function, $\beta=1/T$ is the inverse of the temperature $T$, and $W_t$ is a  $d$-dimensional Wiener processs. If the potential energy function $V$ is confining for $\Omega$ (see, e.g., \cite[Definition 4.2]{bhattacharya2009stochastic}), it can be shown that the equilibrium probability distribution of the Langevin dynamics \eqref{eq:Langevin} is the Boltzmann-Gibbs distribution,
\begin{align}
    {\phi^*}(x) = \frac{1}{Z_\beta}\exp(-\beta V(x))
    \label{eq:boltzmann}
\end{align}
where $Z_\beta=\int_{\Omega} \exp(-\beta V(x))\, d x$ is the partition function. Moreover, the evolution of the distribution of the particle system can be described by the corresponding time-dependent Fokker-Planck equation, 
\begin{align}
    \pfrac{\phi}{t} = \beta^{-1} \Delta \phi + \nabla \cdot (\nabla V \phi)=:-A\phi,\quad \phi(x,0)=\phi_0(x),\quad \|\phi(\cdot,t)\|_1=1
    \label{eq:fokker-planck}
\end{align}
where $\phi_0$ is the initial distribution. $\|\phi(\cdot,t)\|_1=1$ ensures that the $\int \vert\phi(x,t)\vert dx=1$. Therefore, \eqref{eq:fokker-planck} is the counterpart of \eqref{eq:ivp} in our framework.

Now we need to be able to approximate $\exp(-A\delta t)\phi_t$ as particle systems. Assuming the current density $\phi_t$ is a MPS/TT, we use the following procedures to generate a particle approximation $\hat \phi_{t+1}$:
\begin{enumerate}
\item We apply conditional sampling on the current density estimate $\phi_{\theta_t}$ (\secref{sec:conditional_sampling}) to generate $N$ i.i.d. samples $x^1,\ldots, x^N\sim \phi_{\theta_t}$.
\item Then, we simulate the overdamped Langevin dynamics \eqref{eq:Langevin} using Euler-Maruyama method over time interval $\delta t$ for each of the $N$ initial stochastic samples  $x^1,\ldots x^N\sim \phi_{\theta_t}$. By the end of $\delta t$ we have final particle positions $x^1,\ldots,x^N\sim \phi_{t+1}$, and
\begin{align}
    \hat \phi_{t+1}(x) &= \frac{1}{N} \sum_{i=1}^N \delta(x - x^i),
    \label{eq:sampled_fokker_planck}
\end{align}
\end{enumerate}
by standard Monte Carlo approximation. Note that the only difference between this application and the quantum ground-state problem is the conversion of the empirical distribution into an MPS/TT $\phi_{\theta_{t+1}}(x)$. Here, we employ a version of sketching for continuous distributions instead of discrete distributions as used in quantum ground-state energy estimation (\secref{sec:quantum-evolution}). For more detailed information on MPS/TT sketching for continuous distributions, we refer readers to Appendix C of \cite{hur2022generative}.

\section{Convergence Analysis}\label{section:convergence}

In this section, we provide the convergence analysis for the proposed method. We look at the case for the Fokker-Planck equation, in a simplified discretized setting. Let $P_{\delta t}\in \mathbb{R}^{n^d\times n^d}$ be a Markov-transition kernel of a stochastic process on a discrete state space $[n]^d$. Denote the stationary distribution as $\phi^\star$, which satisfies $P_{\delta t} \phi^\star = \phi^\star$. We want to show that Alg.~\ref{alg:the_alg} converges to $ \phi^\star$. To facilitate the discussion, we define a few new notations. For a $d$-tensor $u$ of size $n^d$, we define its ``Frobenius norm'' to be
\begin{equation}
\|u\|_F := \sqrt{\sum_{i_1,\cdots, i_d} u(i_1,\cdots,i_d)^2}.
\end{equation}
Furthermore, when representing $u$ as an MPS/TT with cores $\ccG_1,\cdots,\ccG_d$, we use the notation
\begin{equation}
u = \ccG_1\circ \ccG_2\cdots \circ \ccG_d.
\end{equation}
We are also going to use a standard perturbation theory result for the solution of a linear system.
\begin{lemma}[Theorem 3.48, \cite{wendland2017numerical}]\label{lemma:LS stability}
Suppose $A x^\star = b$, $A\in \mathbb{R}^{n\times n}, b\in  \mathbb{R}^{n\times 1}$. Further $(A+\Delta A)x = (b+\Delta b)$. Then with $\| A^\dagger \|_2\|\Delta A\|_2\leq 1$, we have
\begin{equation}
\frac{\|x-x^\star\|}{\|x^\star\|}\leq \frac{\|A^\dagger\|_2}{1-\| A^\dagger \|_2\|\Delta A\|_2}(\|\Delta A\|_2+  \|\Delta b\|_2/\|x^\star\|).
\end{equation}
\end{lemma}

Our main theorem is stated in Theorem~\ref{thm:main}, which shows that the iterates in Alg.~\ref{alg:the_alg} are contracting towards the true solution, perturbed by some error.  To prove it, we make the following assumption.
\begin{assumption}\label{assumption:spectral gap}
Let $P_{\delta t}\in \mathbb{R}^{n^d\times n^d}$ be a Markov-transition kernel of a stochastic process on a discrete state space $[n]^d$. We assume that it has eigenvalues $1=\lambda_1^{\delta t}>\lambda_2^{\delta t}\geq \lambda_3^{\delta t}\cdots $, and therefore has a unique top eigenvector $P_{\delta t} \phi^\star = \phi^\star$. Furthermore, $\min(\phi^\star)>0$.
\end{assumption}
Such an assumption is important to characterize the contraction rate towards $\phi^\star$ in Alg~\ref{alg:the_alg}:
\begin{lemma}\label{lemma:contraction with error}
Suppose in Alg.~\ref{alg:the_alg}, $\|\phi_{\theta_t}- \phi_{t}\|_F \leq \tilde \nu$, then 
\begin{equation}
\|\phi_{\theta_t} - \phi^\star \|_F \leq a^t \| \phi_{\theta_0} - \phi^\star \|_F + \tilde \nu
\end{equation}
where $a := \frac{\max{\phi^\star}}{\min{\phi^\star}}\lambda_2^{\delta t}$ for $\delta t$ large enough such that $a<1$.
\end{lemma}
\begin{proof}
We first recall several definitions in Alg.~\ref{alg:the_alg}. Let $\phi_{t+1} = P_{\delta t}\phi_{\theta_t}$, $\hat \phi_{t+1}$ be an empirical distribution of $\phi_{t+1}$ with $N$ samples, which satisfies $\mathbb{E}(\hat \phi_{t+1}) = \phi_{t+1}$.  Let $\phi_{\theta_{t+1}}$ be an MPS/TT approximation of $\phi_{t+1}$, obtained from plugging $\hat \phi_{t+1}$ into \eqref{eq:sampled A B}. 

First, under Assumption~\ref{assumption:spectral gap}, standard results in Markov process (see for example \cite{chen2000equivalence}) give 
\begin{equation}
\| \phi_{t+1} - \phi^\star\|_{2,{\phi^\star}^{-1}} =  \|P_{\delta t} \phi_{\theta_t} -\phi^\star\|_{2,{\phi^\star}^{-1}}\leq \lambda_2^{\delta t} \|  \phi_{\theta_t} - \phi^\star \|_{2,{\phi^\star}^{-1}}. 
\end{equation}
This further gives
\begin{equation}\label{eq:markov contract}
\| \phi_{t+1} - \phi^\star\|_{2}\leq \lambda_2^{\delta t} \frac{\max{\phi^\star}}{\min{\phi^\star}}\|  \phi_{\theta_t} - \phi^\star \|_{2}. 
\end{equation}
Using \eqref{eq:markov contract} together with assuming $\|\phi_{\theta_t}- \phi_{t}\|_F \leq \tilde \nu$, we have
\begin{equation}\label{eq:contraction}
\| \phi_{\theta_{t+1}} - \phi^\star \|_{2} \leq \lambda_2^{\delta t} \frac{\max{\phi^\star}}{\min{\phi^\star}} \|\phi_{\theta_{t}} - \phi^\star\|_2 + \tilde \nu.
\end{equation}
By applying induction to \eqref{eq:contraction} we have the desired conclusion.
\end{proof}
In Lemma~\ref{lemma:contraction with error}, the error $\tilde \nu$ is caused by approximation and variance error associated with representing the iterates of Alg.~\ref{alg:the_alg} as MPS/TT. To provide estimates to these errors, we define a notion that characterizes the nearness between a function $u$ and an MPS/TT:
\begin{defn}\label{def:identifiable}
For $u:[n]^d\rightarrow \mathbb{R}$, we say it is $\epsilon$-identifiable with a rank-$r$ MPS/TT $u^\star$ if each of the associated $A_k[u]$ (defined in \eqref{eq:Adef}) can be approximated by $\|A_k[u]-A_k[u^\star]\|_2\leq \epsilon \|A_k[u^\star]\|_2,\|B_k[u]-B_k[u^\star]\|_F\leq \epsilon \|B_k[u^\star]\|_F$. Furthermore, $\|A_k[u^\star]^\dagger\|_2\|A_k[u^\star]\|_2\leq 1$.
\end{defn}
This notion of nearness between a function and an MPS/TT can be turned into a $\|\cdot \|_F$ bound between them.

\begin{lemma}\label{lemma:stability}
Let $u^\star$ be a rank-$r$ MPS/TT, represented by cores $\ccG_1^\star,\cdots,\ccG_d^\star$. Let $\hat u$ be $\epsilon$-identifiable with $u^\star$ as defined in Def.~\ref{def:identifiable}, then solving \eqref{eq:both_sketch_CDE}  with $A_k[\hat u], B_k[\hat u]$ gives
\begin{equation}
\|\ccG_k - \ccG^\star_k\|_F \leq \| \ccG^\star_k\|_F \nu(\epsilon),\quad \|u - u^\star \|_F \leq  d\nu(\epsilon)(1+\nu(\epsilon))^d \prod_{k=1}^d \| \ccG_k^\star \|_F.
\end{equation}
where $\nu(\epsilon) := \frac{2\epsilon}{1-\epsilon}$
\end{lemma}
\begin{proof}
First, 
\begin{eqnarray*}
&\ &\|\ccG_k - \ccG^\star_k \|_F \cr
&\leq& \frac{\| A_k[u^\star]^\dagger \|_2}{1-\| A_k[u^\star]^\dagger\|_2\|A_k[u^\star]-A_k[\hat u]\|_2 }2\|\ccG^\star_k\|_F\max\{\|A_k[u^\star]^\dagger-A_k[\hat u]^\dagger\|_2  , \| B_k[u^\star]-B_k[\hat u]\|_F/\|\ccG^\star_k\|_F\}\cr 
&\leq& \frac{\| A_k[u^\star]^\dagger \|_2}{1-\| A_k[u^\star]^\dagger\|_2\|A_k[u^\star]-A_k[\hat u]\|_2 }2\|\ccG^\star_k\|_F\max\{\|A_k[u^\star]\|_2,\|B_k[u^\star]\|_F/\|\ccG^\star_k\|_F\}\epsilon\cr 
&\leq&\frac{\| A_k[u^\star]^\dagger \|_2}{1-\| A_k[u^\star]^\dagger\|_2\|A_k[u^\star]-A_k[\hat u]\|_2}2\|\ccG^\star_k\|_F\max\{\|A_k[u^\star]\|_2,\frac{\|A_k[u^\star]\|_2\|\ccG^\star_k\|_F}{\|\ccG^\star_k\|_F}\}\epsilon\cr
&\leq&\frac{2\epsilon}{1-\epsilon}\|\ccG_k^\star\|_F
\end{eqnarray*}
The first inequality is due to Lemma~\ref{lemma:LS stability}. The second inequality is due to the assumption of the lemma and Def~\ref{def:identifiable}. The last inequality is due to the Definition~\ref{def:identifiable} where $\|A_k[u^\star]^\dagger\|_2\|A_k[u^\star]\|_2\leq 1$. Then using a telescoping sum, we have
\begin{eqnarray*}
&\ &\|u -u^\star\|_F \cr
&=& \| \ccG_1 \circ  \ccG_2 \circ \cdots \circ \ccG_d -\ccG^\star_1 \circ  \ccG^\star_2\cdots \circ \ccG^\star_d\|_F\cr 
&=&  \|(\ccG_1  -\ccG^\star_1) \circ  \ccG^\star_2\cdots \circ \ccG^\star_d\|_F + \| \ccG_1\circ( \ccG_2-\ccG_2^\star)\circ \ccG^\star_3\circ \cdots \ccG^\star_d\|_F\cdots + \| \ccG_1\circ \cdots \ccG_{d-1}\circ (\ccG_d -\ccG^\star_d)\|_F\cr
&\leq & \| \ccG_1  -\ccG^\star_1\|_F  \| \ccG^\star_2\cdots \circ \ccG^\star_d\|_F + \| \ccG_1\|_F \|\ccG_2-\ccG_2^\star\|_F  \|\ccG^\star_3\circ \cdots \ccG^\star_d\|_F\cdots + \| \ccG_1\circ \cdots \ccG_{d-1}\|_F \|\ccG_d -\ccG^\star_d\|_F\cr 
&\leq & \nu(\epsilon)  \|\ccG_1^\star\|_F \| \ccG^\star_2\circ \cdots \circ \ccG^\star_d\|_F  + \nu(\epsilon)  \|\ccG_1\|_F \|\ccG_2^\star\|_F \| \ccG^\star_3\circ \cdots \circ \ccG^\star_d\|_F  \cdots + \nu(\epsilon)  \|\ccG_1\circ \cdots \circ \ccG_{d-1}\|_F \|G_d^\star\|_F\cr 
&\leq &  d\nu(\epsilon) (1+\nu(\epsilon))^d \prod_{k=1}^d \| \ccG_k^\star \|_F.
\end{eqnarray*}
\end{proof}


We are now ready to provide a proof of convergence for Alg~\ref{alg:the_alg}:
\begin{theorem}\label{thm:main}
Let the sketch $S_k,T_k$'s in \eqref{eq:system_A} \eqref{eq:system_B} be coming from cluster basis (Def.~\ref{def:cluster}). In Alg.~\ref{alg:the_alg}, we assume for all $t$, $\phi_t$ is $\epsilon_1$-identifiable with $\phi^\star_{\theta_t}=\ccG_{1,t}^\star\circ \cdots \circ \ccG_{d,t}^\star$ that is rank-$r$. Then for $t\in [0,T]$, we have
\begin{equation}\label{eq:error propagation}
\|\phi_{\theta_t} - \phi^\star \|_F \leq a^t \| \phi_{\theta_0} - \phi^\star \|_F + 2d \nu(\epsilon_1+\epsilon_2) (1+\nu(\epsilon_1 + \epsilon_2))^d  c_G
\end{equation}
with probability $1-\delta$ where 
\begin{itemize}
\item $a := \frac{\max{\phi^\star}}{\min{\phi^\star}}\lambda_2^{\delta t}$ for $\delta t$ large enough such that $a<1$.
\item $\nu(\epsilon)=\frac{2\epsilon}{1-\epsilon}$.
\item $\epsilon_2 =O \left(\sqrt{\frac{n^{2c+1}\log\left(T {d\choose 2c + 1} n^{2c+1}\right)}{N}}\right).$
\item $c_G = \max_t \prod_{k=1}^d \| \ccG_{k,t}^\star \|_F.$\
\item $\delta = O(1/T {d\choose 2c + 1} n^{2c+1}).$
\end{itemize}
\end{theorem}
\begin{proof}
Assuming  $\phi_{t}$ is $\epsilon_1$-identifiable with some $\phi_{\theta_{t}}^\star = \ccG^\star_{1,t}\circ \cdots \circ \ccG^\star_{d,t}$ that is of rank-$r$. By Lemma~\ref{lemma:stability}, this gives  $\|\phi_{t}-\phi_{\theta_{t}}^\star\|_2\leq d\nu(\epsilon_1) (1+\nu(\epsilon_1))^d c_G$. Suppose for now, the empirical distribution $\hat \phi_{t}$ of $\phi_{t}$ satisfies  $\| A_k[\hat \phi_{t}]-A_k[\phi_{t}]\|_F,\| B_k[\hat \phi_{t}]-B_k[\phi_{t}]\|_F\leq \epsilon_2,\forall k\in [d]$. Then $\| A_k[\hat \phi_{t}]-A_k[\phi_{\theta_{t}}^\star]\|_F\leq \| A_k[ \phi_{t}]-A_k[\phi_{\theta_{t}}^\star]\|_F +  \| A_k[\hat \phi_{t}]-A_k[\phi_{t}]\|_F \leq \epsilon_1+\epsilon_2$ (and also $\| B_k[\hat \phi_{t}]-B_k[\phi_{\theta_{t}}^\star]\|_F\leq \epsilon_1+\epsilon_2$). This means $\hat \phi_{t}$ is $(\epsilon_1+\epsilon_2)$-identifiable with $\phi_{\theta_t}^\star$. Since $\phi_{\theta_{t}}$ is the MPS/TT obtained from $A_k[\hat \phi_{t}],B_k[\hat \phi_{t}], k\in [d]$,  by Lemma~\ref{lemma:stability}, we have $\|\phi_{\theta_{t}}-\phi_{\theta_{t}}^\star\|_2\leq d \nu(\epsilon_1+\epsilon_2) (1+\nu(\epsilon_1 + \epsilon_2))^d  c_G$. By triangle inequality, $\|\phi_{\theta_{t}}- \phi_{t} \|_2 \leq \|\phi_{\theta_{t}}- \phi_{\theta_{t}}^\star \|_2+\|\phi_{\theta_{t}}^\star- \phi_{t} \|_2 \leq 2d \nu(\epsilon_1+\epsilon_2) (1+\nu(\epsilon_1 + \epsilon_2))^d  c_G$. 

In order to complete the proof, we now give the expression of $\epsilon_2$, which is due to the variance of $A_k[\hat \phi_t]$ and  $B_k[\hat \phi_t]$. We first look at $B_k[\hat \phi_t]$ of the form \eqref{eq:system_B}. Notice that each pair of $S_{k-1}(\cdot,\xi_{k-1}),T_{k+1}(\cdot,\gamma_{k})$ involves at most $n^{2c}$ variables due to the choice of cluster basis. Therefore, the variance of a single entry $B_k[\hat \phi_t][\xi_{k-1},x_k,\gamma_{k}]$ comes from a $(2c+1)$-marginal distribution ${\hat \phi}_{t,\mathcal{S}}(x_{\mathcal{S}})$ of ${\hat \phi}$ for a subset $\mathcal{S}\in [d]$ where $\vert \mathcal{S} \vert = 2c+1$. Each samples in the empirical distribution ${\hat \phi}_{t,\mathcal{S}}$ is a Bernoulli random variable with probability ${\phi}_{t,\mathcal{S}}(x_\mathcal{S})$. Therefore, by Hoeffding's inequality, for a fixed $x_\mathcal{S}$,
\begin{equation}
\text{Pr}(\vert {\hat \phi}_{t,\mathcal{S}}(x_\mathcal{S}) -  {\phi}_{t,\mathcal{S}}(x_\mathcal{S})\vert \geq \epsilon)\leq \exp(-N\epsilon^2).
\end{equation}
Now we want to bound $\vert {\hat \phi}_{t,\mathcal{S}}(x_\mathcal{S}) -  {\phi}_{t,\mathcal{S}}(x_\mathcal{S})\vert$ for every $x_\mathcal{S}\in [n]^{2c+1}$. With a union bound over $n^{2c+1}$ entries of ${\hat \phi}_{t,\mathcal{S}}$, this implies 
\begin{equation}
\frac{\|   {\hat \phi}_{t,\mathcal{S}} -  {\phi}_{t,\mathcal{S}}\|_F}{\|   {\phi}_{t,\mathcal{S}}\|_F} \leq \frac{\|   {\hat \phi}_{t,\mathcal{S}}-  {\phi}_{t,\mathcal{S}}\|_\infty}{\|    {\phi}_{t,\mathcal{S}}\|_F} \leq \sqrt{n^{2c+1}}\epsilon\quad \text{with probability} \ 1- n^{2c+1} \exp(-N\epsilon^2).
\end{equation}
Now since $B_k[\hat \phi_{t}]$ consists of applying an orthogonal change of basis to $\hat \phi_{t,S}$ for an $\mathcal{S}\in [d]$ (due to our choice of cluster basis in Def.~\ref{def:cluster}), we also have 
\begin{equation}\label{eq:B hoeffding}
\frac{\|  B_k[\hat \phi_{t}]  -  B_k[\phi_{t}]\|_F}{\|   B_k[\phi_{t}]\|_F} \leq \sqrt{n^{2c+1}}\epsilon \quad \text{with probability} \ 1- n^{2c+1} \exp(-N\epsilon^2).
\end{equation}
A similar bound on $\frac{\|  A_k[\hat \phi_{t}]  -  A_k[\phi_{t}]\|_F}{\|   A_k[\phi_{t}]\|_F}$ can be obtained likewise. Using \eqref{eq:B hoeffding} and applying a union bound over all subsets $\mathcal{S}\in [d]$ that contributes to the construction of $A_k[\hat \phi_{t}],B_k[\hat \phi_{t}] $, and for all time $t$, we get
\begin{equation}
\frac{\|  B_k[\hat \phi_{t}]  -  B_k[\phi_{t}]\|_F}{\|   B_k[\phi_{t}]\|_F}, \frac{\|  A_k[\hat \phi_{t}]  -  A_k[\phi_{t}]\|_F}{\| A_k[\phi_{t}]\|_F}\leq \sqrt{n^{2c+1}}\epsilon, \quad k\in[d], t\in [T]
\end{equation}
with probability $1- T {d\choose 2c + 1} n^{2c+1} \exp(-N\epsilon^2)$. Letting $\epsilon = O \sqrt{\frac{\log\left(T {d\choose 2c + 1} n^{2c+1}\right)}{N}}$ and identifying $\epsilon_2 = \sqrt{n^{2c+1}}\epsilon$ completes the proof.
\end{proof}
\begin{rem}
A few remarks are in order. $\epsilon_1$ is the bias error committed by approximating the ``true'' solution $\phi_t$ by an MPS/TT (in terms of the notion defined in Def~\ref{def:identifiable}), which depends on the underlying physics of the problem. $\epsilon_2$ is the variance error of determining an MPS/TT from empirical distribution $\hat \phi_t$. From \eqref{eq:error propagation}, it seems like we have surmounted the curse of dimensionality for solving a (discretized) high-dimensional Fokker-Planck equation, since the error of determining the true solution only grows linearly in $d$ for $\epsilon_1+\epsilon_2=o(1/d)$. However, while the variance error $\epsilon_2$ can be reduced by increasing the number of samples (and we only need samples $N\sim O(n^{4c+2})$ to have a good approximation), in practice when solving a high-dimensional PDE, the bias (approximation) error $\epsilon_1$ could be difficult to reduce. This is where the curse of dimensionality could enter.
\end{rem}

\section{Numerical Experiments}
\label{sec:experiments}

In this section, we present numerical experiments for the two applications introduced in Section~\ref{sec:applications}, namely the quantum many-body ground-state problem (Sections~\ref{sec:quantum_ground_energy}) and the solution of the Fokker–Planck equation (Section~\ref{sec:fokker-planck-simulations}).

\subsection{Stochastic quantum imaginary-time evolution}
\label{sec:quantum_ground_energy}

\begin{figure}[htb]
    \centering
    \includegraphics[width=0.45\textwidth]{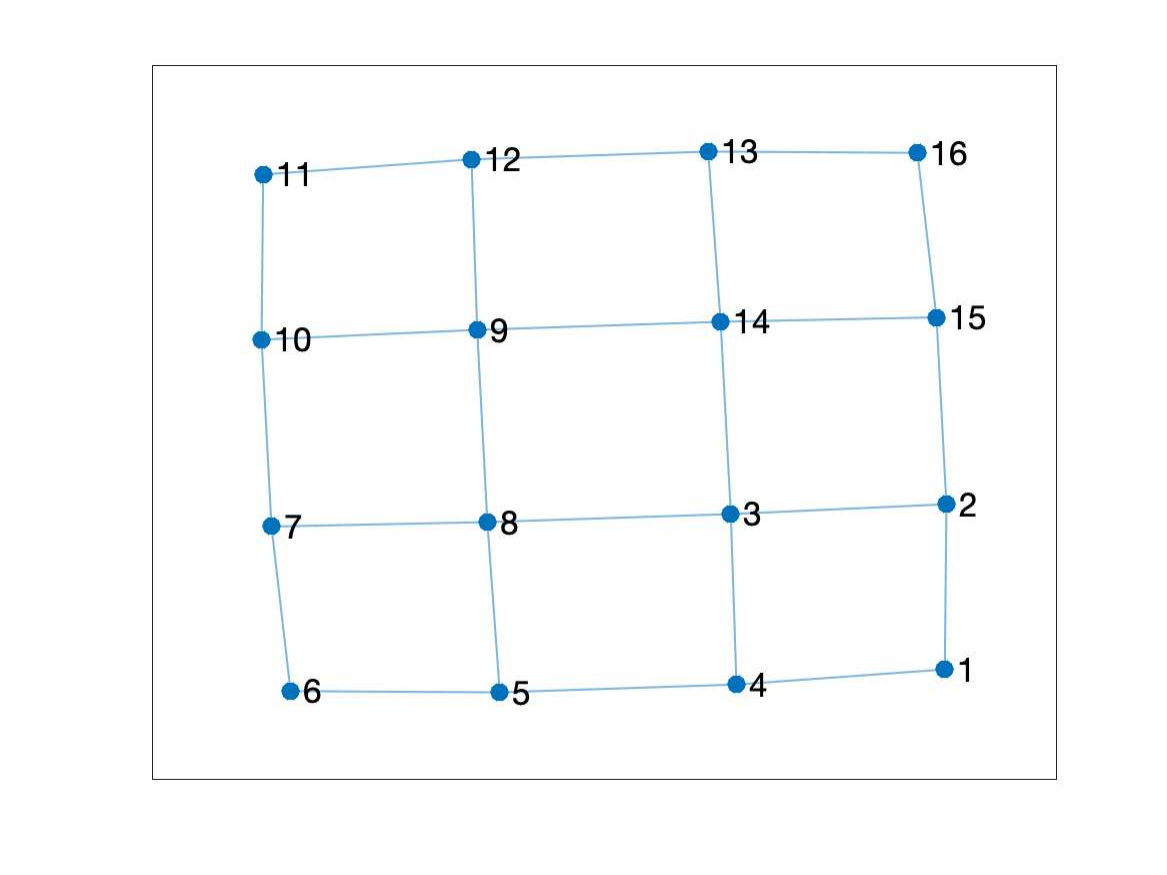}
    \caption{Example of $2$D space-filling curve for Ising model of $4\times 4$ lattices.}
    \label{fig:2D_space_filling}
\end{figure}

In this subsection, we study the ground-state energy estimation problem using the transverse-field Ising model with the following quantum Hamiltonian,
\begin{align}\label{eq:tfi hamiltonian}
    H = -\sum_{i,j=1}^d J_{ij} S^z_i S^z_j - h \sum_{i} S^x_i,
\end{align}
where $S^z_j$, $S^x_j$ are the Pauli matrices \cite{gull1993imaginary},
\begin{align}
    & S^z_j = I_2 \otimes I_2 \otimes \cdots \otimes \underbrace{\begin{pmatrix}1&0\\0&-1\end{pmatrix}}_{\text{j-th dimension}} \otimes \cdots \otimes I_2,\\
    & S^x_j = I_2 \otimes I_2 \otimes \cdots \otimes \underbrace{\begin{pmatrix}0&1\\1&0\end{pmatrix}}_{\text{j-th dimension}} \otimes \cdots \otimes I_2,
\end{align}
and $I_2$ is the $2\times 2$ identity matrix. When $h=1$, and $J$ is the adjacency matrix of a 1D cycle graph, the system undergoes a quantum phase transition. We consider three models under this category: (a) $d=16$ sites, (b) $d=32$ sites, and (c) $d=64$ sites. We also consider a 2D Ising system which is configured as (d) $d=16$ sites with $h=3$, and $J$ the adjacency matrix of a 2D periodic square lattice. For the 1D model, the dimensions of the MPS/TT are naturally ordered according to the sites on the 1D chain. In 2D Ising model case, we use a space-filling curve \cite{sagan2012space} to order the dimensions. For example, we show the space-filling curve and the ordering of MPS/TT dimensions in a $4\times 4$ lattices in \figref{fig:2D_space_filling}.

\begin{figure}[]
    \centering
    \begin{subfigure}{0.4\textwidth}
        \centering
        \includegraphics[width=\textwidth]{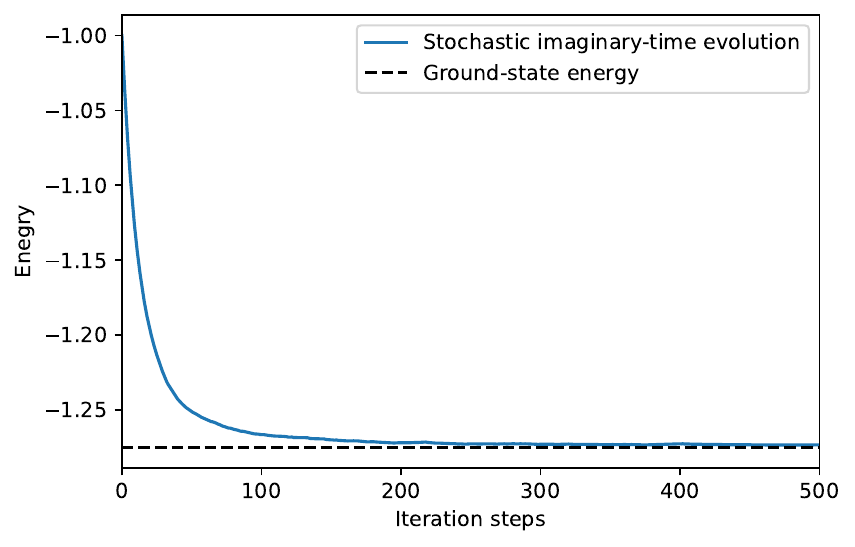}
        \caption{$d=16$ 1D transverse-field Ising symmetric estimator}
        \label{fig:afqmc_1d_16d_energy}
    \end{subfigure}
    \begin{subfigure}{0.4\textwidth}
        \centering
        \includegraphics[width=\textwidth]{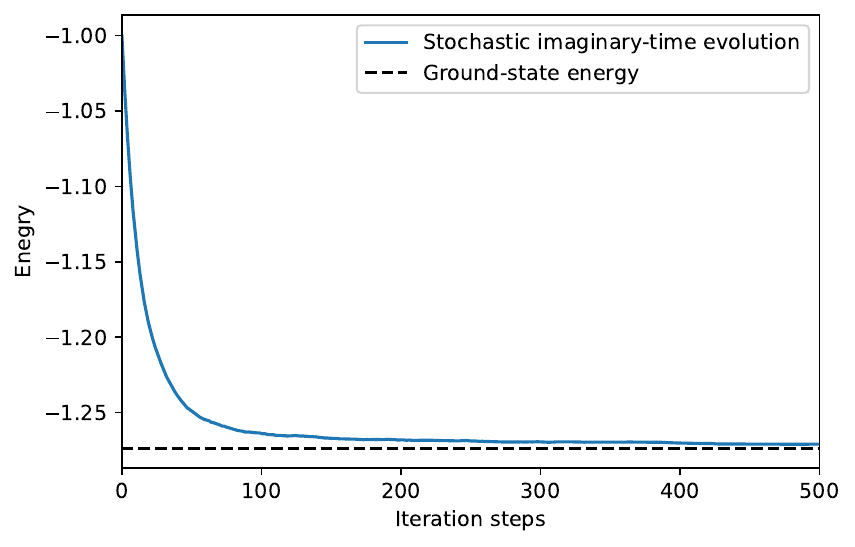}
        \caption{$d=32$ 1D transverse-field Ising symmetric estimator}
        \label{fig:afqmc_1d_32d_energy}
    \end{subfigure}\\
    \begin{subfigure}{0.4\textwidth}
        \centering
        \includegraphics[width=\textwidth]{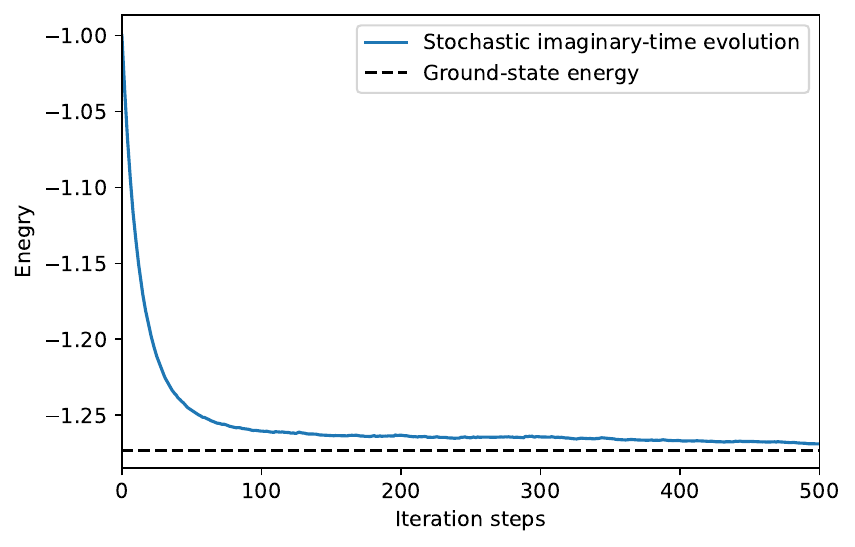}
        \caption{$d=64$ 1D transverse-field Ising symmetric estimator}
        \label{fig:afqmc_1d_64d_energy}
    \end{subfigure}
    \begin{subfigure}{0.4\textwidth}
        \centering
        \includegraphics[width=\textwidth]{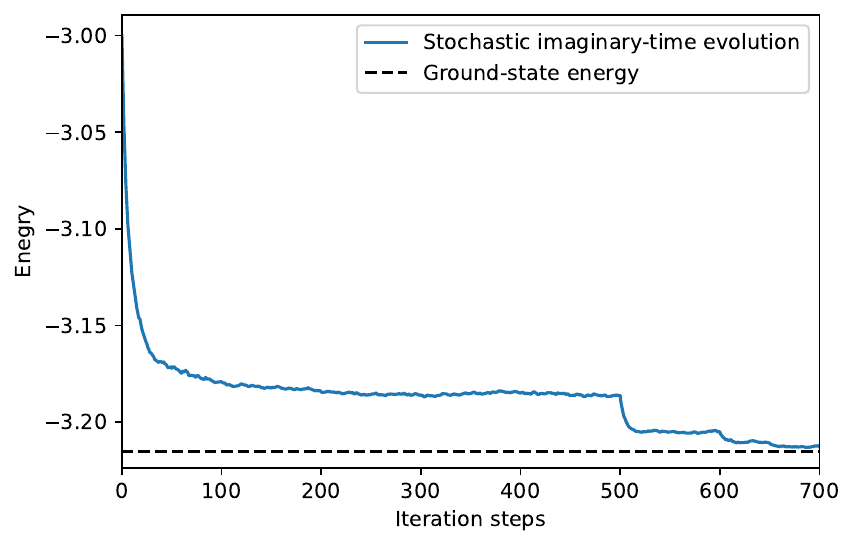}
        \caption{$4\times 4$ 2D transverse-field Ising symmetric estimator}
        \label{fig:afqmc_2d_energy}
    \end{subfigure}
    \caption{Stochastic imaginary-time evolution energy plots. The ground-state energy is shown as horizontal dashed lines.}
    \label{fig:afqmc_enerngy}
\end{figure}

For both stochastic and deterministic imaginary-time evolution, we set the infinitesimal time step $\delta t$ to be $0.01$, and use $2000$ samples in each iteration to approximate the propagator $\exp(-\delta t H)$. The rest of the parameters are set as follows: we use a $\tilde r=60$ (size of $\vert\Gamma_k\vert,\vert\Xi_k\vert$) random tensor sketches (\secref{sec:sketching}) for sketching. The rank of the fitted TT/MPS is gradually increased as the imaginary time grows, up to $16$. 
We initialize the trial wavefunction to be the fully disordered state
\begin{equation}\label{eq:initial trial}
    \Psi_\mathrm{tr} = \frac{1}{\sqrt{2^d}} \cdot
    \left(
    1,1
    \right)^T \otimes \cdots \otimes 
    \left(
    1,1
    \right)^T, 
\end{equation}
which can also be interpreted as the eigenvector corresponding to the smallest eigenvalue of $H_1[O]$. 
The imaginary-time evolution energy is shown in \figref{fig:afqmc_enerngy}. Here we use the symmetric energy estimator given by 
\begin{align}
    E_{\text{symmetric}} = \frac{\langle \phi_t, H, \phi_t \rangle}{\langle \phi_t, \phi_t \rangle},
\end{align}
where $\phi_t$ is the wavefunction of the $t$-th iteration. Theoretically, the energy given by the symmetric estimator can only be larger than the ground-state energy. Often in quantum Monte Carlo, mixed estimators
\begin{align}
    E_{\text{mixed}} = \frac{\langle \phi, H, \phi_t \rangle}{\langle \phi, \phi_t \rangle},
\end{align}
where $\phi$ is a fixed reference wavefunction is used to reduce the bias error originating from the variance in $\phi_t$ \cite{zhang2013auxiliary,shi2021some}. The variance can be further reduced by taking the average of the mixed energy estimators over several iterations.

\begin{figure}[]
    \centering
    \begin{subfigure}{0.4\textwidth}
        \centering
        \includegraphics[width=\textwidth]{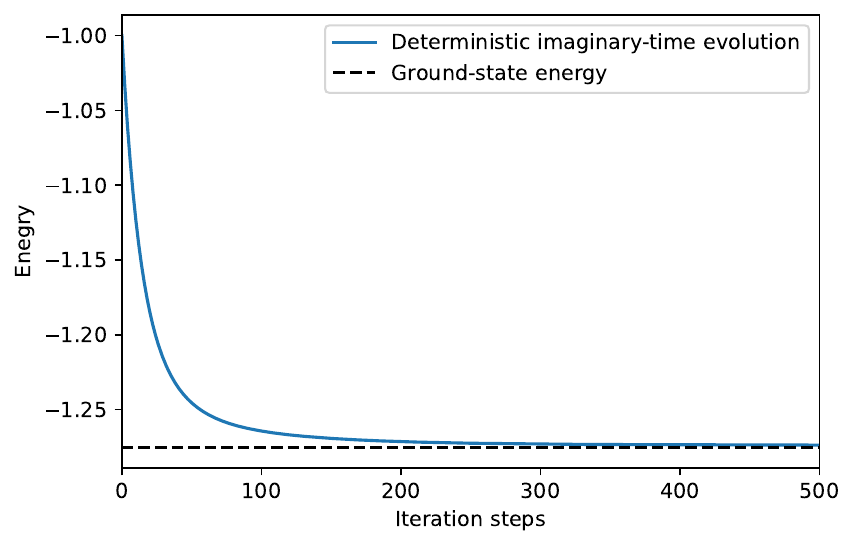}
        \caption{$d=16$ 1D transverse-field Ising symmetric estimator}
        \label{fig:deter afqmc_1d_16d_energy}
    \end{subfigure}
    \begin{subfigure}{0.4\textwidth}
        \centering
        \includegraphics[width=\textwidth]{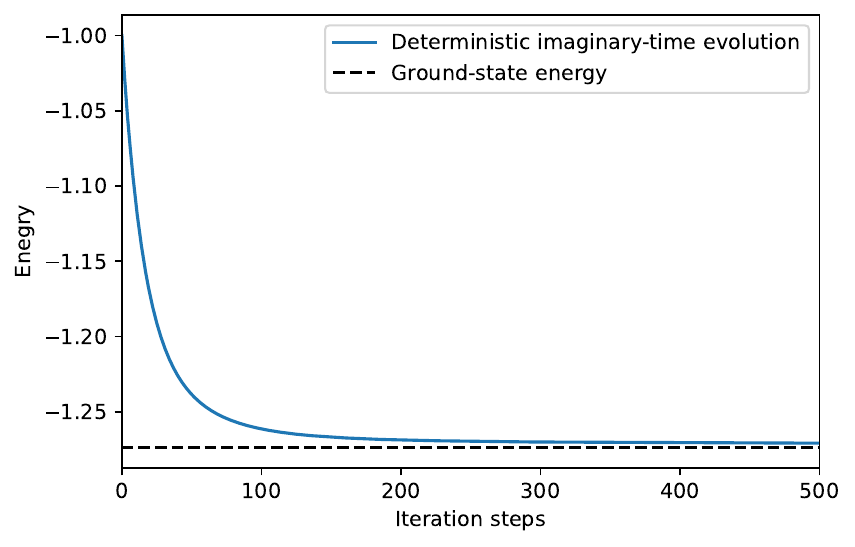}
        \caption{$d=32$ 1D transverse-field Ising symmetric estimator}
        \label{fig:deter afqmc_1d_32d_energy}
    \end{subfigure}\\
    \begin{subfigure}{0.4\textwidth}
        \centering
        \includegraphics[width=\textwidth]{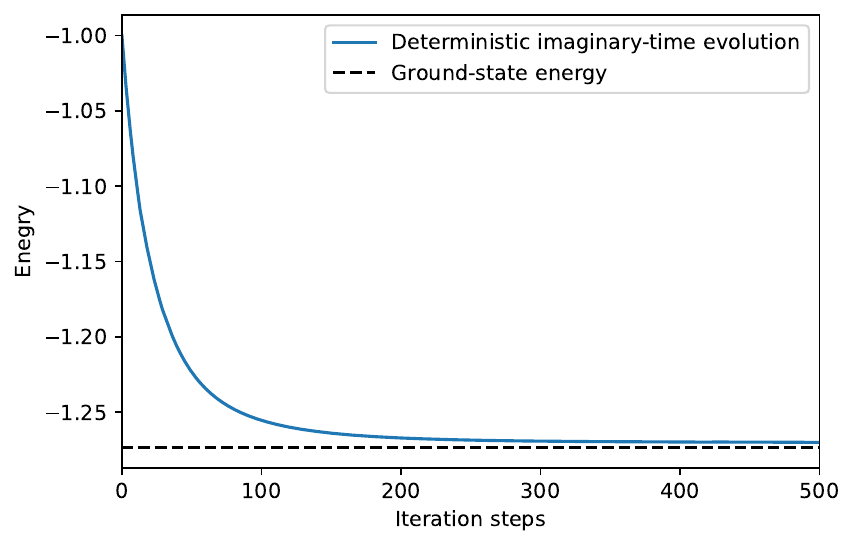}
        \caption{$d=64$ 1D transverse-field Ising symmetric estimator}
        \label{fig:deter afqmc_1d_64d_energy}
    \end{subfigure}
    \begin{subfigure}{0.4\textwidth}
        \centering
        \includegraphics[width=\textwidth]{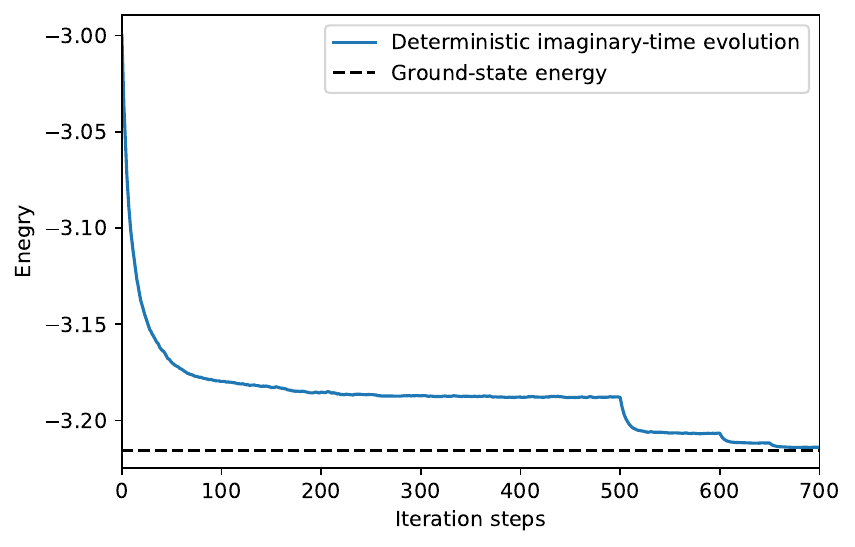}
        \caption{$4\times 4$ 2D transverse-field Ising symmetric estimator}
        \label{fig:deter afqmc_2d_energy}
    \end{subfigure}
    \caption{Deterministic imaginary-time evolution energy plots. The ground-state energy is shown as horizontal dashed lines.}
    \label{fig:deter afqmc_enerngy}
\end{figure}

\figref{fig:afqmc_enerngy} and \figref{fig:deter afqmc_enerngy} illustrate the energy convergence based on the symmetric estimator for 1D Ising model with $16$, $32$, and $64$ sites at $h = 1.0$, and for 2D Ising model with a $4 \times 4$ lattice at $h = 3.0$. The results correspond to stochastic and deterministic quantum imaginary-time evolution, respectively. 
The ground-state energy of 1D Ising model is reported in \cite{sandvik2007variational}. For the $4\times 4$ 2D Ising model, we are still able to store the Hamiltonian exactly in memory so we solve the ground-state energy by exact eigen-decomposition. 

We firstly discuss the numerical results for stochastic quantum imaginary-time evolution. For 1D Ising model with $16$ sites, the ground-state energy per spin is $-1.2753$ and our approach converges to energy $-1.2751$, with a relative error of $1.33\times 10^{-4}$. For 1D Ising model with $32$ sites, the ground-state energy per spin is $-1.2738$ and our approach converges to energy $-1.2731$, with a relative error of $5.31\times 10^{-4}$. For 1D Ising model with $64$ sites, the ground-state energy per spin is $-1.2734$ and our approach converges to energy $-1.2712$, with a relative error of $1.68\times 10^{-3}$. In the 2D case, the ground-state energy per spin is $-3.2155$ and our approach converges to energy $-3.2132$, with a relative error of $7.21\times 10^{-4}$. Our approach has already achieved stable convergence and very accurate ground-state energy estimation with only the symmetric estimators. 

Next, we fix $J$ and test the performance of our proposed algorithm with a range of magnetic fields $h$. The results for 1D Ising model with $d=32$, $d=64$ sites and 2D Ising model with $d=16$ sites are summarized in Table~\ref{table:ising_1d_32}, Table~\ref{table:ising_1d_64} and Table~\ref{table:ising_2d}, respectively. The ground-state energy for 1D models can be exactly computed following \cites{burkhardt1985finite,lieb1961two} and the energy for the 2D Ising model is computed by exact eigen-decomposition. For all magnetic fields, our algorithm achieves stable convergence (not shown) and relative error in energy of order $\leq O(10^{-3})$. 

For the deterministic quantum imaginary-time evolution, the performance is even better. In the 1D Ising model, as shown in Table~\ref{table:ising_1d_32 deter} and Table~\ref{table:ising_1d_64 deter}, the relative error in energy at the critical point $h = 1.0$ is on the order of $O(10^{-5})$ or less, and is significantly smaller for other magnetic field values. In the 2D Ising model with a $4 \times 4$ lattice, the relative energy error is on the order of $O(10^{-4})$, comparable to the results obtained with stochastic evolution.

\begin{table}[]
\centering

\begin{subtable}[t]{\textwidth}
\centering
\begin{tabular}{|c|c|c|c|c|c|}
\hline
                             & $h=0.2$              & $h=0.6$              & $h=1.0$              & $h=1.4$              & $h=1.8$              \\ \hline
Ground-state energy per site & -1.0100            & -1.0922            & -1.2738            & -1.5852            & -1.9418            \\ \hline
Stochastic evolution         & -1.0100            & -1.0922            & -1.2731            & -1.5844            & -1.9411            \\ \hline
Relative error               & $4.43\times 10^{-6}$ & $4.85\times 10^{-5}$ & $5.31\times 10^{-4}$ & $4.78\times 10^{-4}$ & $3.77\times 10^{-4}$ \\ \hline
\end{tabular}
\caption{1D periodic Ising model with $d=32$ sites.}
\label{table:ising_1d_32}
\end{subtable}

\vspace{1em}

\begin{subtable}[t]{\textwidth}
\centering
\begin{tabular}{|c|c|c|c|c|c|}
\hline
                             & $h=0.2$              & $h=0.6$              & $h=1.0$              & $h=1.4$              & $h=1.8$              \\ \hline
Ground-state energy per site & -1.0100            & -1.0922            & -1.2734            & -1.5852            & -1.9418            \\ \hline
Stochastic evolution         & -1.0100            & -1.0922            & -1.2712            & -1.5806            & -1.9382            \\ \hline
Relative error               & $3.50\times 10^{-6}$ & $4.47\times 10^{-5}$ & $1.68\times 10^{-3}$ & $2.88\times 10^{-3}$ & $1.88\times 10^{-3}$ \\ \hline
\end{tabular}
\caption{1D periodic Ising model with $d=64$ sites.}
\label{table:ising_1d_64}
\end{subtable}

\vspace{1em}

\begin{subtable}[t]{\textwidth}
\centering
\begin{tabular}{|c|c|c|c|c|c|}
\hline
                             & $h=2.2$              & $h=2.6$              & $h=3.0$              & $h=3.4$              & $h=3.8$              \\ \hline
Ground-state energy per site & -2.6238            & -2.8990            & -3.2155            & -3.5744            & -3.9487            \\ \hline
Stochastic evolution         & -2.6227            & -2.8885           & -3.2132            & -3.5723            & -3.9468            \\ \hline
Relative error               & $4.44\times 10^{-4}$ & $8.64\times 10^{-4}$ & $7.21\times 10^{-4}$ & $5.92\times 10^{-4}$ & $4.92\times 10^{-4}$ \\ \hline
\end{tabular}
\caption{2D Ising model on a $4\times4$ lattice ($d=16$). The ground truth is computed by exact diagonalization.}
\label{table:ising_2d}
\end{subtable}

\caption{Estimated ground-state energy per site for stochastic quantum imaginary-time evolution across different Ising model configurations. The maximal TT rank is $16$. }
\label{table:ising_all}
\end{table}

\begin{table}[H]
\centering

\begin{subtable}[t]{\textwidth}
\centering
\begin{tabular}{|c|c|c|c|c|c|}
\hline
                             & $h=0.2$              & $h=0.6$              & $h=1.0$              & $h=1.4$              & $h=1.8$              \\ \hline
Ground-state energy per site & -1.0100              & -1.0922              & -1.2738              & -1.5852              & -1.9418              \\ \hline
Deterministic evolution      & -1.0100              & -1.0922              & -1.2738              & -1.5852              & -1.9418              \\ \hline
Relative error               & $6.82\times 10^{-15}$ & $2.60\times 10^{-9}$ & $2.70\times 10^{-6}$ & $2.63\times 10^{-8}$ & $9.42\times 10^{-10}$ \\ \hline
\end{tabular}
\caption{1D periodic Ising model with $d=32$ sites.}
\label{table:ising_1d_32 deter}
\end{subtable}

\vspace{1em}

\begin{subtable}[t]{\textwidth}
\centering
\begin{tabular}{|c|c|c|c|c|c|}
\hline
                             & $h=0.2$              & $h=0.6$              & $h=1.0$              & $h=1.4$              & $h=1.8$              \\ \hline
Ground-state energy per site & -1.0100              & -1.0922              & -1.2734              & -1.5852              & -1.9418              \\ \hline
Deterministic evolution      & -1.0100              & -1.0922              & -1.2733              & -1.5852              & -1.9418              \\ \hline
Relative error               & $2.26\times 10^{-14}$ & $4.50\times 10^{-9}$ & $2.11\times 10^{-5}$ & $3.08\times 10^{-8}$ & $9.10\times 10^{-10}$ \\ \hline
\end{tabular}
\caption{1D periodic Ising model with $d=64$ sites.}
\label{table:ising_1d_64 deter}
\end{subtable}

\vspace{1em}

\begin{subtable}[t]{\textwidth}
\centering
\begin{tabular}{|c|c|c|c|c|c|}
\hline
                             & $h=2.2$              & $h=2.6$              & $h=3.0$              & $h=3.4$              & $h=3.8$              \\ \hline
Ground-state energy per site & -2.6238              & -2.8990              & -3.2155              & -3.5744              & -3.9487              \\ \hline
Deterministic evolution      & -2.6230              & -2.8893              & -3.2140              & -3.5731              & -3.9477              \\ \hline
Relative error               & $3.03\times 10^{-4}$ & $5.63\times 10^{-4}$ & $4.66\times 10^{-4}$ & $3.31\times 10^{-4}$ & $2.28\times 10^{-4}$ \\ \hline
\end{tabular}
\caption{2D Ising model on a $4\times4$ lattice ($d=16$). The ground truth is computed by exact diagonalization.}
\label{table:ising_2d deter}
\end{subtable}

\caption{Estimated ground-state energy per site for deterministic quantum imaginary-time evolution across different Ising model configurations. The maximal TT rank is $16$. }
\label{table:ising_all_deter}
\end{table}

\subsection{Fokker-Planck Equation}
\label{sec:fokker-planck-simulations}

In this numerical experiment, we solve the Fokker-Planck equation by parameterizing the density with MPS/TT. We consider two systems in this subsection, one with a simple double-well potential that takes a separable form, which is intrinsically a 1D potential, and the other one being the Ginzburg-Landau potential, where we only compare the obtained marginals with the ground truth marginals since the true density is exponentially sized.

\subsubsection{Double-well Potential}

\begin{figure}[htb]
    \centering
    \begin{subfigure}{0.50\textwidth}
        \centering
        \includegraphics[width=\textwidth]{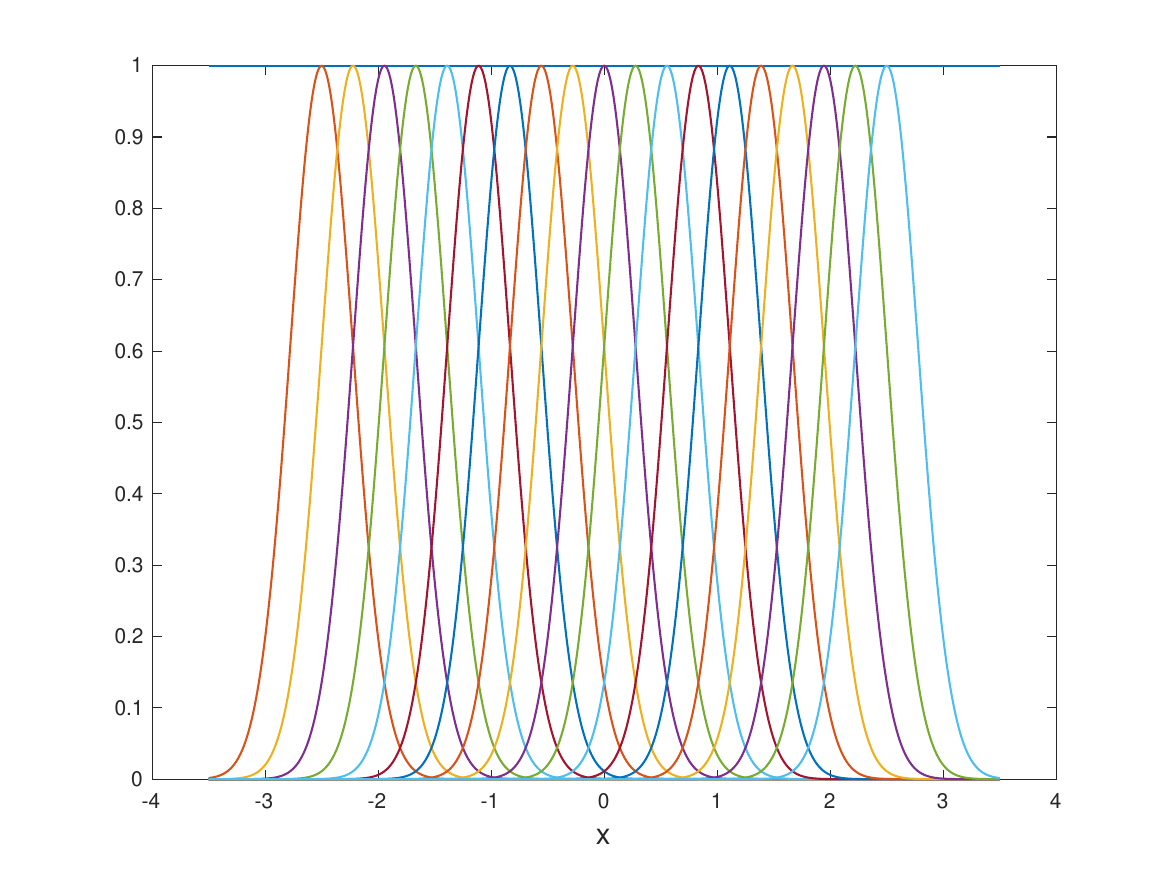}
    \end{subfigure}
    \caption{Visualization of univariate basis functions for each dimension. Here we use univariate Gaussian kernel functions as our basis functions.}
    \label{fig:DW_basis}
\end{figure}

We consider the following double-well potential
\begin{align}
    V(x) = (x_1^2 - 1)^2 + 0.3\sum_{j=2}^d x_j^2, 
    \label{eq:DW-potential2}
\end{align}
and the particle dynamics governed by overdamped Langevin equation \eqref{eq:Langevin}. Since the potential function is easily separable, the equilibrium Boltzmann density is a product of univariate densities for each dimension, i.e.
\begin{align}
    \frac{1}{Z_{\beta}} \exp (-\beta V(\vx)) = \frac{1}{Z_{\beta}} \exp \left(-\beta (x_1^2 - 1)^2 \right) \prod_{j=2}^d \exp \left(-0.3\beta x_j^2 \right).
\end{align}
In this example, we use $\beta=1$ and $d=10$. The support of the domain is a hypercube $[-M, M]^d$ where $M=2.5$. To obtain a continuous MPS/TT approximation, we use the Gaussian kernel function as univariate basis functions $\{b_l\}_{l=1}^{20}$, where
\begin{align}
    b_l(\cdot) = \exp \left( -\frac{(\ \cdot \ + M - (l-1)\Delta x)^2}{2\Delta x^2} \right), \ l=1,\dots,20,
\end{align}
$\Delta x=5/18$ to form cluster basis for spanning the PDE solution, as mentioned in \secref{sec:sketching}. In terms of the sketching function, we use cluster basis functions with $c=1,2$ as mentioned in \ref{def:cluster}. This results in $\tilde r=100$ tensor sketches. After sketching, we obtain a continuous analog of \eqref{eq:both_sketch_CDE} (where each $\ccG_k$ is a set of continuous univariate functions), and we solve this linear least-squares in function space using the same set of univariate basis functions. We visualize all the basis functions in \figref{fig:DW_basis}.

We start from the uniform distribution over the hypercube $[-M, M]^d$ and evolve the distribution towards equilibrium. To approximate the given solution at each time $\phi_{\theta_t}$ as an MPS/TT, we first sample from this distribution via a conditional sampling \secref{sec:conditional_sampling}. Then we simulate the overdamped Langevin process forward for the sampled particles up to time $\delta t$, as detailed in \secref{sec:fokker-planck-evolution}. Then we estimate a new  MPS/TT representation $\phi_{\theta_{t+1}}$. We choose $\delta t = 0.02$ time and we use $N=10^4$ samples for all iterations. 

\begin{figure}[htb]
    \centering
    \begin{subfigure}{0.40\textwidth}
        \centering
        \includegraphics[width=\textwidth]{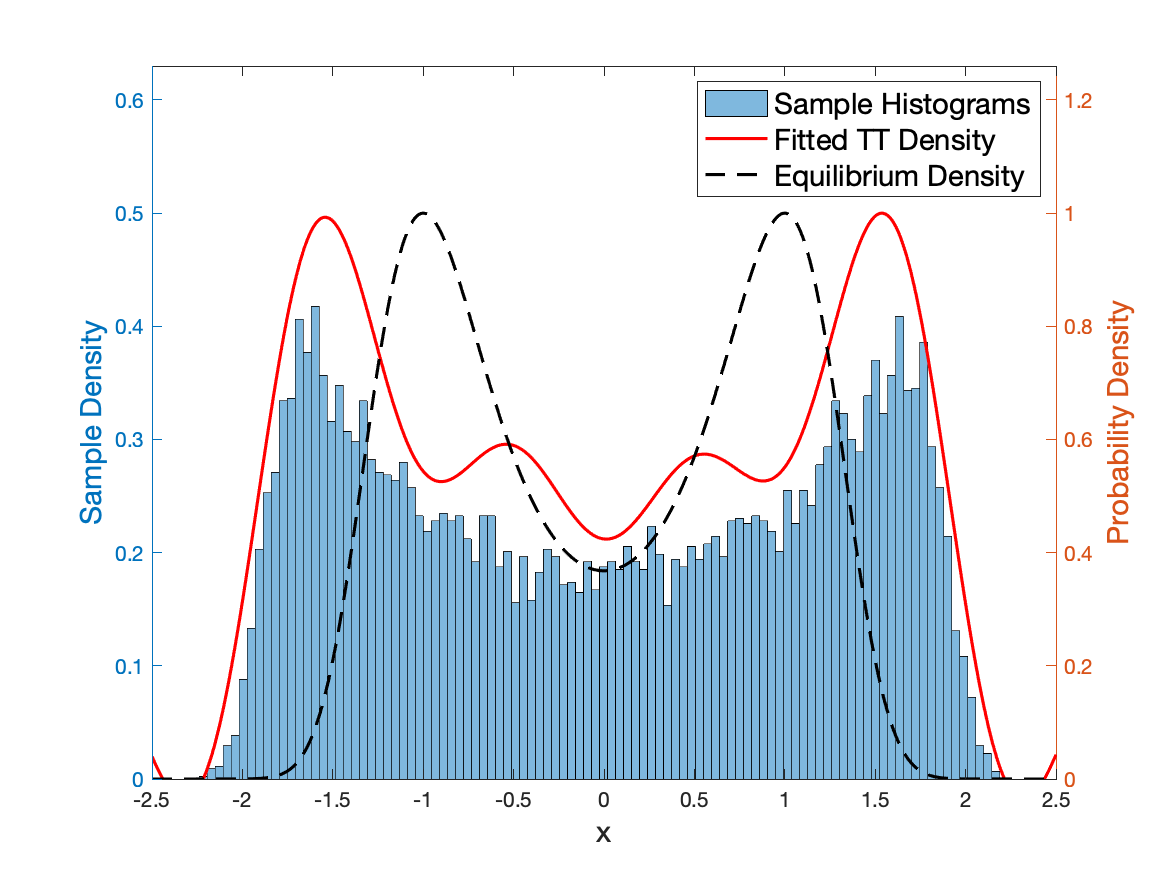}
        \caption{Iteration 1}
        \label{fig:DW_step1}
    \end{subfigure}\quad \quad
    \begin{subfigure}{0.40\textwidth}
        \centering
        \includegraphics[width=\textwidth]{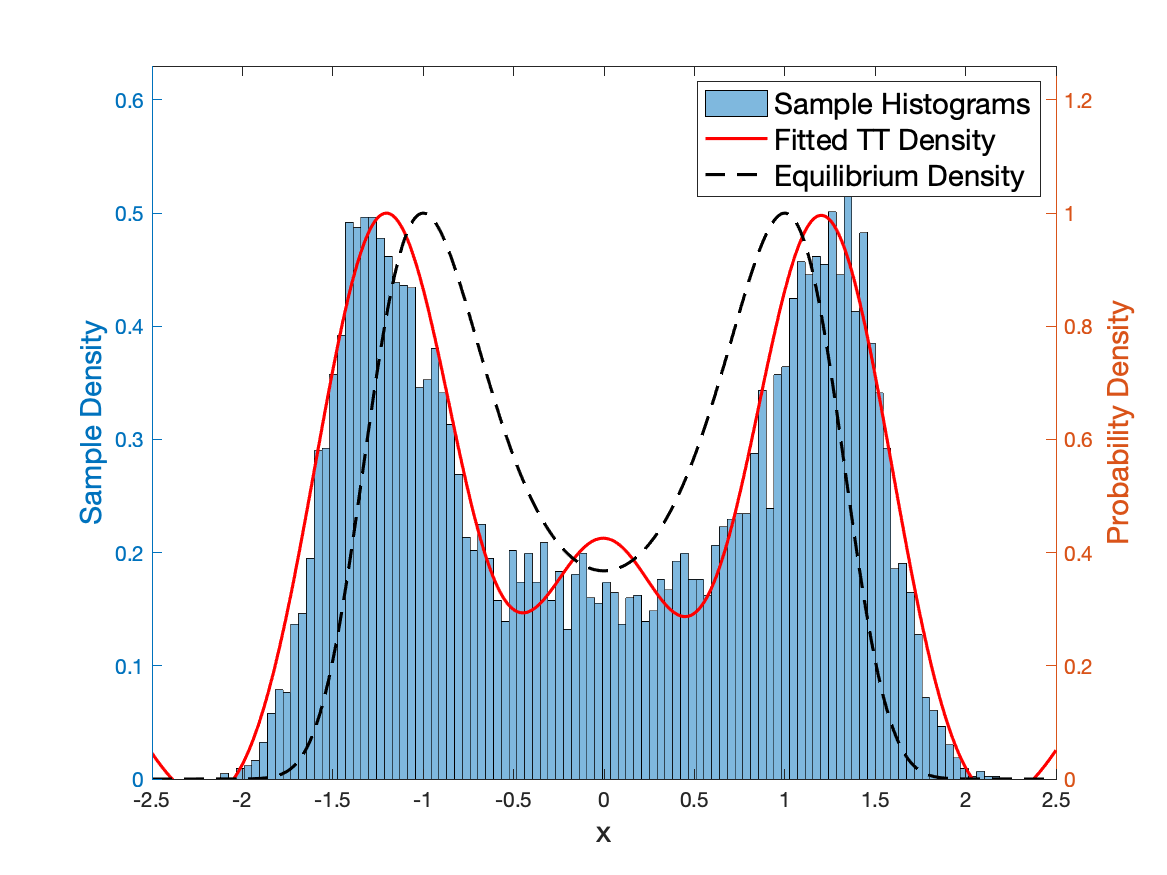}
        \caption{Iteration 3}
        \label{fig:DW_step2}
    \end{subfigure}\\
    \begin{subfigure}{0.40\textwidth}
        \centering
        \includegraphics[width=\textwidth]{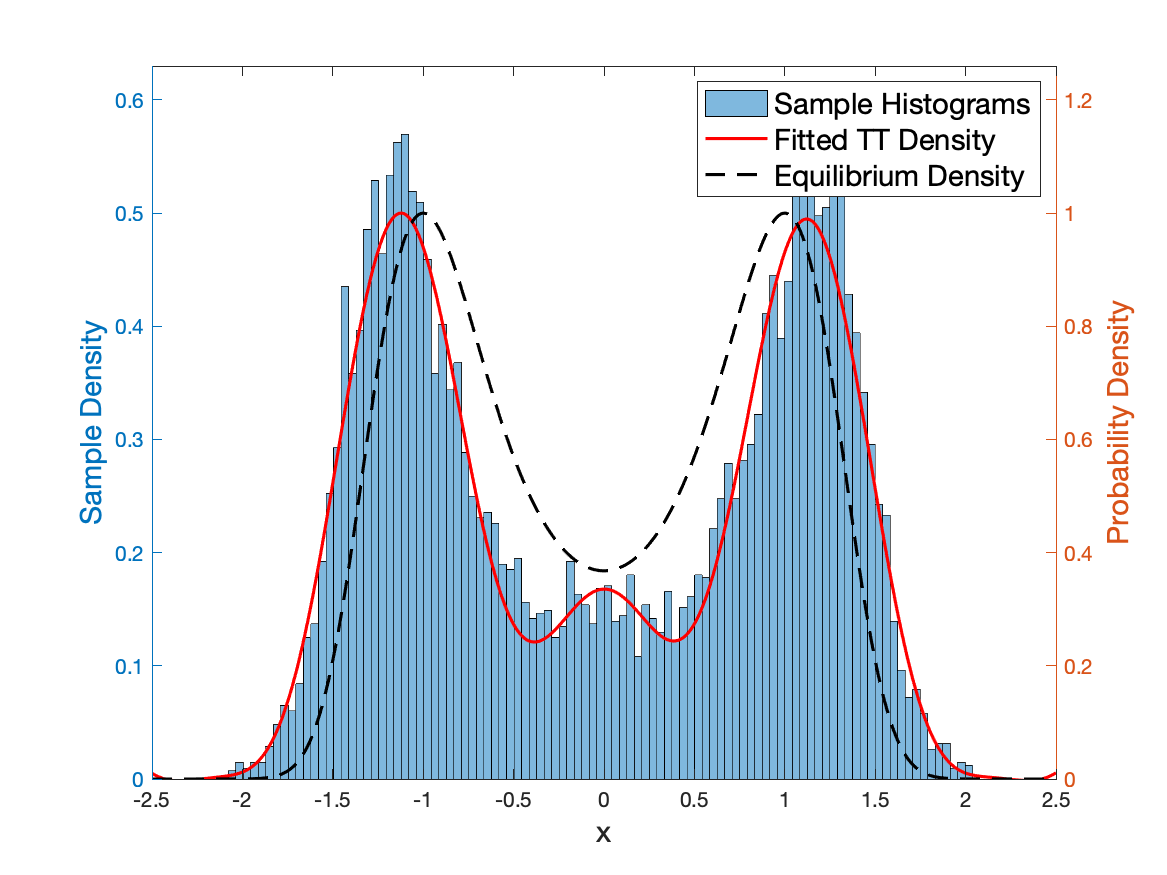}
        \caption{Iteration 5}
        \label{fig:DW_step3}
    \end{subfigure}\quad \quad
    \begin{subfigure}{0.40\textwidth}
        \centering
        \includegraphics[width=\textwidth]{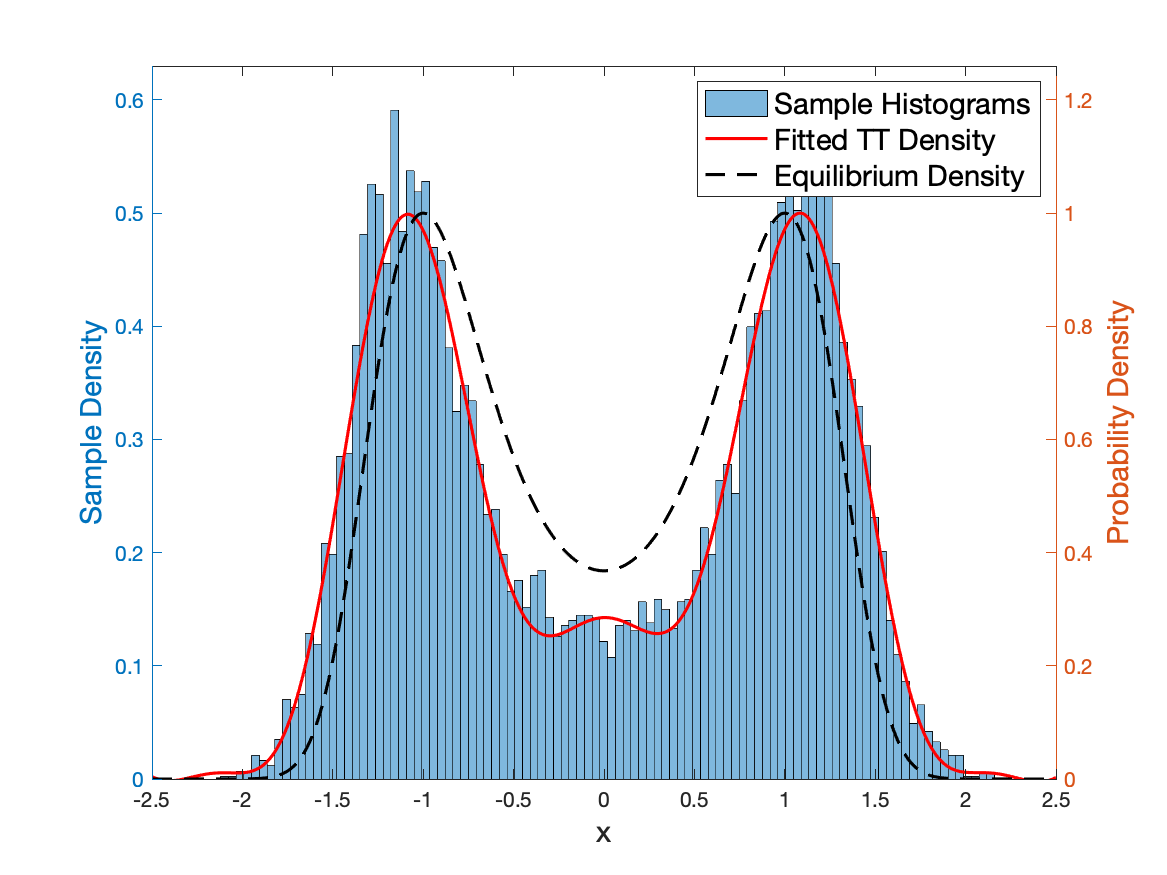}
        \caption{Iteration 7}
        \label{fig:DW_step4}
    \end{subfigure}\\
    \begin{subfigure}{0.40\textwidth}
        \centering
        \includegraphics[width=\textwidth]{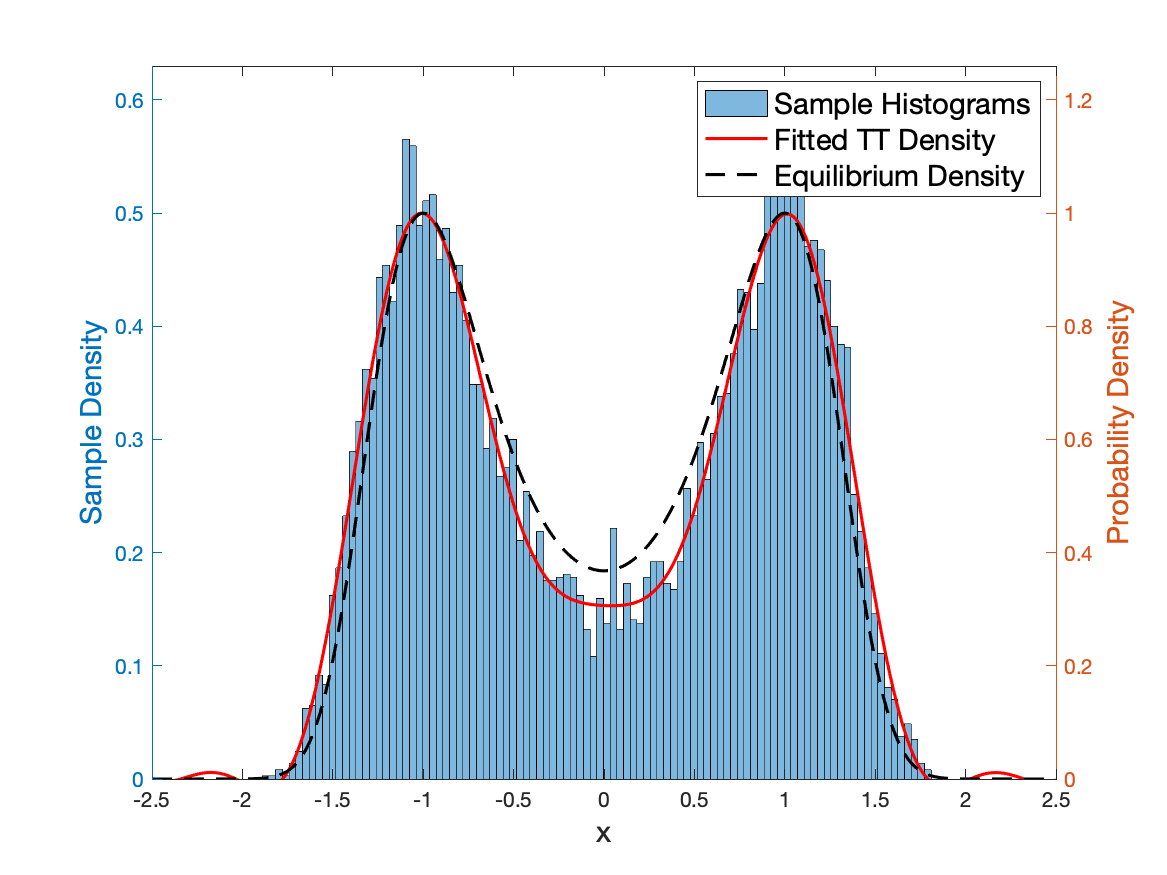}
        \caption{Iteration 20}
        \label{fig:DW_step5}
    \end{subfigure}\quad \quad
    \begin{subfigure}{0.40\textwidth}
        \centering
        \includegraphics[width=\textwidth]{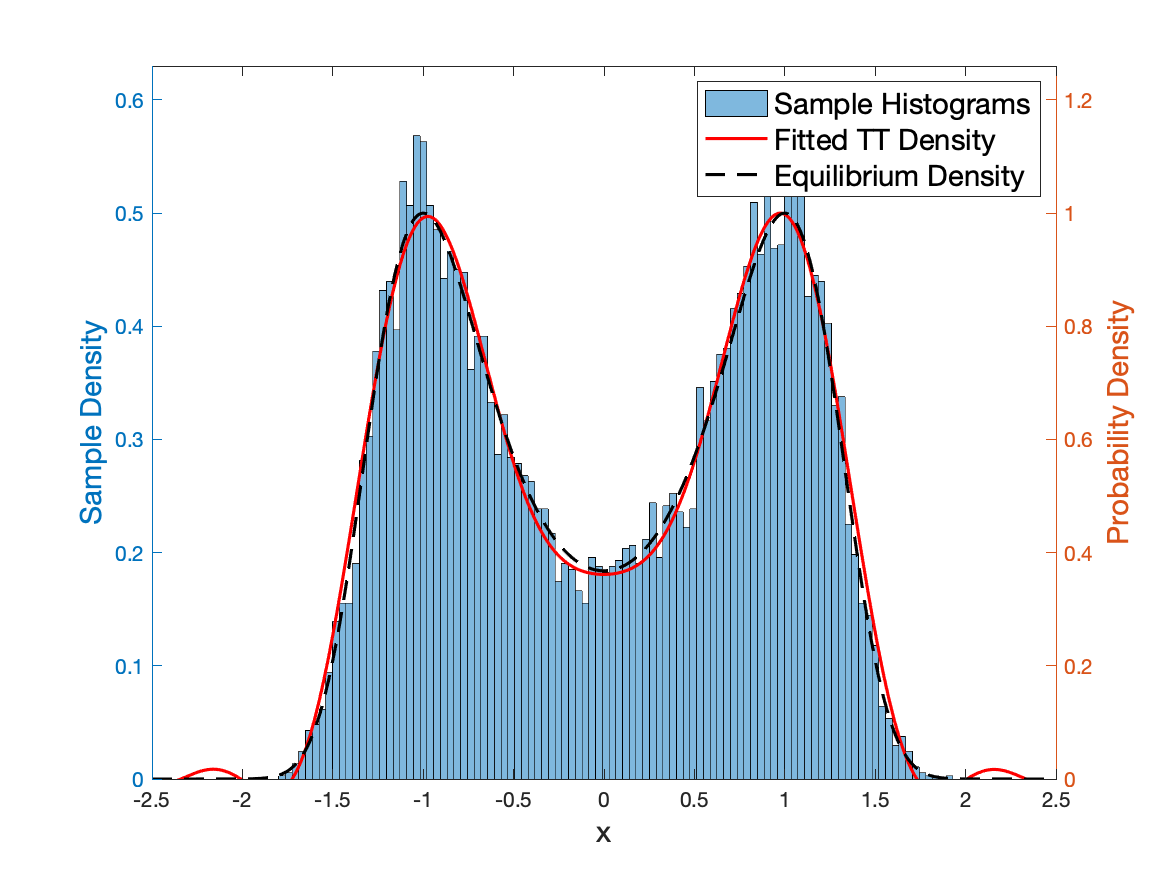}
        \caption{Iteration 30}
        \label{fig:DW_step6}
    \end{subfigure}
    \caption{Visualization of the evolution of first marginal for the double-well potential. The blue histograms correspond to the sample histograms after Langevin simulations at each iteration ($\hat \phi_{t+1}$ in \eqref{eq:ivp}). The estimated continuous MPS/TT density $\hat \phi_{\theta_{t+1}}$ in \eqref{eq:ivp} and the target equilibrium density $\phi^*$ are represented with red solid lines and black dashed lines, respectively. }
    \label{fig:DW_evolution}
\end{figure}

In \figref{fig:DW_evolution} we visualize the simulated Langevin particles, the fitted continuous MPS/TT density and the target equilibrium density for the first dimension at iteration $1,3,5,7,20$ and $30$. We can observe that the particle distribution gets evolved effectively by the Langevin dynamics and the fitted continuous MPS/TT density accurately captures the histograms of the particle samples. The low-complexity continuous MPS/TT format also serves as an extra regularization and as a result, is not prone to overfitting. To quantify the performance of our algorithm, we evaluate the relative error metric $E = \|{\phi_1}^* -{\phi_1}_{\theta_j}\| / \|{\phi_1}^* \|$ where ${\phi_1}^*$ and ${\phi_1}_{\theta_j}$ are the first marginal distribution of the ground truth and MPS/TT represented distribution, respectively. At iteration $30$, the relative error $E=3.8\times 10^{-2}$. We can further improve the performance of the algorithm by choosing more basis functions and generating more stochastic samples.

\subsubsection{1D Ginzburg-Landau Potential}
\label{sec:1D-GL}

The Ginzburg-Landau theory was developed to provide a mathematical description of phase transition \cite{GL-model}. In this numerical example, we consider a simplified Ginzburg-Landau model, in which the potential energy is defined as 
\begin{align}
    V(U) := \sum_{i=1}^{d+1} \frac{\lambda}{2} \left(\frac{U_i - U_{i-1}}{h}\right)^2 + \frac{1}{4\lambda} (1 - U_i^2)^2,
    \label{eq:discrete-GL}
\end{align}
where $h = 1/(d+1)$, $U_0=U_{d+1}=0$. We fix $d=16$, $\lambda=0.03$ and the temperature $\beta=1/8$. We use the same set of $20$ basis function for all dimensions as shown in \figref{fig:DW_basis}. We use the cluster basis with $c=1,2$ as sketching tensor, which results in a tensor rank $\tilde r=240$. For this example, we solve the Fokker-Planck equation starting from the initial uniform distribution over the hypercube $[-M, M]^d$ and evolve the distribution with $\delta t = 0.002$ time and $N=10^4$ samples.

\begin{figure}[htb]
    \centering
    \begin{subfigure}{0.40\textwidth}
        \centering
        \includegraphics[width=\textwidth]{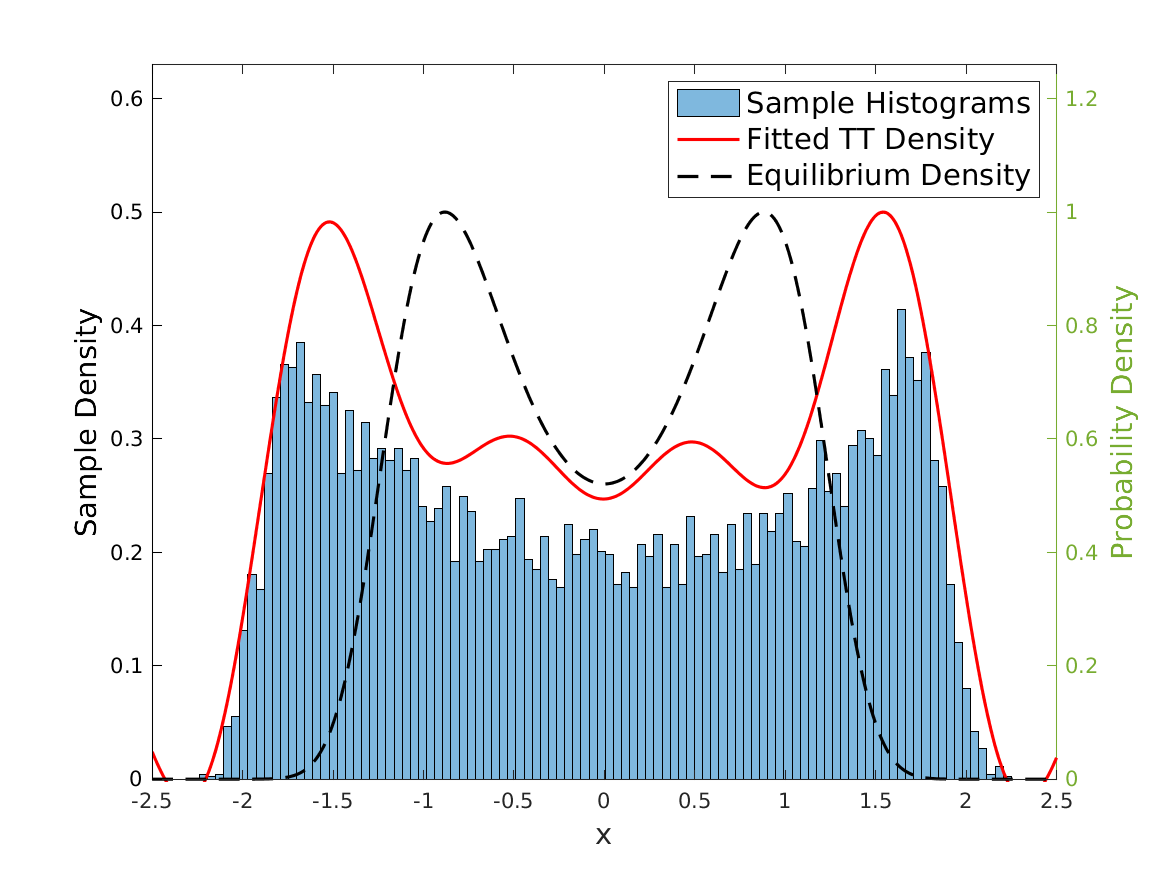}
        \caption{Iteration 1}
        \label{fig:GL_step1}
    \end{subfigure}\quad \quad
    \begin{subfigure}{0.40\textwidth}
        \centering
        \includegraphics[width=\textwidth]{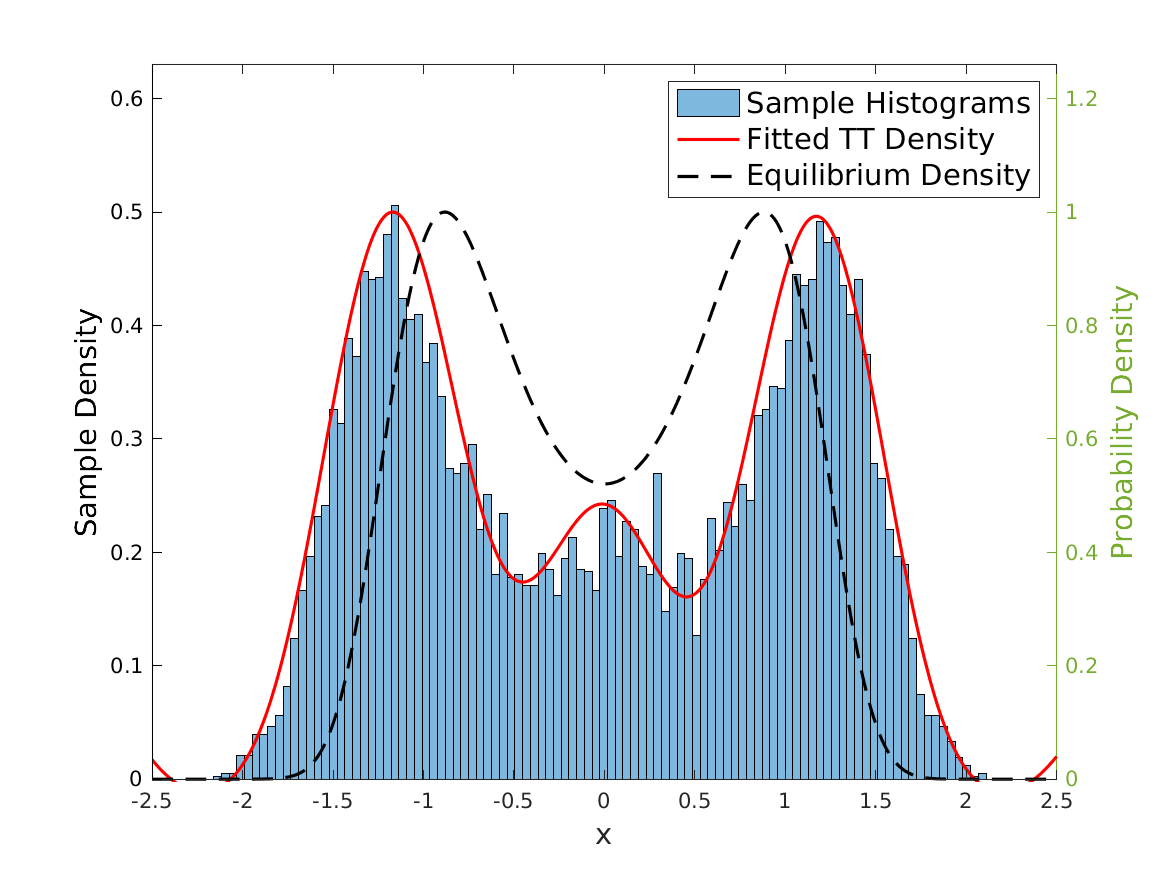}
        \caption{Iteration 3}
        \label{fig:GL_step2}
    \end{subfigure}\\
    \begin{subfigure}{0.40\textwidth}
        \centering
        \includegraphics[width=\textwidth]{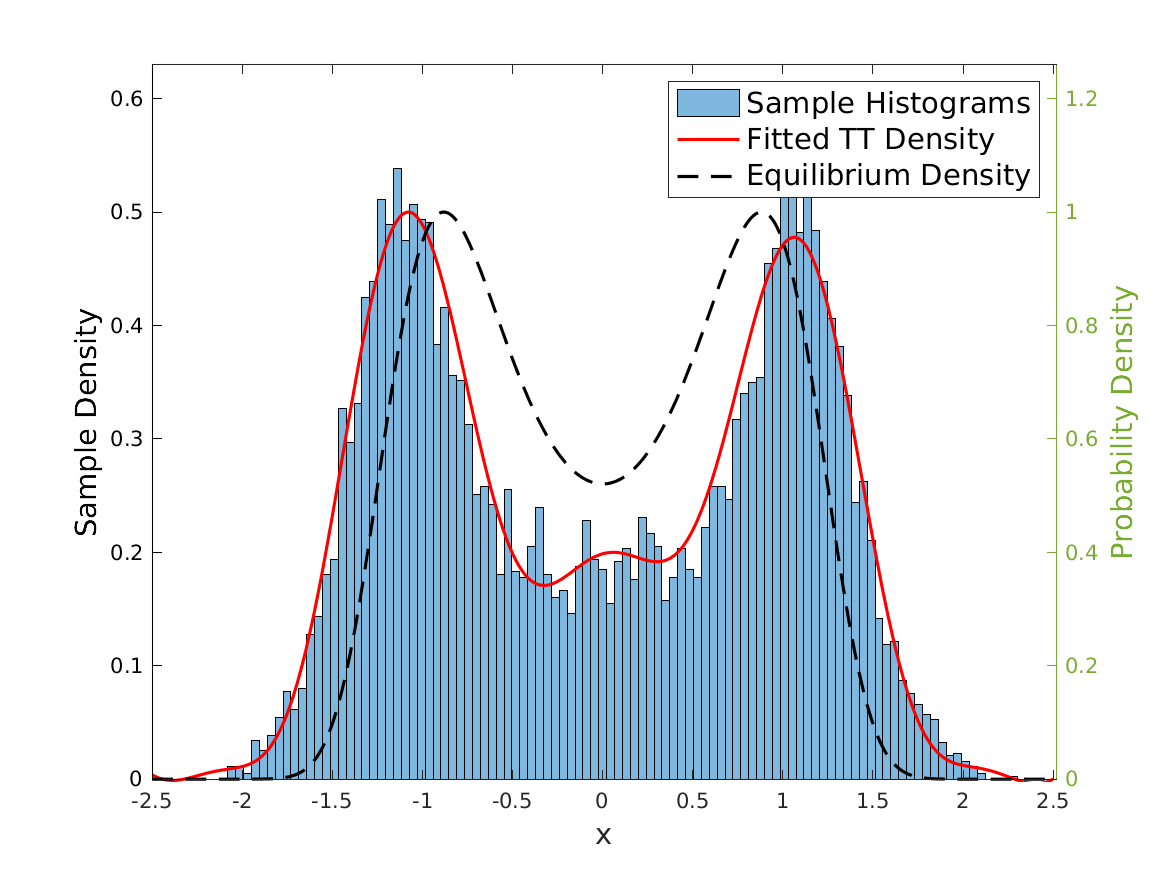}
        \caption{Iteration 5}
        \label{fig:GL_step3}
    \end{subfigure}\quad \quad
    \begin{subfigure}{0.40\textwidth}
        \centering
        \includegraphics[width=\textwidth]{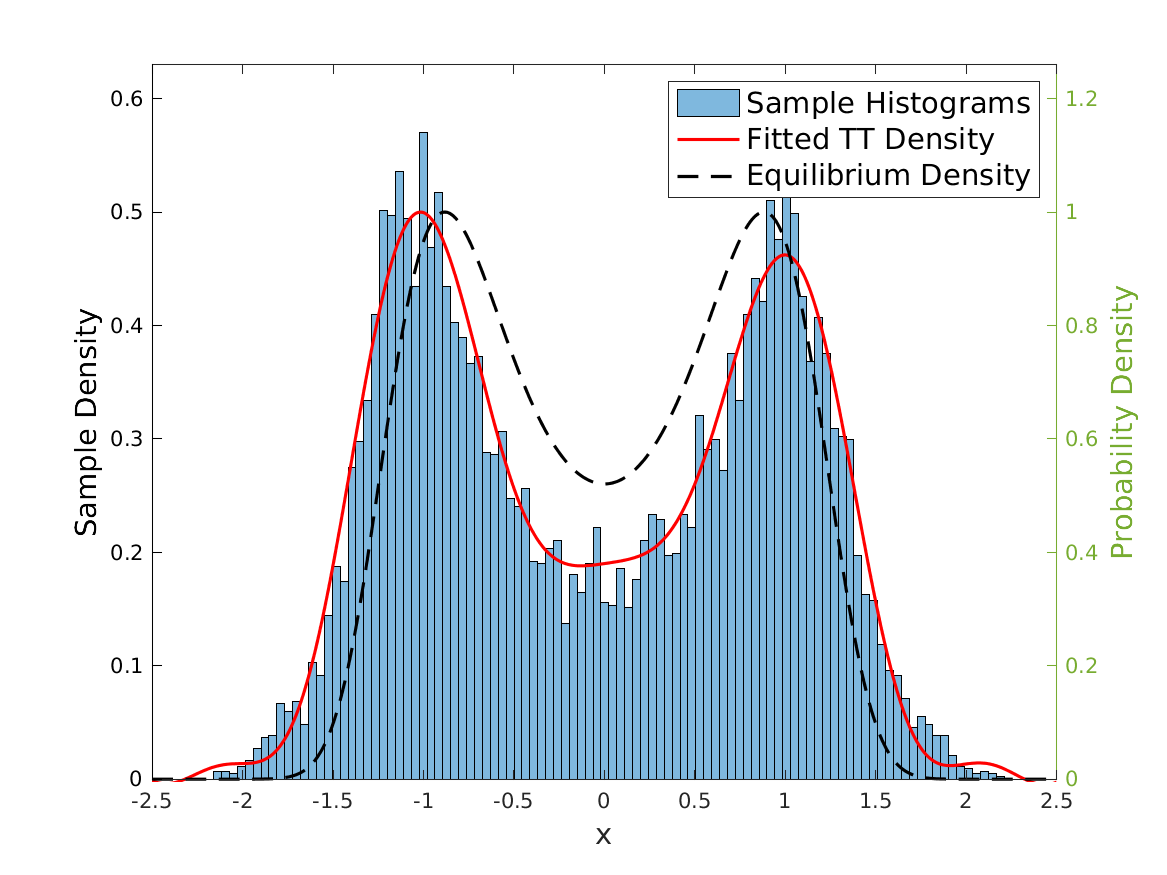}
        \caption{Iteration 7}
        \label{fig:GL_step4}
    \end{subfigure}\\
    \begin{subfigure}{0.40\textwidth}
        \centering
        \includegraphics[width=\textwidth]{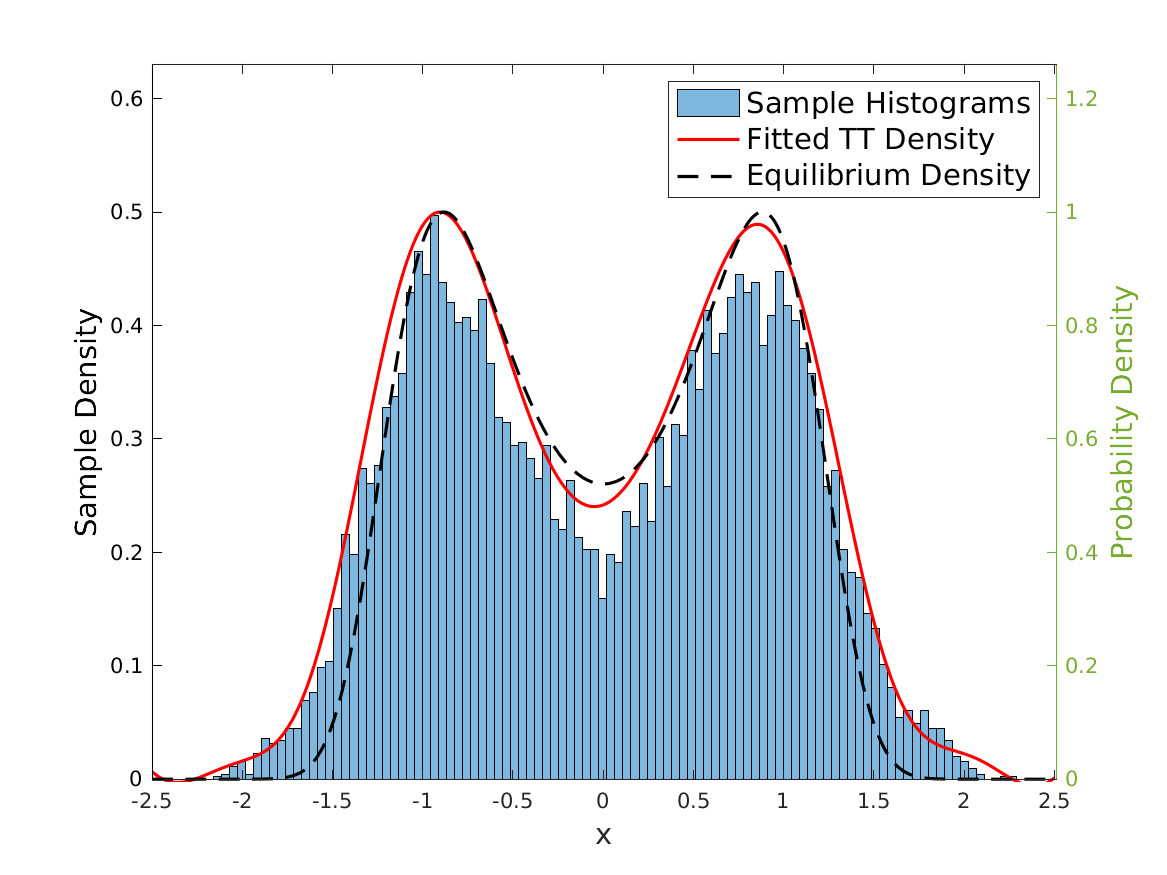}
        \caption{Iteration 20}
        \label{fig:GL_step5}
    \end{subfigure}\quad \quad
    \begin{subfigure}{0.40\textwidth}
        \centering
        \includegraphics[width=\textwidth]{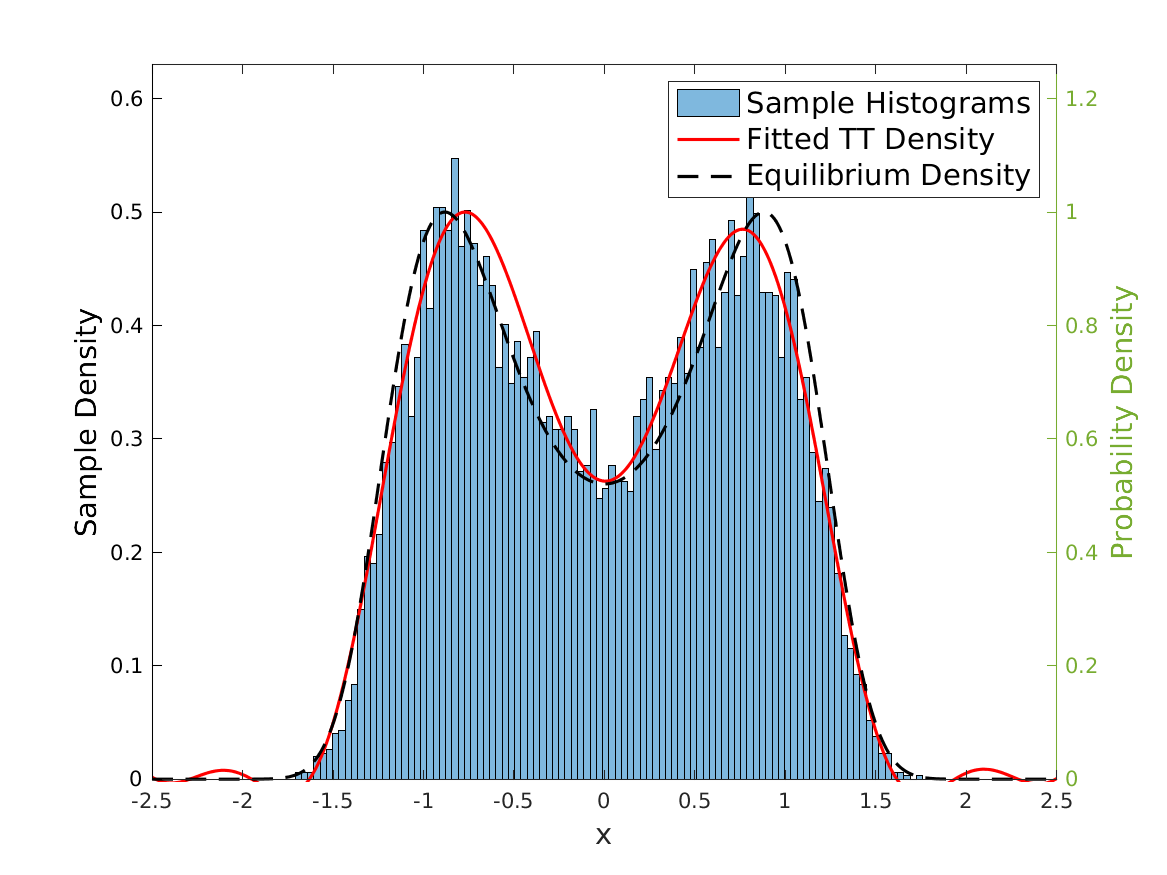}
        \caption{Iteration 100}
        \label{fig:GL_step6}
    \end{subfigure}
    \caption{Visualization of the evolution of the 8-th marginal distribution for the 1D Ginzburg-Landau potential. The blue histograms correspond to the sample histograms after Langevin simulations at each iteration ($\hat \phi_{t+1}$. in \eqref{eq:ivp}) The estimated continuous MPS/TT density $\hat \phi_{\theta_{t+1}}$ in \eqref{eq:ivp} and the target equilibrium density $\phi^*$ are represented with red solid lines and black dashed lines, respectively. }
    \label{fig:GL_evolution}
\end{figure}

In \figref{fig:GL_evolution} we visualize the 8-th marginal distribution of the particle dynamics, the MPS/TT density, and the equilibrium density, at iteration $1,3,5,7,20$ and $100$. At iteration $100$, the relative error of the 8-th marginal distribution is $E=9.8\times 10^{-2}$. 

\subsubsection{2D Ginzburg-Landau Potential}
\label{sec:2D-GL}

In this section, we consider an analogous Ginzburg-Landau-type model on a two-dimensional $(\tilde{d}+2, \tilde{d}+2)$ square lattice space. The 2D Ginzburg-Landau example cannot be easily done with the traditional tensor-network method, which requires the compression of the semigroup operator into an MPO. Such a compression is difficult when the variables cannot be ordered along a 1D line. This, however, is not an issue for us, since we never compress any operator in our framework. Similarly, we define the value of the scalar field at lattice points as $U_{i,j}$, $i,j=0,\dots, \tilde{d}+1$, and define the potential energy,

\begin{align}
    V(U) := \frac{\lambda}{2} \sum_{i=1}^{\tilde{d}+1} \sum_{j=1}^{\tilde{d}+1} \left[\left(\frac{U_{i,j} - U_{i-1,j}}{h}\right)^2 + \left(\frac{U_{i,j} - U_{i,j-1}}{h}\right)^2 \right] + \sum_{i=1}^{\tilde{d}} \sum_{j=1}^{\tilde{d}} \frac{1}{4\lambda} (1 - U_{i,j}^2)^2,
    \label{eq:discrete-GL-2D}
\end{align}
with boundary conditions 
\begin{align}
    U_{0,:}=U_{\tilde{d}+1,:} = \mathbf{1}, \ U_{:,0}=U_{:,\tilde{d}+1} = -\mathbf{1}.
\end{align}
Here the total number of dimensionality is $d=\tilde{d}^2$. We set $\tilde{d}=4$, $h=1/(d+1)$, $\lambda=0.03$ and the inverse temperature $\beta=1/10$ for this example. All other settings remain the same as in the 1D case (\secref{sec:2D-GL}). We use the same set of $20$ basis functions for all dimensions as shown in \figref{fig:DW_basis} and we use the same cluster basis with $\tilde r=240$ as in the previous section. We solve the Fokker-Planck equation starting from the initial uniform distribution over the hypercube $[-M, M]^d$ and evolve the distribution with $\delta t = 0.002$ time and $N=10^4$ samples. To order the dimensions in a 2D lattice into a chain-like structure, we use the $d=16$ space-filling curve (\figref{fig:2D_space_filling}). 

\begin{figure}[htb]
    \centering
    \begin{subfigure}{0.40\textwidth}
        \centering
        \includegraphics[width=\textwidth]{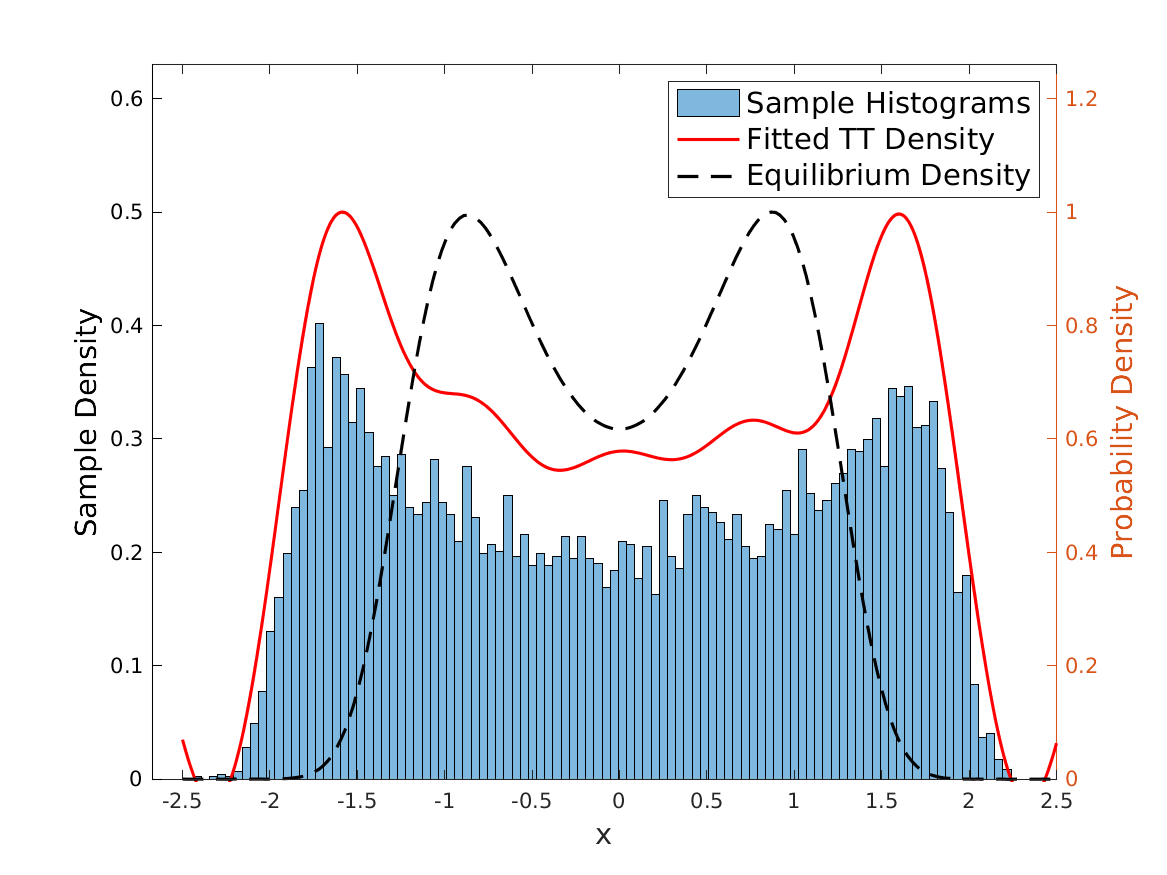}
        \caption{Iteration 1}
        \label{fig:GL_2D_step1}
    \end{subfigure}\quad \quad
    \begin{subfigure}{0.40\textwidth}
        \centering
        \includegraphics[width=\textwidth]{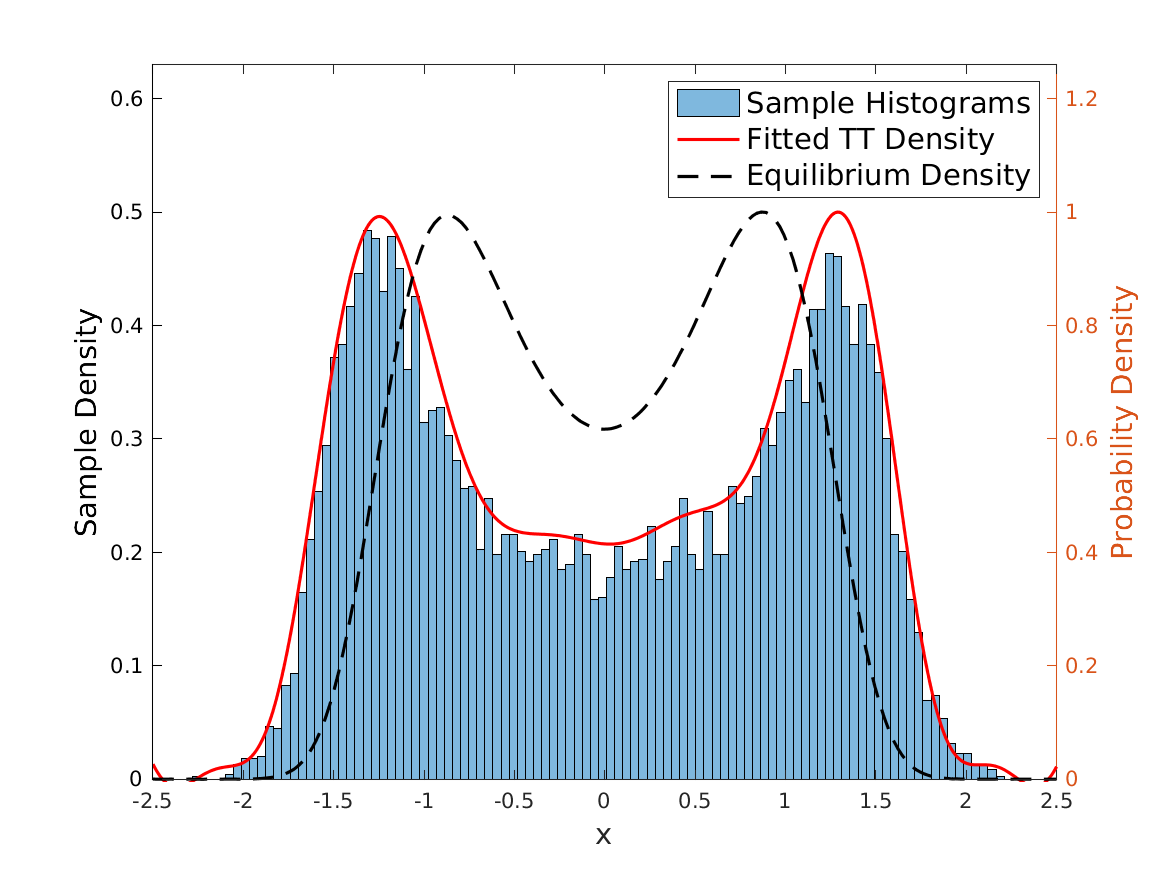}
        \caption{Iteration 3}
        \label{fig:GL_2D_step2}
    \end{subfigure}\\
    \begin{subfigure}{0.40\textwidth}
        \centering
        \includegraphics[width=\textwidth]{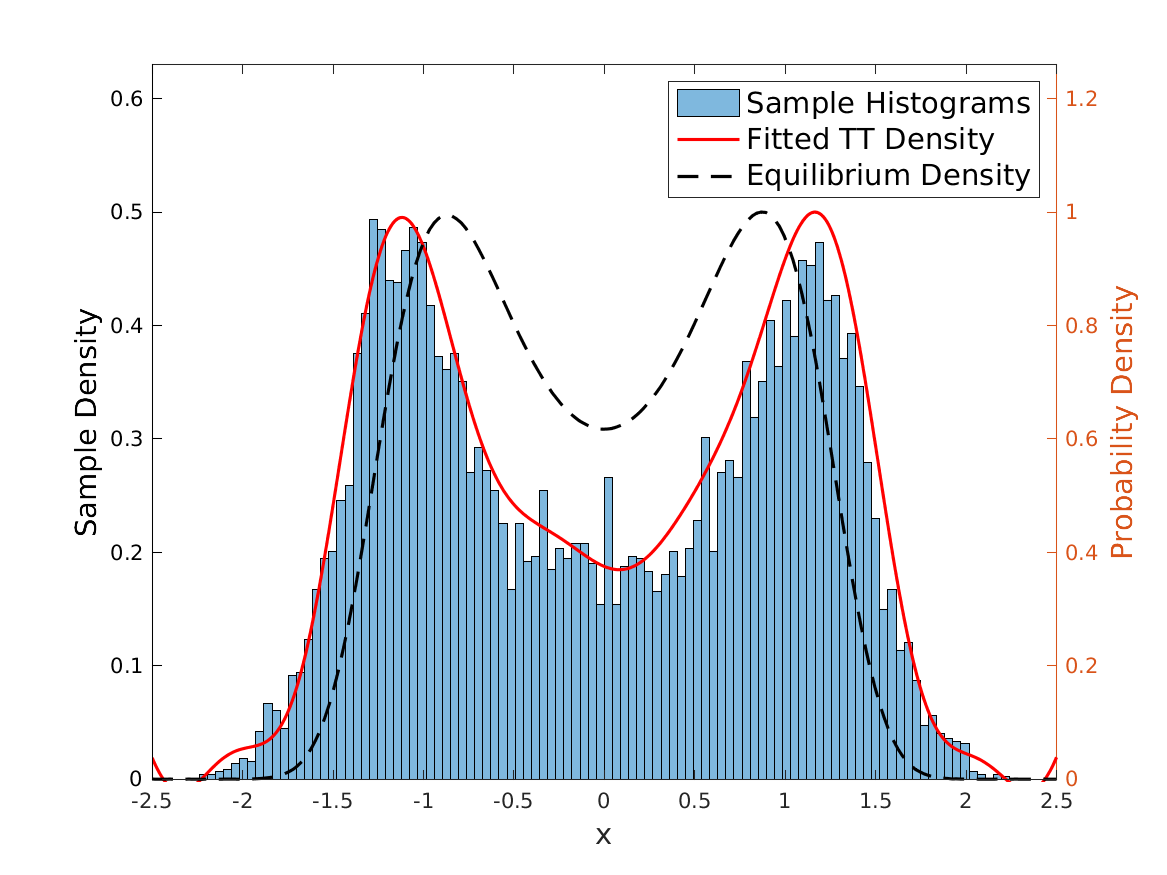}
        \caption{Iteration 5}
        \label{fig:GL_2D_step3}
    \end{subfigure}\quad \quad
    \begin{subfigure}{0.40\textwidth}
        \centering
        \includegraphics[width=\textwidth]{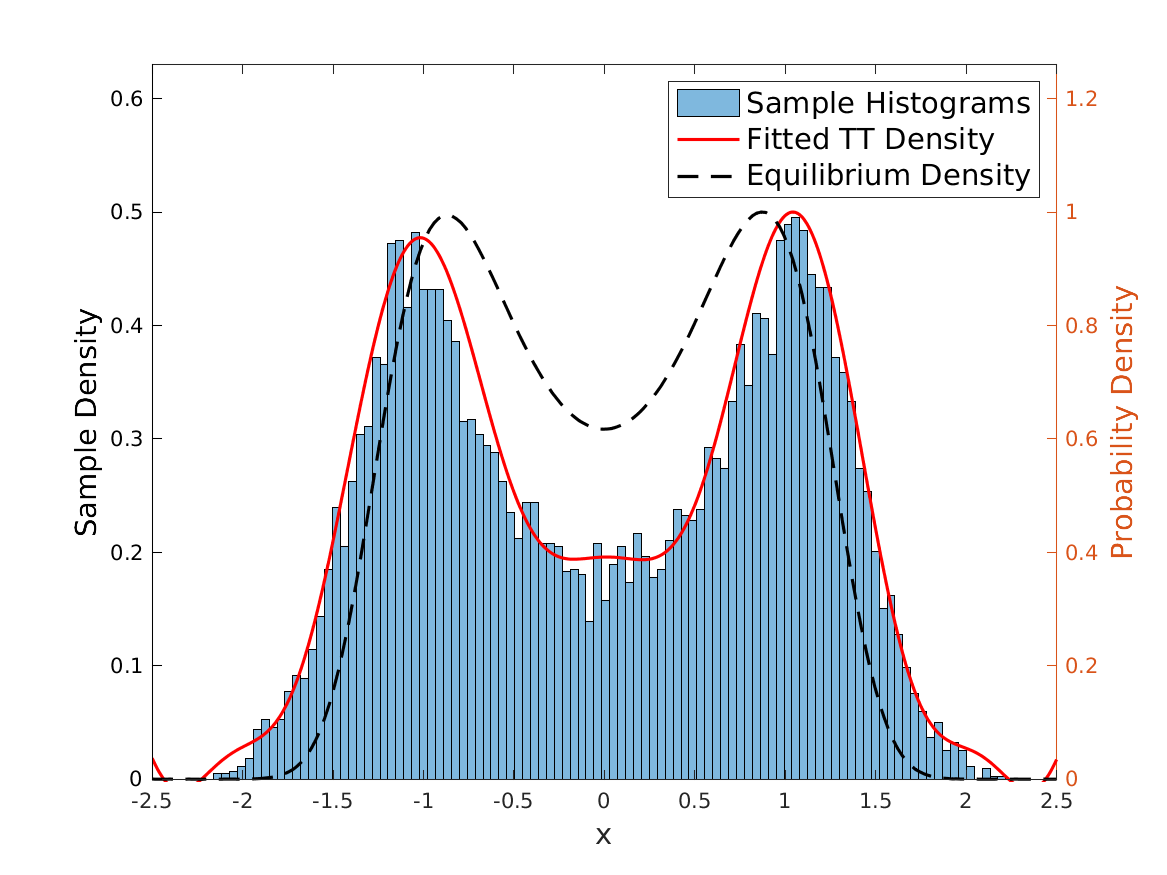}
        \caption{Iteration 7}
        \label{fig:GL_2D_step4}
    \end{subfigure}\\
    \begin{subfigure}{0.40\textwidth}
        \centering
        \includegraphics[width=\textwidth]{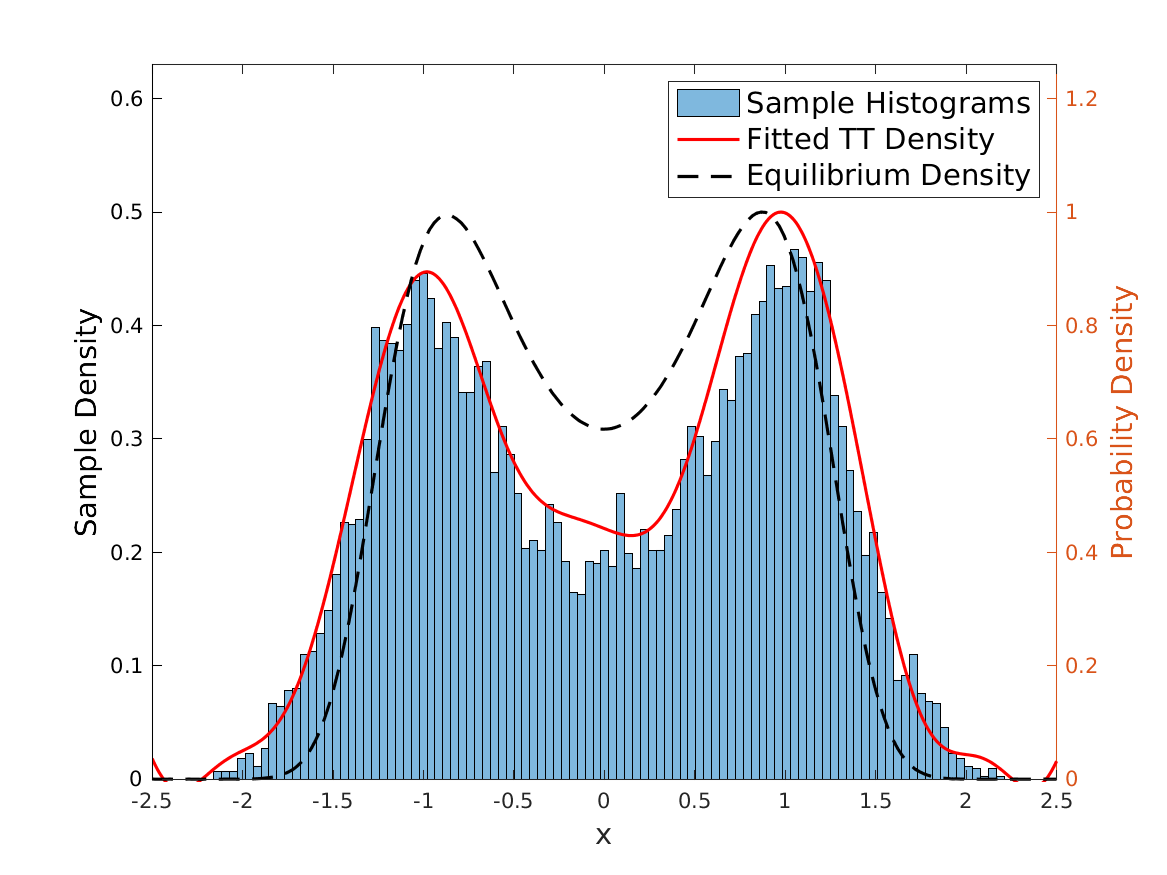}
        \caption{Iteration 20}
        \label{fig:GL_2D_step5}
    \end{subfigure}\quad \quad
    \begin{subfigure}{0.40\textwidth}
        \centering
        \includegraphics[width=\textwidth]{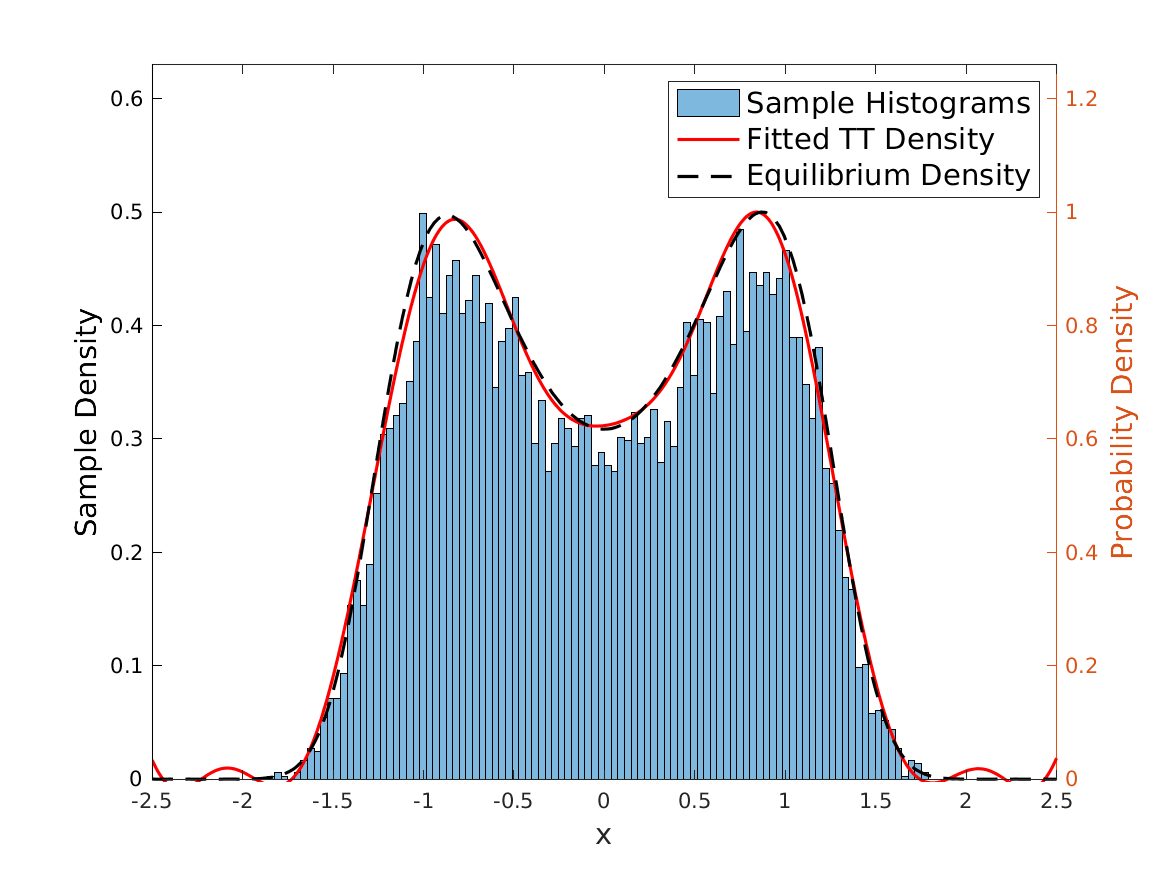}
        \caption{Iteration 100}
        \label{fig:GL_2D_step6}
    \end{subfigure}
    \caption{Visualization of the evolution of the 8-th marginal distribution for the 2D Ginzburg-Landau potential. The blue histograms correspond to the sample histograms after Langevin simulations at each iteration ($\hat \phi_{t+1}$ in \eqref{eq:ivp}). The estimated continuous MPS/TT density $\hat \phi_{\theta_{t+1}}$ in \eqref{eq:ivp} and the target equilibrium density $\phi^*$ are represented with red solid lines and black dashed lines, respectively. }
    \label{fig:GL_2D_evolution}
\end{figure}

In \figref{fig:GL_2D_evolution} we visualize the 8-th marginal distribution of the particle dynamics, the MPS/TT density, and the equilibrium density, at iteration $1,3,5,7,20$ and $100$. At iteration $100$, the relative error of the 8-th marginal distribution is $E=2.9\times 10^{-2}$. 

\section{Conclusion}
\label{sec:discussions}

In this paper, we propose a novel and general framework that combines Monte Carlo simulation with an MPS/TT ansatz. By leveraging the advantages of both approaches, our method offers an efficient way to apply the semigroup/time-evolution operator and control the variance and population of random walkers using tensor-sketching techniques.

The performance of our algorithm is determined by two factors: the number of randomized sketches and the number of samples used. Our algorithm is expected to succeed when we can employ samples to determine a low-rank MPS/TT representation based on estimating certain low-order moments.  Hence, it is crucial to investigate the function space under time evolution and understand how the MPS/TT representation of the solution can be efficiently determined using Monte Carlo method in a statistically optimal manner.

\section{Acknowledgements}
Y.C. and Y.K. acknowledge partial support from NSF Award No.\ DMS-211563 and DOE Award No.\ DE-SC0022232. The authors also thank  Michael Lindsey and Shiwei Zhang regarding discussions on potential improvements for the proposed method. 

\newpage
\bibliographystyle{plain}
\bibliography{ref}

@article{qin2016coupling,
  title={Coupling quantum Monte Carlo and independent-particle calculations: Self-consistent constraint for the sign problem based on the density or the density matrix},
  author={Qin, Mingpu and Shi, Hao and Zhang, Shiwei},
  journal={Physical Review B},
  volume={94},
  number={23},
  pages={235119},
  year={2016},
  publisher={APS}
}

@article{chen2000equivalence,
  title={Equivalence of exponential ergodicity and L2-exponential convergence for Markov chains},
  author={Chen, Mu-Fa},
  journal={Stochastic processes and their applications},
  volume={87},
  number={2},
  pages={281--297},
  year={2000},
  publisher={Elsevier}
}

@book{wendland2017numerical,
  title={Numerical linear algebra: An introduction},
  author={Wendland, Holger},
  volume={56},
  year={2017},
  publisher={Cambridge University Press}
}

@book{bhattacharya2009stochastic,
  title={Stochastic processes with applications},
  author={Bhattacharya, Rabi N and Waymire, Edward C},
  year={2009},
  publisher={SIAM}
}

@article{suzuki1976relationship,
  title={Relationship between d-dimensional quantal spin systems and (d+ 1)-dimensional ising systems: Equivalence, critical exponents and systematic approximants of the partition function and spin correlations},
  author={Suzuki, Masuo},
  journal={Progress of theoretical physics},
  volume={56},
  number={5},
  pages={1454--1469},
  year={1976},
  publisher={Oxford University Press}
}

@article{trotter1959product,
  title={On the product of semi-groups of operators},
  author={Trotter, Hale F},
  journal={Proceedings of the American Mathematical Society},
  volume={10},
  number={4},
  pages={545--551},
  year={1959}
}

@article{ulmke2000auxiliary,
  title={Auxiliary-field Monte Carlo for quantum spin and boson systems},
  author={Ulmke, Martin and Scalettar, RT},
  journal={Physical Review B},
  volume={61},
  number={14},
  pages={9607},
  year={2000},
  publisher={APS}
}

@Article{Ceperley1995,
  Title                    = {Path integrals in the theory of condensed helium},
  Author                   = {D. M. Ceperley},
  Journal                  = {Rev. Mod. Phys.},
  Year                     = {1995},
  Pages                    = {279},
  Volume                   = {67}
}

@article{tt-decomposition,
  title={Tensor-train decomposition},
  author={Oseledets, Ivan V},
  journal={SIAM Journal on Scientific Computing},
  volume={33},
  number={5},
  pages={2295--2317},
  year={2011},
  publisher={SIAM}
}

@book{GL-model,
  title={Ginzburg-Landau phase transition theory and superconductivity},
  author={Hoffmann, K-H and Tang, Qi},
  volume={134},
  year={2012},
  publisher={Birkh{\"a}user}
}

@article{dolgov2011fast,
  title={Fast solution of multi-dimensional parabolic problems in the TT/QTT-format with initial application to the Fokker-Planck equation},
  author={Dolgov, Sergey and Khoromskij, Boris N and Oseledets, Ivan V},
  year={2011}
}

@article{dolgov2020approximation,
  title={Approximation and sampling of multivariate probability distributions in the tensor train decomposition},
  author={Dolgov, Sergey and Anaya-Izquierdo, Karim and Fox, Colin and Scheichl, Robert},
  journal={Statistics and Computing},
  volume={30},
  pages={603--625},
  year={2020},
  publisher={Springer}
}

@article{tang2022generative,
  title={Generative modeling via tree tensor network states},
  author={Tang, Xun and Hur, Yoonhaeng and Khoo, Yuehaw and Ying, Lexing},
  journal={arXiv preprint arXiv:2209.01341},
  year={2022}
}

@article{hur2022generative,
  title={Generative modeling via tensor train sketching},
  author={Hur, Yoonhaeng and Hoskins, Jeremy G and Lindsey, Michael and Stoudenmire, E Miles and Khoo, Yuehaw},
  journal={arXiv preprint arXiv:2202.11788},
  year={2022}
}

@article{ruthotto2021introduction,
  title={An introduction to deep generative modeling},
  author={Ruthotto, Lars and Haber, Eldad},
  journal={GAMM-Mitteilungen},
  volume={44},
  number={2},
  pages={e202100008},
  year={2021},
  publisher={Wiley Online Library}
}

@article{tt_committor,
  title={Committor functions via tensor networks},
  author={Chen, Yian and Hoskins, Jeremy and Khoo, Yuehaw and Lindsey, Michael},
  journal={Journal of Computational Physics},
  volume={472},
  pages={111646},
  year={2023},
  publisher={Elsevier}
}

@article{han2018solving,
  title={Solving high-dimensional partial differential equations using deep learning},
  author={Han, Jiequn and Jentzen, Arnulf and E, Weinan},
  journal={Proceedings of the National Academy of Sciences},
  volume={115},
  number={34},
  pages={8505--8510},
  year={2018},
  publisher={National Acad Sciences}
}

@article{yu2018deep,
  title={The deep Ritz method: a deep learning-based numerical algorithm for solving variational problems},
  author={Yu, Bing and others},
  journal={Communications in Mathematics and Statistics},
  volume={6},
  number={1},
  pages={1--12},
  year={2018},
  publisher={Springer}
}

@inproceedings{zhai2022deep,
  title={A deep learning method for solving Fokker-Planck equations},
  author={Zhai, Jiayu and Dobson, Matthew and Li, Yao},
  booktitle={Mathematical and Scientific Machine Learning},
  pages={568--597},
  year={2022},
  organization={PMLR}
}

@article{ambartsumyan2020hierarchical,
  title={Hierarchical matrix approximations of Hessians arising in inverse problems governed by PDEs},
  author={Ambartsumyan, Ilona and Boukaram, Wajih and Bui-Thanh, Tan and Ghattas, Omar and Keyes, David and Stadler, Georg and Turkiyyah, George and Zampini, Stefano},
  journal={SIAM Journal on Scientific Computing},
  volume={42},
  number={5},
  pages={A3397--A3426},
  year={2020},
  publisher={SIAM}
}

@article{chen2021scalable,
  title={Scalable Gaussian Process Analysis for Implicit Physics-Based Covariance Models},
  author={Chen, Yian and Anitescu, Mihai},
  journal={International Journal for Uncertainty Quantification},
  volume={11},
  number={6},
  year={2021},
  publisher={Begel House Inc.}
}

@article{chen2023scalable,
  title={Scalable Physics-based Maximum Likelihood Estimation using Hierarchical Matrices},
  author={Chen, Yian and Anitescu, Mihai},
  journal={arXiv preprint arXiv:2303.10102},
  year={2023}
}

@article{ho2016hierarchical,
  title={Hierarchical interpolative factorization for elliptic operators: differential equations},
  author={Ho, Kenneth L and Ying, Lexing},
  journal={Communications on Pure and Applied Mathematics},
  volume={69},
  number={8},
  pages={1415--1451},
  year={2016},
  publisher={Wiley Online Library}
}

@article{zhang2013auxiliary,
  title={Auxiliary-Field Quantum Monte Carlo for Correlated Electron Systems. Emergent Phenomena in Correlated Matter: Autumn School Organized by the Forschungszentrum J{\"u}lich and the German Research School for Simulation Sciences at Forschungszentrum J{\"u}lich 23-27 September 2013},
  author={Zhang, Shiwei},
  journal={Lecture Notes of the Autumn School Correlated Electrons},
  volume={3},
  pages={2013},
  year={2013}
}

@article{shi2021some,
  title={Some recent developments in auxiliary-field quantum Monte Carlo for real materials},
  author={Shi, Hao and Zhang, Shiwei},
  journal={The Journal of Chemical Physics},
  volume={154},
  number={2},
  pages={024107},
  year={2021},
  publisher={AIP Publishing LLC}
}

@article{afmc2,
  title={Computation within the auxiliary field approach},
  author={Baeurle, SA},
  journal={Journal of Computational Physics},
  volume={184},
  number={2},
  pages={540--558},
  year={2003},
  publisher={Elsevier}
}

@article{afqmc1,
  title={Monte Carlo calculations of coupled boson-fermion systems. I},
  author={Blankenbecler, Richard and Scalapino, DJ and Sugar, RL},
  journal={Physical Review D},
  volume={24},
  number={8},
  pages={2278},
  year={1981},
  publisher={APS}
}

@article{gull1993imaginary,
  title={Imaginary numbers are not real—the geometric algebra of spacetime},
  author={Gull, Stephen and Lasenby, Anthony and Doran, Chris},
  journal={Foundations of Physics},
  volume={23},
  number={9},
  pages={1175--1201},
  year={1993},
  publisher={Springer}
}

@article{hubbard1959calculation,
  title={Calculation of partition functions},
  author={Hubbard, John},
  journal={Physical Review Letters},
  volume={3},
  number={2},
  pages={77},
  year={1959},
  publisher={APS}
}

@article{carlson2011auxiliary,
  title={Auxiliary-field quantum Monte Carlo method for strongly paired fermions},
  author={Carlson, J and Gandolfi, Stefano and Schmidt, Kevin E and Zhang, Shiwei},
  journal={Physical Review A},
  volume={84},
  number={6},
  pages={061602},
  year={2011},
  publisher={APS}
}

@article{qin2016benchmark,
  title={Benchmark study of the two-dimensional Hubbard model with auxiliary-field quantum Monte Carlo method},
  author={Qin, Mingpu and Shi, Hao and Zhang, Shiwei},
  journal={Physical Review B},
  volume={94},
  number={8},
  pages={085103},
  year={2016},
  publisher={APS}
}

@article{lee2022twenty,
  title={Twenty Years of Auxiliary-Field Quantum Monte Carlo in Quantum Chemistry: An Overview and Assessment on Main Group Chemistry and Bond-Breaking},
  author={Lee, Joonho and Pham, Hung Q and Reichman, David R},
  journal={Journal of Chemical Theory and Computation},
  volume={18},
  number={12},
  pages={7024--7042},
  year={2022},
  publisher={ACS Publications}
}

@article{ren2022high,
  title={High-dimensional density estimation with tensorizing flow},
  author={Ren, Yinuo and Zhao, Hongli and Khoo, Yuehaw and Ying, Lexing},
  journal={arXiv preprint arXiv:2212.00759},
  year={2022}
}

@book{sagan2012space,
  title={Space-filling curves},
  author={Sagan, Hans},
  year={2012},
  publisher={Springer Science \& Business Media}
}

@article{sandvik2007variational,
  title={Variational quantum Monte Carlo simulations with tensor-network states},
  author={Sandvik, Anders W and Vidal, Guifre},
  journal={Physical review letters},
  volume={99},
  number={22},
  pages={220602},
  year={2007},
  publisher={APS}
}

@article{chertkov2021solution,
  title={Solution of the Fokker--Planck equation by cross approximation method in the tensor train format},
  author={Chertkov, Andrei and Oseledets, Ivan},
  journal={Frontiers in Artificial Intelligence},
  volume={4},
  pages={668215},
  year={2021},
  publisher={Frontiers Media SA}
}

@article{oseledets2011tensor,
  title={Tensor-train decomposition},
  author={Oseledets, Ivan V},
  journal={SIAM Journal on Scientific Computing},
  volume={33},
  number={5},
  pages={2295--2317},
  year={2011},
  publisher={SIAM}
}

@book{howie1995fundamentals,
  title={Fundamentals of semigroup theory},
  author={Howie, John Mackintosh},
  number={12},
  year={1995},
  publisher={Oxford University Press}
}

@article{clifford1961algebraic,
  title={The algebraic theory of semigroups, vol. 1},
  author={Clifford, Alfred H and Preston, Gordon B},
  journal={AMS surveys},
  volume={7},
  pages={1967},
  year={1961}
}

@article{hollings2009early,
  title={The early development of the algebraic theory of semigroups},
  author={Hollings, Christopher},
  journal={Archive for history of exact sciences},
  volume={63},
  pages={497--536},
  year={2009},
  publisher={Springer}
}

@article{wang2015fast,
  title={Fast and guaranteed tensor decomposition via sketching},
  author={Wang, Yining and Tung, Hsiao-Yu and Smola, Alexander J and Anandkumar, Anima},
  journal={Advances in neural information processing systems},
  volume={28},
  year={2015}
}

@article{ahle2019almost,
  title={Almost optimal tensor sketch},
  author={Ahle, Thomas D and Knudsen, Jakob BT},
  journal={arXiv preprint arXiv:1909.01821},
  year={2019}
}

@article{peng2023generative,
  title={Generative Modeling via Hierarchical Tensor Sketching},
  author={Peng, Yifan and Chen, Yian and Stoudenmire, E Miles and Khoo, Yuehaw},
  journal={arXiv preprint arXiv:2304.05305},
  year={2023}
}

@article{chan2011density,
  title={The density matrix renormalization group in quantum chemistry},
  author={Chan, Garnet Kin-Lic and Sharma, Sandeep},
  journal={Annual review of physical chemistry},
  volume={62},
  pages={465--481},
  year={2011},
  publisher={Annual Reviews}
}

@article{white1992density,
  title={Density matrix formulation for quantum renormalization groups},
  author={White, Steven R},
  journal={Physical review letters},
  volume={69},
  number={19},
  pages={2863},
  year={1992},
  publisher={APS}
}

@article{cao2018stochastic,
  title={Stochastic dynamical low-rank approximation method},
  author={Cao, Yu and Lu, Jianfeng},
  journal={Journal of Computational Physics},
  volume={372},
  pages={564--586},
  year={2018},
  publisher={Elsevier}
}

@inproceedings{lee2022unbiasing,
  title={Unbiasing Fermionic Quantum Monte Carlo with a Quantum Computer},
  author={Lee, Joonho},
  booktitle={APS March Meeting Abstracts},
  volume={2022},
  pages={N40--007},
  year={2022}
}

@article{lim2017fast,
  title={Fast randomized iteration: Diffusion Monte Carlo through the lens of numerical linear algebra},
  author={Lim, Lek-Heng and Weare, Jonathan},
  journal={SIAM Review},
  volume={59},
  number={3},
  pages={547--587},
  year={2017},
  publisher={SIAM}
}

@article{greene2019beyond,
  title={Beyond walkers in stochastic quantum chemistry: Reducing error using fast randomized iteration},
  author={Greene, Samuel M and Webber, Robert J and Weare, Jonathan and Berkelbach, Timothy C},
  journal={Journal of chemical Theory and Computation},
  volume={15},
  number={9},
  pages={4834--4850},
  year={2019},
  publisher={ACS Publications}
}

@article{burkhardt1985finite,
  title={Finite-size scaling of the quantum Ising chain with periodic, free, and antiperiodic boundary conditions},
  author={Burkhardt, Th W and Guim, IHNSOUK},
  journal={Journal of Physics A: Mathematical and General},
  volume={18},
  number={1},
  pages={L33},
  year={1985},
  publisher={IOP Publishing}
}

@article{lieb1961two,
  title={Two soluble models of an antiferromagnetic chain},
  author={Lieb, Elliott and Schultz, Theodore and Mattis, Daniel},
  journal={Annals of Physics},
  volume={16},
  number={3},
  pages={407--466},
  year={1961},
  publisher={Elsevier}
}

\end{document}